%% file: wpca,arxiv.tex
\documentclass{article}

\usepackage{fullpage}
\usepackage[dvipsnames]{xcolor}
\usepackage{ifpdf}

\usepackage{caption}
\usepackage{hyperref}
\colorlet{inlinkcolor}{green!50!black}
\colorlet{exlinkcolor}{red!50!black}
\hypersetup{
  colorlinks=true,
  allcolors = inlinkcolor,
  urlcolor = exlinkcolor,
  frenchlinks=false,
  pdfborder={0 0 0},
  naturalnames=false,
  hypertexnames=false,
  breaklinks
}

\usepackage{amsmath,amsthm}
\usepackage[capitalize,nameinlink]{cleveref}
\crefname{section}{section}{sections}
\crefname{subsection}{subsection}{subsections}
\Crefname{section}{Section}{Sections}
\Crefname{subsection}{Subsection}{Subsections}
\Crefname{figure}{Figure}{Figures}
\crefformat{equation}{\textup{#2(#1)#3}}
\crefrangeformat{equation}{\textup{#3(#1)#4--#5(#2)#6}}
\crefmultiformat{equation}{\textup{#2(#1)#3}}{ and \textup{#2(#1)#3}}
{, \textup{#2(#1)#3}}{, and \textup{#2(#1)#3}}
\crefrangemultiformat{equation}{\textup{#3(#1)#4--#5(#2)#6}}%
{ and \textup{#3(#1)#4--#5(#2)#6}}{, \textup{#3(#1)#4--#5(#2)#6}}{, and \textup{#3(#1)#4--#5(#2)#6}}
\Crefformat{equation}{#2Equation~\textup{(#1)}#3}
\Crefrangeformat{equation}{Equations~\textup{#3(#1)#4--#5(#2)#6}}
\Crefmultiformat{equation}{Equations~\textup{#2(#1)#3}}{ and \textup{#2(#1)#3}}
{, \textup{#2(#1)#3}}{, and \textup{#2(#1)#3}}
\Crefrangemultiformat{equation}{Equations~\textup{#3(#1)#4--#5(#2)#6}}%
{ and \textup{#3(#1)#4--#5(#2)#6}}{, \textup{#3(#1)#4--#5(#2)#6}}{, and \textup{#3(#1)#4--#5(#2)#6}}
\crefdefaultlabelformat{#2\textup{#1}#3}

\theoremstyle{plain}
\newtheorem{theorem}{Theorem}[section]
\newtheorem{lemma}[theorem]{Lemma}
\newtheorem{corollary}[theorem]{Corollary}
\newtheorem{proposition}[theorem]{Proposition}
\theoremstyle{definition}
\newtheorem{example}[theorem]{Example}
\newtheorem{remark}[theorem]{Remark}

\numberwithin{equation}{section}

\newcommand{\insupp}[1]{\cref{#1}} % References to the appendix
\newcommand{\Insupp}[1]{\Cref{#1}}
\newcommand{\supporapp}{appendix}

\newcommand{\email}[1]{\protect\href{mailto:#1}{#1}}
\newenvironment{keywords}{\begin{small}\noindent\textbf{Key words.}}{\end{small}}
\newenvironment{AMS}{\begin{small}\noindent\textbf{AMS subject classifications.}}{\end{small}}

\usepackage{amsfonts}
\usepackage{graphicx}
\usepackage{algorithmic}
\ifpdf
  \DeclareGraphicsExtensions{.eps,.pdf,.png,.jpg}
\else
  \DeclareGraphicsExtensions{.eps}
\fi
\usepackage{xspace}

\usepackage{amsmath,amssymb,amsopn}
\usepackage{dsfont,mathtools}
\usepackage{altvars}

\usepackage{pgfplots}
\usetikzlibrary{math,cd}
\usepgfplotslibrary{colorbrewer,fillbetween}
\pgfplotsset{compat=1.12}
\usepackage[caption=false]{subfig}

\ifpdf
  \usepackage{crossreftools}
\fi

\usepackage{enumitem}
\newlist{propenum}{enumerate}{1}
\setlist[propenum]{label=(\alph*)}
\crefname{propenumi}{Property}{Properties}

\title{Optimally Weighted PCA for\\High-Dimensional Heteroscedastic Data%
  \thanks{%
    D. Hong was supported in part by
    NSF Graduate Research Fellowship DGE 1256260,
    NSF Grant ECCS-1508943,
    NSF BIGDATA grant IIS 1837992,
    the Dean's Fund for Postdoctoral Research of the Wharton School,
    and NSF Mathematical Sciences Postdoctoral Research Fellowship DMS 2103353.
    F. Yang was supported in part by the Wharton Dean's Fund for Postdoctoral Research.
    J. A. Fessler was supported in part by
    the UM-SJTU data science seed fund,
    NSF grant IIS 1838179,
    and NIH Grant U01 EB 018753.
    L. Balzano was supported in part by
    DARPA-16-43-D3M-FP-037, by NSF CAREER award CCF-1845076,
    by NSF Grant ECCS-1508943
    and by ARO YIP Award W911NF1910027.
  }%
}

\author{%
  David Hong%
  \thanks{%
    Department of Statistics and Data Science,
    Wharton School,
    University of Pennsylvania,
    Philadelphia, PA, 19104 USA
    (\email{dahong67@wharton.upenn.edu}).
  }
  \and
  Fan Yang%
  \thanks{%
    Yau Mathematical Sciences Center, Tsinghua University,
    Beijing, 100084 China
    (\email{fyangmath@mail.tsinghua.edu.cn}).
  }
  \and
  Jeffrey A. Fessler%
  \thanks{%
    Department of Electrical Engineering and Computer Science,
    University of Michigan,
    Ann Arbor, MI, 48109 USA
    (\email{fessler@umich.edu}, \email{girasole@umich.edu}).
  }
  \and
  Laura Balzano%
  \footnotemark[4]
}

\DeclareMathOperator*{\diag}{diag}
\DeclareMathOperator*{\blockdiag}{blockdiag}
\DeclareMathOperator*{\tr}{tr}
\DeclareMathOperator{\proj}{proj}
\newcommand{\CT}{\mathsf{H}}
\newcommand{\frob}{\mathrm{F}}
\newcommand{\opnorm}{\mathrm{op}}

\DeclareMathOperator*{\aslim}{aslim}
\newcommand{\asto}{\overset{\mathrm{a.s.}}{\longrightarrow}}
\newcommand{\unifto}[1]{\underset{#1}{\rightrightarrows}}
\newcommand{\asunifto}[1]{\overset{\mathrm{a.s.}}{\unifto{#1}}}
\newcommand{\BigO}{\mathrm{O}}

\DeclareMathOperator{\supp}{supp}
\DeclareMathOperator*{\argmax}{argmax}
\newcommand{\diff}[2]{\frac{\partial #1}{\partial #2}}
\newcommand{\idiff}[2]{\partial #1/\partial #2}
\newcommand{\iid}{\mathrm{IID}}
\newcommand{\ith}{$i$th\xspace}

\newcommand{\Sedge}{S_{\mathrm{edge}}}
\newcommand{\Sout}{S_{\mathrm{out}}}
\newcommand{\mplaw}{m_{\mathrm{MP}}}

\newcommand{\opt}{\star}
\newcommand{\inv}{\mathrm{(inv)}}
\newcommand{\full}{\mathrm{(full)}}

\renewcommand{\Re}{\operatorname{Re}}
\renewcommand{\Im}{\operatorname{Im}}

\newcommand{\revone}[1]{#1}

\begin{document}

\maketitle

\begin{abstract}
  \input{content/abstract}
\end{abstract}

\begin{keywords}
principal component analysis,
large-dimensional data,
heterogeneous quality,
optimal weighting
\end{keywords}

\begin{AMS}
62H25
\end{AMS}

\input{content/intro}
\input{content/main,result}
\input{content/simulation}
\input{content/comparison}
\input{content/ext,comp,more,methods}
\input{content/derivation}
\input{content/ext,signal,var,est}

\input{content/experiment}
\input{content/conclusion}

\input{acknowledgement}
\bibliographystyle{siamplain}
\bibliography{refs}

\clearpage
\appendix

\input{content/ext,het,signal}
\input{content/supp,derivation,calc}
\input{content/supp,ext,signal,var}
\input{content/experiment,preprocess}

\end{document}

%% file: content/abstract.tex
%!TEX root = ../wpca.tex
% 
Modern data are increasingly both high-dimensional
and heteroscedastic.
This paper considers the challenge
of estimating underlying principal components
from high-dimensional data
with noise that is heteroscedastic across samples,
i.e., some samples are noisier than others.
Such heteroscedasticity naturally arises,
e.g., when combining data from diverse sources or sensors.
A natural way to account for this heteroscedasticity
is to give noisier blocks of samples less weight in PCA
by using the leading eigenvectors of a weighted sample covariance matrix.
We consider the problem of choosing weights
to optimally recover the underlying components.
In general, one cannot know these optimal weights
since they depend on the underlying components we seek to estimate.
However, we show that under some natural statistical assumptions
the optimal weights converge
to a simple function of the signal and noise variances
for high-dimensional data.
Surprisingly, the optimal weights are not the inverse noise variance weights
commonly used in practice.
We demonstrate the theoretical results through numerical simulations
and comparisons with existing weighting schemes.
Finally, we briefly discuss how estimated signal and noise variances
can be used when the true variances are unknown,
and we illustrate the optimal weights on real data from astronomy.

%% file: content/intro.tex
%!TEX root = ../wpca.tex

\section{Introduction} \label{sct:intro}%
Principal Component Analysis (PCA) is a fundamental technique
for discovering underlying components in data
and is a workhorse method
for analyzing modern \emph{high-dimensional} data.
However, conventional PCA
does not recover underlying principal components well
when the data has \emph{heteroscedastic} noise,
as is common in practice.
In particular,
its performance can degrade substantially
when the noise is heteroscedastic across samples,
i.e., some samples are noisier than others.
PCA suffers from treating all the samples uniformly
with performance held back by the noisiest samples,
as was rigorously characterized in \cite{hong2018apo}.
Weighted PCA addresses this shortcoming
by giving less weight to lower quality samples.
This naturally raises
a crucial question:
how should the weights be chosen?
Namely, \emph{what are the optimal weights?}

This paper addresses this question
by rigorously deriving optimal weights
that are simple functions
of the signal and noise variances.
Surprisingly, they are not the inverse noise variance weights
that are commonly used in practice.
We now elaborate in more detail.

\subsection{High-dimensional and heteroscedastic data}%
Modern applications of PCA span numerous and diverse areas
across all of engineering and the sciences,
ranging from medical imaging~\cite{ardekani1999adi,pedersen2009ktp}
to cancer data classification~\cite{sharma2015and},
genetics~\cite{leek2010acs},
and environmental sensing~\cite{papadimitriou2005spd,wagner1996sdu},
to name just a few.
Increasingly, the number of features measured is
comparable to or even larger than the number of samples,
i.e., the data are \emph{high-dimensional}.
Traditional asymptotic analysis
of the performance of methods
as the number of samples grows
(with a fixed number of features)
do not apply well to such settings.
Modern data analysis needs new theory and methods
for the high-dimensional regime
where both the number of features
and number of samples are large \cite{johnstone2009sco}.

Modern datasets are also frequently composed of samples
with heteroscedastic
(i.e., heterogeneous)
noise.
In particular,
we consider noise that is \emph{heteroscedastic across samples},
namely,
some samples are noisier than others.
Such data arises naturally
when samples are obtained
at varying times or by varying means or equipment.
For example,
in the field of analytical chemistry,
\cite{cochran1977swp} considers spectrophotometric data
obtained from averages taken over
varying windows of time;
samples from shorter windows
are noisier.
As another example,
in the field of air quality monitoring,
samples come from various sources:
governments agencies provide low-noise data
obtained from carefully operated instruments,
while
individuals provide noisier data
obtained from cheaper and easy to setup sensors
\cite{epa2021aqs,purpleair2021rta}.
As a final example,
in the field of astronomy,
measurements of astronomical objects
such as stars and quasars
can have various levels of noise
due to atmospheric and detector effects
that vary from object to object
\cite{bailey2012pca,tamuz2005cse,tsalmantza2012add}.
More generally,
modern big data analysis is often performed using
datasets built up by combining myriad sources,
so one can expect that data with heteroscedastic noise
will be the norm.
Modern data analysis needs PCA methods
that effectively account for this type of heteroscedasticity.
Indeed, such methods may also unlock new opportunities
to effectively leverage new sources of data
with heteroscedastic noise.

\subsection{Weighted PCA}
\label{sct:wpca}%
Weighted PCA accounts for heteroscedastic noise
by giving smaller weight to noisier samples.
Analogous to unweighted PCA,
the principal components
are the leading eigenvectors of
the \emph{weighted sample covariance matrix}
\begin{equation*}
  \bhtSigma_{\bmw}
  \coloneqq
  \sum_{\ell=1}^L
  w_\ell \bmY_\ell \bmY_\ell^\CT
  ,
\end{equation*}
where $\bmY_1,\dots,\bmY_L$ are
$L$ blocks of samples
with associated noise variances $v_1,\dots,v_L > 0$,
the superscript $\CT$ denotes the Hermitian transpose,
and $w_1,\dots,w_L \geq 0$ are the weights.
Existing choices
for the weights include:
\begin{itemize}

\item \textbf{uniform weights}
($w_\ell = 1$):
these weights correspond to unweighted PCA
and can be a natural choice
when the noise is close to homoscedastic.
However, its performance degrades
with increasing noise heteroscedasticity,
as was shown in \cite[Theorem 2]{hong2018apo}.

\item \textbf{binary weights}
($w_\ell = 1$ for less noisy blocks
and $w_\ell = 0$ for the rest):
these weights correspond to performing unweighted PCA
using only less noisy blocks of samples
and are a natural choice
when some samples are much noisier than the rest.
The idea is to exclude noisier samples
that do more harm than good.
However, doing so also omits any useful information
that was in the excluded samples.
How to decide if a block of samples
is better to include or exclude
can be unclear.

\item \textbf{inverse noise variance weights}
($w_\ell = 1/v_\ell$):
these weights whiten the noise,
making it homoscedastic,
and can be interpreted as
a maximum likelihood weighting \cite{young1941mle}.
They are a natural way to account for noise heteroscedasticity
while using all the samples
and are commonly used in practice.

\end{itemize}
It has been unclear
which of these existing options to choose
and whether any are optimal.

\subsection{Contribution of this paper}%
The main contribution of this paper
is \emph{optimal weights}
for high-dimensional heteroscedastic data,
that are rigorously derived
under some natural statistical assumptions.
Roughly put,
we show that
when both dimensions of the data are large,
the optimal weights converge to
the following simple \emph{asymptotic optimal weights}:
\begin{equation*}
  w_\ell
  =
  \frac{1}{v_\ell}\frac{1}{1 + v_\ell/\lambda}
  ,
\end{equation*}
where $v_\ell$ is the noise variance and $\lambda$ is the signal variance.
Notably,
these weights are inverse noise variance weights
scaled by a simple term that depends on the noise-to-signal ratio $v_\ell/\lambda$.
See \cref{sec:opt:weight}
for the precise statement of the result
and \cref{proof:main:result}
for its proof.

Naturally, one wonders how well these results
apply for data with finitely many samples and features.
Numerical simulations in \cref{sec:opt:weight:sim}
illustrate that
the optimal weights in finite dimensions
(which are a function of the random signal coefficients and noise)
concentrate around the asymptotic optimal weights
as the data grows in size.
As a result,
the asymptotic optimal weights
are often close to the optimal weights in finite dimensions
when the dimensions are large enough.

We also compare the asymptotic optimal weights
with the existing weights above:
uniform,
binary,
and inverse noise variance weights.
In particular, we consider
how close they are to
the optimal weights in finite dimensions
(\cref{sec:weight:comp:weights}),
how well they perform in finite dimensions
(\cref{sec:weight:comp:rec}),
and in what regimes they achieve positive asymptotic recovery
(\cref{sec:weight:comp:phase}).
Overall, the asymptotic optimal weights
outperform the existing weights.

One also wonders how to calculate the asymptotic optimal weights
when the signal and noise variances are unknown.
Naturally, one might consider simply using estimators
of these variances.
We explain that the resulting estimated weights
are also asymptotically optimal
as long as the estimators are consistent,
and we give an example of such estimators
(\cref{sec:signal:var:est}).

Finally, we illustrate the asymptotic optimal weights on real data
(\cref{sec:exp:astro}).
The data are quasar spectra
measured by the Sloan Digital Sky Survey
and have heteroscedastic noise.
The example exhibits some of the main themes of the paper
and illustrates the potential for optimally weighted PCA
to improve performance in real data.

\subsection{Related works}%
\label{sec:related:works}%
Previous work on PCA for noise that is heteroscedastic
(whether across samples or otherwise)
have addressed various important questions,
as elaborated below.
However,
to the best of our knowledge,
the important question of optimal weighting
was not previously considered.
This paper rigorously answers the question of optimal weighting
for noise that is heteroscedastic across samples.
It will be interesting for future works to consider
this question for other forms of heteroscedastic noise.

Some of our other papers considered
various aspects of
noise that is heteroscedastic across samples.
In particular,
\cite{hong2016tat,hong2018apo}
derive the asymptotic performance of unweighted PCA
and characterize the impact of heteroscedasticity.
See \cite[Sections~1.3 and~2.3]{hong2018apo} for a discussion of the
connections to previous analyses of PCA for homoscedastic noise
(such as \cite{johnstone2009oca,nadler2008fsa,paul2007aos}),
and see \cite[Section~S1]{hong2018apo} for a discussion
of the connections to spiked covariance models.
Alternatively, \cite{hong2019ppf,hong2021hpp} consider a probabilistic PCA approach,
where the noise heteroscedasticity
is modeled via the statistical likelihood.
The resulting method is not a weighted PCA\@.
Instead,
one must solve a challenging optimization problem,
and \cite{hong2019ppf,hong2021hpp}
develop several algorithms for this purpose.

A closely related model for heterogeneous data
arises in the context
of low-rank clutter estimation
for RADAR.
In this setting,
the noise is homoscedastic
but the clutter signal
has heterogeneous strengths,
i.e.,
the clutter covariances
are a common low-rank matrix
scaled by heterogeneous power factors.
Maximum likelihood estimation
of the common low-rank matrix and the power factors
involves solving a challenging optimization problem,
and \cite{breloy2015cse,breloy2016rcm,sun2016lca,collas2021ppf}
develop efficient algorithms for this purpose.
The estimation performance is analyzed in \cite{besson2016bfa},
and \cite{abdallah2020bss} considers maximum a posteriori estimation.
This heterogeneous signal strength model is related
to the heteroscedastic noise model in the present paper
through a straightforward rescaling of the data:
scaling each sample by the inverse of its power factor
yields the model in the present paper.
Thus,
the optimally weighted PCA
developed here
can be straightforwardly modified to apply to the heterogeneous signal strength model;
see \insupp{sec:het:signal} for details.

Several recent works develop PCA variants
for high-dimensional data
with noise that is heteroscedastic across features.
In contrast to the samplewise heteroscedasticity we consider,
featurewise heteroscedasticity produces a
nonuniform bias along the diagonal of the covariance matrix
that skews its eigenvectors
even with infinitely many samples.
An approach based on
spectral shrinkage with noise whitening
(to make it homoscedastic)
is developed in \cite{leeb2021oss};
the noise is whitened by weighting the features
by their inverse noise variance.
Whitening both the features and samples
is considered in \cite{leeb2021mdf},
and whitening in the context of linearly transformed signals
is considered in \cite{dobriban2020opi}.
Alternatively,
\cite{zhang2022hpa} addresses the bias in the covariance matrix
by iteratively replacing its biased diagonal entries
using low-rank approximation.
Estimating the number of underlying principal components
is another important problem in this setting,
and recent works \cite{hong2020stn:arxiv:v1,ke2021eot,landa2021brt:arxiv:v2}
have developed new methods for tackling this challenge
under heteroscedastic noise.
\revone{%
Estimated principal components
are also often combined with
estimates of associated signal variances
to obtain estimates
of an underlying signal matrix or covariance.
For homoscedastic noise,
existing works have made tremendous progress
on how to estimate these signal variances
to optimize various objectives,
typically by applying a carefully designed shrinkage
to the eigenvalues of the sample covariance matrix;
see, e.g., \cite{donoho2018oso} and the references therein.
A few recent works address this question
in the context of heteroscedastic noise:
\cite{leeb2021oss} derives optimal shrinkages
for use with whitening,
\cite{leeb2021mdf} derives optimal spectral denoisers,
and
\cite{nadakuditi2014oaa} derives
an optimal data-driven shrinkage.
}%

Many works have considered weighted PCA methods in general;
see \cite[Section~14.2.1]{jolliffe2002pca}
for a survey of some of these works.
For example,
\cite[Sections~5.4--5.5]{deville1983dai}
discusses weighting features by inverse noise variance weights
to account for featurewise heteroscedasticity.
Weighting both samples and features
is proposed in \cite{cochran1977swp}
for analyzing spectrophotometric data
from scanning wavelength kinetics experiments;
the weights are again inverse noise variance.
Similar schemes have also been proposed
in metabolomics \cite{jansen2004aol}
and astronomy \cite{bailey2012pca,tamuz2005cse},
to name just a few areas.
Weighting data by inverse noise variance weights
has been a recurring theme.

Overall,
previous work on PCA for heteroscedastic noise
made significant progress on various important questions.
However, the important question of optimal weighting
was not previously considered.
This paper addresses that question
for noise that is heteroscedastic across samples.

\subsection{Organization of the paper}%
\Cref{sec:opt:weight} states our main result:
optimal weights and performance
for high-dimensional data with heteroscedastic noise.
\Cref{sec:opt:weight:sim} performs numerical simulations in finite dimensions,
and \cref{sec:weight:comp} compares the asymptotic optimal weights
with existing weighting schemes:
inverse noise variance weighted PCA,
PCA using only a single block of the data,
and unweighted PCA\@.
\revone{%
\Cref{sec:comp:more:methods} compares
optimally weighted PCA with some additional methods.
}%
\Cref{proof:main:result} proves the main result.
The optimal weights depend on the signal and noise variances.
\Cref{sec:signal:var:est} describes how estimates of these variances
can be used when the true variances are unknown.
\Cref{sec:exp:astro} illustrates optimally weighted PCA
on real data coming from astronomy.
Codes for reproducing
the figures in this paper are available online at:
\url{https://gitlab.com/dahong/optimally-weighted-pca-heteroscedastic-data}

For readers mostly interested in understanding
the underlying theory and proofs of the main result,
we suggest starting with
\cref{sec:opt:weight,proof:main:result}.
For readers mostly interested in
using optimally weighted PCA,
we suggest starting with
\revone{\cref{sec:opt:weight,sec:opt:weight:sim,sec:weight:comp,sec:comp:more:methods,sec:signal:var:est,sec:exp:astro}}.

%% file: content/main,result.tex
%!TEX root = ../wpca.tex

\section{Main result: optimal weights and performance}
\label{sec:opt:weight}%
We begin by making precise
the notion of optimal weights and optimal performance.
Consider a dataset $\bmY$
having $k$ underlying orthonormal components $\bmu_1,\dots,\bmu_k$,
where $\bmY$ is made of $L$ blocks $\bmY_1,\dots,\bmY_L$
of samples with heteroscedastic noise.
Then,
given weights $w_1,\dots,w_L \geq 0$,
the $\bmw$-weighted PCA estimate of the \ith component $\bmu_i$ from $\bmY$,
denoted $\bhtu_i(\bmw,\bmY)$, is
\begin{equation}
  \bhtu_i(\bmw,\bmY)
  \coloneqq
  \text{\ith leading eigenvector of the weighted sample covariance }
  \sum_{\ell = 1}^L w_\ell \bmY_\ell \bmY_\ell^\CT
  .
\end{equation}
A natural way to measure the performance of the estimate,
i.e., how well $\bhtu_i(\bmw,\bmY)$ recovers the \ith component $\bmu_i$,
is by the square inner product $r_i(\bmw,\bmY)$
given by
\begin{equation} \label{eq:rec}
  r_i(\bmw,\bmY) \coloneqq |\bmu_i^\CT \bhtu_i(\bmw,\bmY)|^2
  .
\end{equation}
Finally, optimal weights $\bmw^{\opt}_i(\bmY)$
and the optimal performance $r^{\opt}_i(\bmY)$
for the \ith component
are defined by
\begin{align} \label{eq:opt:weight:oracle}
  \bmw^{\opt}_i(\bmY)
  &\in
  \argmax_{\bmw}
  \; r_i(\bmw,\bmY)
  , &
  r^{\opt}_i(\bmY)
  &=
  \max_{\bmw}
  \; r_i(\bmw,\bmY)
  .
\end{align}

Note that the performance \cref{eq:rec}
depends on the underlying component $\bmu_i$.
However, in practice, $\bmu_i$ is of course unknown
so the optimization \cref{eq:opt:weight:oracle} cannot be done.
Fortunately,
as our main result below show\revone{s},
\revone{the optimal weights $\bmw^{\opt}_i(\bmY)$ and optimal performance $r^{\opt}_i(\bmY)$}
can be predicted when
the data
\ref{assump:stat} satisfies some natural statistical assumptions,
and
\ref{assump:asymp} grows large in size,
i.e., under the following setting.

\medskip
\noindent
\textbf{Setting:}
We will assume the following setting throughout the remainder of the paper.
\begin{enumerate}[label=(\alph*)]

\item \label{assump:stat}
The noisy data blocks
$\bmY_1 \in \bbC^{d \times n_1},\dots,\bmY_L \in \bbC^{d \times n_L}$ are generated
from the components $\bmu_1,\dots,\bmu_k \in \bbC^d$
with corresponding signal variances $\lambda_1 > \cdots > \lambda_k > 0$ as follows:
\begin{equation} \label{eq:model:blocks}
  \bmY_\ell = \bmF \bmZ_\ell + \bmE_\ell \in \bbC^{d \times n_\ell}
  , \quad \text{for } \ell = 1,\dots,L
  ,
\end{equation}
where
\begin{itemize}
  \item $\bmF \coloneqq [\sqrt{\lambda_1}\bmu_1,\dots,\sqrt{\lambda_k}\bmu_k] \in \bbC^{d \times k}$
  is a deterministic factor matrix common to all the blocks,
  \item $\bmZ_\ell \in \bbC^{k \times n_\ell}$
  is a coefficient matrix with IID entries having zero mean and unit variance,
  \item $\bmE_\ell \in \bbC^{d \times n_\ell}$
  is a noise matrix with IID entries having zero mean and variance $v_\ell > 0$,
\end{itemize}
and the noise entries further satisfy
a technical condition:
bounded $a$-th moment with $a > 4$,
i.e., $\exists_{a > 4} \text{ s.t. } \bbE|(\bmE_\ell)_{i,j}|^a < \infty$.%
\footnote{This technical condition on the noise
is satisfied by numerous distributions
including the sub-Gaussian and sub-Exponential families
\cite[Propositions 2.5.2 and 2.7.1]{vershynin2018hdp}.}
Note that this model also includes real-valued data
with real-valued coefficients and noise.

\item \label{assump:asymp}
The number of features $d$
and numbers of samples $n_1,\dots,n_L$
all grow towards infinity
but with fixed aspect ratios $n_\ell / d = c_\ell > 0$.
This asymptotic regime captures datasets where
the number of features and samples are roughly comparable,
as is common in modern big data settings.

\end{enumerate}
Note that under the model \cref{eq:model:blocks},
the optimal weights and performance \cref{eq:opt:weight:oracle}
are random quantities
so their convergence
will be probabilistic.
Specifically,
the convergence holds with probability one,
i.e., it is \emph{almost sure convergence},
which we will denote by $\asto$.

\strut\\\noindent
We are now ready to state the main result on the optimal weights and performance.

\begin{theorem}[Asymptotic optimal weights and performance]%
  \label{thm:opt:weight:asymp}%
  The optimal weights $\bmw^{\opt}_i(\bmY)$
  and corresponding optimal performance
  $r^{\opt}_i(\bmY)$
  converge as
  \begin{alignat}{4}
    \label{eq:opt:weight:asymp}
    \bmw^{\opt}_i(\bmY)
    &\quad&
    &\asto
    &\quad&
    &
    \bbrw^{\opt}_i
    &
    \coloneqq
    \bigg(
      \frac{1}{v_1}\frac{1}{1 + v_1/\lambda_i},
      \dots,
      \frac{1}{v_L}\frac{1}{1 + v_L/\lambda_i}
    \bigg)
    \qquad
    \text{up to scaling}
    , \\[0.5em]
    \label{eq:opt:rec:asymp}
    r^{\opt}_i(\bmY)
    &&
    &\asto
    &&
    &
    \brr^{\opt}_i
    &
    \coloneqq
    \text{the unique solution $x \in (0,1)$ of }
    \sum_{\ell=1}^L
      \frac{c_\ell}{v_\ell/\lambda_i}
      \,\frac{1-x}{v_\ell/\lambda_i+x}
    = 1
    ,
  \end{alignat}
  except when
  $\sum_{\ell=1}^L c_\ell (\lambda_i /v_\ell)^2 \leq 1$,
  in which case
  $\bmu_i$ is asymptotically unrecoverable
  by any weighted PCA,
  i.e., $r_i(\bmw,\bmY) \asto 0$ for all weights $\bmw$.
\end{theorem}

\begin{figure}[h] \centering
  \begin{tikzpicture}
    \begin{axis}[
      xlabel={$v_2/v_1$},
      xmin=1, xmax=4, domain=1:4, samples=40,
      xtick={1,...,4},
      minor xtick={1,1.25,...,4},
      ylabel={$w_2/w_1$}, ylabel style={rotate=-90},
      ymin=0, ymax=1,
      ytick={0,0.25,...,1}, yticklabels={$0$,$1/4$,$1/2$,$3/4$,$1$},
      minor ytick={0,0.0625,...,1},
      width=11cm, height=5.25cm,
      legend pos=north east,
      legend cell align=left,
      legend style={draw=none},
      cycle list/Set2,
      clip=false,
    ]
      % Optimal region
      \begin{scope}[font=\footnotesize]
        \addplot[forget plot,draw=none,name path=A] {1/x};
        \addplot[forget plot,draw=none,name path=B] {1/x^2};
        \addplot[forget plot,Gray,opacity=0.15] fill between [of=A and B];

        \draw[Gray,thick,decorate,decoration={brace,raise=1pt,amplitude=5pt}] (4.01,1/4-0.005) -- (4.01,1/16+0.005)
          node [pos=0.5,anchor=west,xshift=6pt] {\shortstack[l]{range of asymptotic\\optimal weights}};
      \end{scope}

      % Optimal weights
      \begin{scope}[very thick]
        \addplot {1/(x*(1-(1-x)/(1+4.00)))}; \addlegendentry{$\lambda_i/v_1 = 4$};
        \addplot {1/(x*(1-(1-x)/(1+1.00)))}; \addlegendentry{$\lambda_i/v_1 = 1$};
        \addplot {1/(x*(1-(1-x)/(1+0.25)))}; \addlegendentry{$\lambda_i/v_1 = 1/4$};
      \end{scope}

      % Other weights
      \begin{scope}[Black,dashed,ultra thick,font=\footnotesize]
        \addplot[forget plot] {1/x}   node [pos=0.30,anchor=south west,xshift=-0.3em]                {inverse noise variance};
        \addplot[forget plot] {1/x^2} node [pos=0.40,anchor=north east,xshift=+0.5em,yshift=+0.25em] {square inv. noise var.};
        \addplot[forget plot] {1}     node [pos=0.50,anchor=north]                                   {unweighted};
        \addplot[forget plot] {0}     node [pos=0.30,anchor=south,yshift=-0.25em]                    {discard $\bmY_2$ (noisier block)};
      \end{scope}

    \end{axis}
  \end{tikzpicture}
  \caption{%
    Relative weight $w_2/w_1$
    given by the optimal weights \cref{eq:opt:weight:asymp}
    as a function of the relative noise variance $v_2/v_1$
    for various signal-to-noise ratios $\lambda_i/v_1$.
    The optimal weights downweight noisier data
    more aggressively than inverse noise variance weights,
    but also do not discard noisier data.
    They lie in the region between inverse and square inverse noise variance weights.
  }
  \label{fig:opt:weights:asymp}
\end{figure}

\begin{remark}[Optimal weights downweight more than inverse noise variance weights]%
  \label{rem:opt:weight:range}%
  The optimal weights \cref{eq:opt:weight:asymp} are
  not the inverse noise variance weights that are commonly used.
  As illustrated in \cref{fig:opt:weights:asymp},
  optimal weights downweight noisier data more aggressively,
  but never discard data.
  Specifically, when the signal-to-noise ratio $\lambda_i/v_\ell$ is small,
  the optimal weights are square inverse noise variance weights up to scale.
  As $\lambda_i/v_\ell$ grows,
  the optimal weights gradually become less aggressive
  and approach inverse noise variance weights.
\end{remark}

\begin{remark}[Using estimated signal and noise variances]%
  \label{rem:signal:var:est}%
  The optimal weights \cref{eq:opt:weight:asymp}
  depend on the noise variances $\bmv$
  and on the signal component variance $\lambda_i$.
  In some settings,
  these parameters are known,
  e.g., from calibration data.
  When they are unknown,
  one may estimate them using existing ideas and approaches,
  then plug them in to obtain estimated weights.
  As discussed in \cref{sec:signal:var:est},
  these estimated weights are also asymptotically optimal
  as long as the variance estimates are consistent.
\end{remark}

\begin{remark}[Heterogeneous signal strengths]%
  \label{rem:het:signal}%
  The result may at first appear to be limited to data with homogeneous signal strengths.
  However, it generalizes straightforwardly to the case
  where the signal in each block is scaled by an associated signal strength,
  as arises, e.g., in RADAR applications
  \cite{breloy2015cse,breloy2016rcm,sun2016lca}.
  Simply preweight the data
  to recover the model \cref{eq:model:blocks};
  see \insupp{sec:het:signal} for a detailed description.
\end{remark}

\begin{remark}[Handling potentially degenerate cases]%
  \label{rem:undef:limsup}%
  A careful reader may note the subtle and technical point that
  there may exist degenerate choices of $\bmw$
  for which $r_i(\bmw,\bmY)$ is not well-defined,
  e.g., if the \ith leading eigenvector
  becomes undefined due to eigenvalue multiplicity.
  At such points, we define $r_i(\bmw,\bmY)$ by its $\limsup$ over $\bmw$.
  Doing so makes $r_i(\bmw,\bmY)$ upper semi-continuous in $\bmw$
  and avoids degenerate situations where its maximum does not exist.
\end{remark}

\revone{%
\begin{remark}[Nonorthogonality of estimated components]%
  \label{rem:nonorthogonality}%
  Since the optimal weights \cref{eq:opt:weight:asymp} are component-specific,
  the components
  $\bhtu_1(\bbrw^{\opt}_1,\bmY),\dots,\bhtu_k(\bbrw^{\opt}_k,\bmY)$
  estimated by optimally weighted PCA
  may not be orthogonal in practice.
  In applications where orthogonality is crucial,
  one option is to sacrifice component-wise optimality
  and use a single set of weights,
  e.g.,
  that just optimizes recovery of the weakest component
  or
  that optimizes some appropriate overall metric of performance.
  Alternatively,
  in many such cases,
  the principal subspace is of greater interest
  than the individual components;
  in these cases,
  one could orthogonalize the components,
  e.g., via Gram-Schmidt.
\end{remark}

\begin{remark}[Phase transition]%
  \label{rem:phase:transition}%
  Analogous to unweighted PCA under homoscedastic noise,
  optimally weighted PCA exhibits a phase transition between settings
  with zero asymptotic performance
  and those with nonzero asymptotic performance.
  As described in \cref{thm:opt:weight:asymp},
  optimally weighted PCA
  has nonzero asymptotic performance
  when $\sum_{\ell=1}^L c_\ell (\lambda_i /v_\ell)^2 > 1$
  (or in other words,
  $\lambda_i > (\sum_{\ell=1}^L c_\ell/v_\ell^2)^{-1/2}$).
  Notably,
  if any weighting scheme has nonzero asymptotic performance,
  then optimally weighted PCA does too,
  as illustrated in \cref{sec:weight:comp:phase}.
\end{remark}
}%

Before proving the main result (\cref{thm:opt:weight:asymp})
in \cref{proof:main:result},
we provide some more intuition
about it
through numerical simulations in finite dimensions (\cref{sec:opt:weight:sim})
and comparisons with existing weighting schemes (\cref{sec:weight:comp}).

%% file: content/simulation.tex
%!TEX root = ../wpca.tex

\section{Numerical simulation} \label{sec:opt:weight:sim}%
This section performs numerical simulations in finite dimensions.
Specifically, we generate
$L=2$ blocks of data
$\bmY_1 \in \bbR^{d \times n_1}$ and $\bmY_2 \in \bbR^{d \times n_2}$
according to the model \cref{eq:model:blocks}
with
\begin{itemize}
  \item $k=1$ component $\bmu_1 \in \bbR^d$ uniformly drawn from the unit sphere,
  \item Gaussian coefficients $(\bmZ_\ell)_{ij} \overset{\iid}{\sim} \clN(0,1)$
  and noise entries $(\bmE_\ell)_{ij} \overset{\iid}{\sim} \clN(0,v_\ell)$,
  \item component variance $\lambda_1 = 1$
  and noise variances $\bmv = (1,3)$.
\end{itemize}
\Cref{fig:opt:weight:sim}
shows the nonasymptotic empirical distributions
of the optimal weights $\bmw^{\opt}_1(\bmY)$
and the corresponding optimal performance $r^{\opt}_i(\bmY)$
from \cref{eq:opt:weight:oracle}
obtained using the true underlying component $\bmu_1$.
Since the weights are meaningful only up to scale,
we show the optimal relative weight $w^{\opt}_{1,2}(\bmY)/w^{\opt}_{1,1}(\bmY)$,
where $w^{\opt}_{i,\ell}(\bmY)$ is the $\ell$th entry of $\bmw^{\opt}_i(\bmY)$.
Similarly, $\brw^{\opt}_{i,\ell}$ denotes the $\ell$th entry of $\bbrw^{\opt}_i$.

\begin{figure}[h] \centering
  \subfloat[Nonasymptotic distributions for $d = 100$, $n_1 = 400$ and $n_2 = 800$. \label{fig:opt:weight:sim:100}]{%
    \begin{tikzpicture}
      \tikzmath{
        \vrat=3; \lmb=1;
        \wlim=(1*(1+\lmb))/(\vrat*(\vrat+\lmb));
        \rlim=0.6782867492531689;
      };
      \begin{axis}[
        xlabel={Optimized relative weight $w^{\opt}_{1,2}(\bmY)/w^{\opt}_{1,1}(\bmY)$},
        xmin=0, xmax=0.55,
        xtick={0,0.16667,0.33333,0.5},
        xticklabels={0,1/6,1/3,1/2},
        ymin=0, ymax=25,
        ytick={0,5,...,25},
        width=0.5\linewidth, height=0.2\linewidth,
        font=\footnotesize,
        clip=false,
      ]
        \addplot[ybar interval,fill opacity=0.5,fill=RoyalBlue]
          table {figs/opt,weight,sim,weights,c-4-8,v-1-3,lambda-1,d-100,ntrials-10000,ratiores-0.01.dat};
        \addplot[ultra thick,Dark2-B,dashed] coordinates {(\wlim,0) (\wlim,26)}
          node[anchor=south,xshift=1em] {limit $\brw^{\opt}_{1,2}/\brw^{\opt}_{1,1}$ from \cref{eq:opt:weight:asymp}};
        \end{axis}
    \end{tikzpicture}
    \qquad
    \begin{tikzpicture}
      \tikzmath{
        \vrat=3; \lmb=1;
        \wlim=(1*(1+\lmb))/(\vrat*(\vrat+\lmb));
        \rlim=0.6782867492531689;
      };
      \begin{axis}[
        xlabel={Optimized performance $r^{\opt}_1(\bmY)$},
        xmin=0, xmax=1,
        xtick={0,0.25,0.5,0.75,1},
        xticklabels={0,1/4,1/2,3/4,1},
        ymin=0, ymax=25,
        ytick={0,5,...,25},
        width=0.5\linewidth, height=0.2\linewidth,
        font=\footnotesize,
        clip=false,
      ]
        \addplot[ybar interval,fill opacity=0.5,fill=RoyalBlue]
          table {figs/opt,weight,sim,rec,c-4-8,v-1-3,lambda-1,d-100,ntrials-10000,ratiores-0.01.dat};
        \addplot[ultra thick,Dark2-B,dashed] coordinates {(\rlim,0) (\rlim,26)}
          node[anchor=south] {limit $\brr^{\opt}_1$ from \cref{eq:opt:rec:asymp}};
      \end{axis}
    \end{tikzpicture}%
  }

  \subfloat[Nonasymptotic distributions for $d = 1000$, $n_1 = 4000$ and $n_2 = 8000$. \label{fig:opt:weight:sim:1000}]{%
    \begin{tikzpicture}
      \tikzmath{
        \vrat=3; \lmb=1;
        \wlim=(1*(1+\lmb))/(\vrat*(\vrat+\lmb));
        \rlim=0.6782867492531689;
      };
      \begin{axis}[
        xlabel={Optimized relative weight $w^{\opt}_{1,2}(\bmY)/w^{\opt}_{1,1}(\bmY)$},
        xmin=0, xmax=0.55,
        xtick={0,0.16667,0.33333,0.5},
        xticklabels={0,1/6,1/3,1/2},
        ymin=0, ymax=25,
        ytick={0,5,...,25},
        width=0.5\linewidth, height=0.2\linewidth,
        font=\footnotesize,
        clip=false,
      ]
        \addplot[ybar interval,fill opacity=0.5,fill=RoyalBlue]
          table {figs/opt,weight,sim,weights,c-4-8,v-1-3,lambda-1,d-1000,ntrials-500,ratiores-0.01.dat};
        \addplot[ultra thick,Dark2-B,dashed] coordinates {(\wlim,0) (\wlim,26)}
          node[anchor=south,xshift=1em] {limit $\brw^{\opt}_{1,2}/\brw^{\opt}_{1,1}$ from \cref{eq:opt:weight:asymp}};
      \end{axis}
    \end{tikzpicture}
    \qquad
    \begin{tikzpicture}
      \tikzmath{
        \vrat=3; \lmb=1;
        \wlim=(1*(1+\lmb))/(\vrat*(\vrat+\lmb));
        \rlim=0.6782867492531689;
      };
      \begin{axis}[
        xlabel={Optimized performance $r^{\opt}_1(\bmY)$},
        xmin=0, xmax=1,
        xtick={0,0.25,0.5,0.75,1},
        xticklabels={0,1/4,1/2,3/4,1},
        ymin=0, ymax=25,
        ytick={0,5,...,25},
        width=0.5\linewidth, height=0.2\linewidth,
        font=\footnotesize,
        clip=false,
      ]
        \addplot[ybar interval,fill opacity=0.5,fill=RoyalBlue]
          table {figs/opt,weight,sim,rec,c-4-8,v-1-3,lambda-1,d-1000,ntrials-500,ratiores-0.01.dat};
        \addplot[ultra thick,Dark2-B,dashed] coordinates {(\rlim,0) (\rlim,26)}
          node[anchor=south] {limit $\brr^{\opt}_1$ from \cref{eq:opt:rec:asymp}};
      \end{axis}
    \end{tikzpicture}%
  }
  \caption{Nonasymptotic empirical distributions of optimal weights $\bmw^{\opt}_1(\bmY)$
    and optimal performance $r^{\opt}_1(\bmY)$
    from \cref{eq:opt:weight:oracle}
    for an illustrative example with two blocks of data
    $\bmY_1 \in \bbR^{d \times n_1}$
    and
    $\bmY_2 \in \bbR^{d \times n_2}$
    generated with noise variances $v_1 = 1$ and $v_2 = 3$,
    and with one underlying component having variance $\lambda_1 = 1$.}
  \label{fig:opt:weight:sim}
\end{figure}

Note first that the nonasymptotic distributions
for both the optimal weights $\bmw^{\opt}_1(\bmY)$
and the optimal performance $r^{\opt}_1(\bmY)$
are generally centered around their respective theoretical limits
$\bbrw^{\opt}_1$ and $\brr^{\opt}_1$ from
\cref{eq:opt:weight:asymp,eq:opt:rec:asymp}.
Moreover, they concentrate as the data grows in size from
\cref{fig:opt:weight:sim:100} to \cref{fig:opt:weight:sim:1000}.
This illustrates the almost sure convergence of the
nonasymptotic optimal weights and performance
to their limits.

Naturally, one also wonders whether the asymptotic results
of \cref{thm:opt:weight:asymp}
can be used to choose optimal weights or predict optimal performance
for real data,
which are finite-dimensional.
These experiments demonstrate
that the asymptotic optimal weights \cref{eq:opt:weight:asymp}
and performance \cref{eq:opt:rec:asymp}
can indeed be applied to choose weights
that are often close to optimal
for finite-dimensional data
and to predict their corresponding performance.

%% file: content/comparison.tex
%!TEX root = ../wpca.tex

\section{Comparison with existing weighting schemes}
\label{sec:weight:comp}%
This section compares the asymptotic optimal weights of \cref{eq:opt:weight:asymp}
with existing weighting schemes.
To ease discussion,
we will focus on the case with only two blocks $\bmY_1$ and $\bmY_2$,
i.e., $L=2$,
but the same insights apply more generally.
The weighting schemes considered are:
\begin{alignat}{2}
  \label{eq:inv:weight}
  & \textnormal{inverse noise variance:}
  &\qquad
  \bmw &= (1/v_1,1/v_2)
  , \\
  \label{eq:only1:weight}
  & \textnormal{only use $\bmY_1$:}
  &
  \bmw &= (1,0)
  , \\
  \label{eq:only2:weight}
  & \textnormal{only use $\bmY_2$:}
  &
  \bmw &= (0,1)
  , \\
  \label{eq:unweight:weight}
  & \textnormal{unweighted:}
  &
  \bmw &= (1,1)
  ,
\end{alignat}
where the weights are, as always, meaningful only up to scale.
Note that using only $\bmY_1$ or $\bmY_2$
are both special cases of general binary weights
that discard blocks of data.

\subsection{Comparison of weights}
\label{sec:weight:comp:weights}%
This section compares how close the various weighting schemes
are to the distribution of the actual empirically optimized weights
$\bmw^{\opt}_i(\bmY)$ from \cref{eq:opt:weight:oracle}
in an illustrative example.
As in \cref{sec:opt:weight:sim},
$\bmY_1 \in \bbR^{d \times n_1}$ and $\bmY_2 \in \bbR^{d \times n_2}$
are generated from the model \cref{eq:model:blocks}
with Gaussian coefficients and noise,
and a single signal component $\bmu_1 \in \bbR^d$
drawn uniformly from the unit sphere.
The noise variances are $\bmv = (1,3)$,
and the data sizes are $d = 1000$, $n_1 = 4000$ and $n_2 = 8000$.

\begin{figure}[h] \centering
  \subfloat[Moderate signal component variance $\lambda_1 = 1$. \label{fig:weight:comp:weights:lam1}]{%
    \begin{tikzpicture}
      \tikzmath{ \vrat=3; \lmb=1; \winv=1/\vrat; \wopt=(1*(1+\lmb))/(\vrat*(\vrat+\lmb)); }
      \begin{axis}[
        xlabel={Optimized relative weight $w^{\opt}_{1,2}(\bmY)/w^{\opt}_{1,1}(\bmY)$},
        xmin=0, xmax=1,
        xtick={0,0.16667,0.33333,0.5,0.66667,0.83333,1},
        xticklabels={0,1/6,1/3,1/2,2/3,5/6,1},
        ymin=0, ymax=25,
        ytick={0,5,...,25},
        width=0.47\linewidth, height=0.2\linewidth,
        font=\footnotesize,
        cycle list/Dark2,
        clip=false,
      ]
        \addplot[ybar interval,fill=RoyalBlue,fill opacity=0.5,forget plot] table {figs/opt,weight,sim,weights,c-4-8,v-1-\vrat,lambda-\lmb,d-1000,ntrials-500,ratiores-0.01.dat};
        \begin{scope}[ultra thick,dashed]
          \addplot coordinates {(\winv,0) (\winv,26)} node[anchor=south west,shift={(-1.5em,-0.1em)}] {inverse noise var.};
          \addplot coordinates {(\wopt,0) (\wopt,32)} node[anchor=south west,shift={(-1.0em,-0.2em)}] {asymptotic optimal};
          \addplot coordinates {(    1,0) (    1,26)} node[anchor=south east,shift={( 1.5em,-0.3em)}] {unweighted};
          \addplot coordinates {(    0,0) (    0,26)} node[anchor=south     ,shift={(-1.0em,-0.1em)}] {\shortstack{only use $\bmY_1$\\(cleaner block)}};
        \end{scope}
      \end{axis}
    \end{tikzpicture}%
  }
  \hfill
  \subfloat[Large signal component variance $\lambda_1 = 30$. \label{fig:weight:comp:weights:lam30}]{%
    \begin{tikzpicture}
      \tikzmath{ \vrat=3; \lmb=30; \winv=1/\vrat; \wopt=(1*(1+\lmb))/(\vrat*(\vrat+\lmb)); }
      \begin{axis}[
        xlabel={Optimized relative weight $w^{\opt}_{1,2}(\bmY)/w^{\opt}_{1,1}(\bmY)$},
        xmin=0, xmax=1,
        xtick={0,0.16667,0.33333,0.5,0.66667,0.83333,1},
        xticklabels={0,1/6,1/3,1/2,2/3,5/6,1},
        ymin=0, ymax=25,
        ytick={0,5,...,25},
        width=0.47\linewidth, height=0.2\linewidth,
        font=\footnotesize,
        cycle list/Dark2,
        clip=false,
      ]
        \addplot[ybar interval,fill=RoyalBlue,fill opacity=0.5,forget plot] table {figs/opt,weight,sim,weights,c-4-8,v-1-\vrat,lambda-\lmb,d-1000,ntrials-500,ratiores-0.01.dat};
        \begin{scope}[ultra thick,dashed]
          \addplot coordinates {(\winv,0) (\winv,26)} node[anchor=south west,shift={(-0.5em,-0.1em)}] {inverse noise var.};
          \addplot coordinates {(\wopt,0) (\wopt,32)} node[anchor=south     ,shift={( 2.0em,-0.2em)}] {asymptotic optimal};
          \addplot coordinates {(    1,0) (    1,26)} node[anchor=south east,shift={( 1.5em,-0.3em)}] {unweighted};
          \addplot coordinates {(    0,0) (    0,26)} node[anchor=south     ,shift={(-0.2em,-0.1em)}] {\shortstack{only use $\bmY_1$\\(cleaner block)}};
        \end{scope}
      \end{axis}
    \end{tikzpicture}%
  }
  \caption{Comparison of weighting schemes
    for an illustrative example with two blocks of data
    $\bmY_1 \in \bbR^{d \times n_1}$
    and
    $\bmY_2 \in \bbR^{d \times n_2}$
    generated with noise variances $v_1 = 1$ and $v_2 = 3$,
    where $d = 1000$, $n_1 = 4000$ and $n_2 = 8000$.
    The nonasymptotic distributions
    of empirically optimized weights $\bmw^{\opt}_1(\bmY)$
    from \cref{eq:opt:weight:oracle}
    are shown as histograms,
    with the asymptotic optimal weights $\bbrw^{\opt}_1$
    from \cref{eq:opt:weight:asymp}
    and the existing weights \cref{eq:inv:weight,eq:only1:weight,eq:only2:weight,eq:unweight:weight}
    overlaid as lines
    (the weights \cref{eq:only2:weight} that only use $\bmY_2$
    do not appear since they are at infinity).
  }
  \label{fig:weight:comp:weights}
\end{figure}

\Cref{fig:weight:comp:weights} shows the nonasymptotic distribution
of the empirically optimized weights \cref{eq:opt:weight:oracle}
with the asymptotic optimal weights $\bbrw^{\opt}_1$
from \cref{eq:opt:weight:asymp}
and the existing weights \cref{eq:inv:weight,eq:only1:weight,eq:only2:weight,eq:unweight:weight}
overlaid.
\Cref{fig:weight:comp:weights:lam1}
shows a case with a moderate signal component variance $\lambda_1 = 1$,
and \cref{fig:weight:comp:weights:lam30} shows a case with a strong signal component $\lambda_1 = 30$.
Since the weights are meaningful only up to scale,
we show the relative weight
$w^{\opt}_{1,2}(\bmY)/w^{\opt}_{1,1}(\bmY)$
given to the noisier block $\bmY_2$,
where $w^{\opt}_{i,\ell}(\bmY)$ is the $\ell$th entry of $\bmw^{\opt}_i(\bmY)$,
as in \cref{fig:opt:weight:sim}.

We make the following observations
from \cref{fig:weight:comp:weights}:
\begin{itemize}
  \item As before, the optimal weights are centered
    around the asymptotic optimal weights.
  \item When the signal component variance is moderate,
    the optimal weights do not overlap with
    the existing weighting schemes.
    They more aggressively downweight
    the noisier block than inverse noise variance weights
    but also do not discard the noisier block.
  \item When the signal component variance is large,
    the optimal weights overlap with
    inverse noise variance weights.
\end{itemize}

\subsection{Comparison of performance}
\label{sec:weight:comp:rec}%
This section compares the various weighting schemes
in terms of their performance $r_i(\bmw,\bmY)$ in finite dimensions.
As in \cref{sec:opt:weight:sim},
$\bmY_1 \in \bbR^{d \times n_1}$ and $\bmY_2 \in \bbR^{d \times n_2}$
are generated from the model \cref{eq:model:blocks}
with Gaussian coefficients and noise,
and a single signal component $\bmu_1 \in \bbR^d$ drawn uniformly from the unit sphere.
The signal component variance is $\lambda_1 = 2$,
and the data sizes are $d = 10^3$, $n_1 = 10^3$ and $n_2 = 10^4$.
The first noise variance is $v_1 = 1$
and the second noise variance ranges from $v_2=1$ to $v_2=20$.
\Cref{fig:weight:comp:rec} shows the nonasymptotic distribution of performance
for the asymptotic optimal weights \cref{eq:opt:weight:asymp}
and the existing weights \cref{eq:inv:weight,eq:only1:weight,eq:only2:weight,eq:unweight:weight}.

\begin{figure}[h] \centering
  \begin{tikzpicture}
    \begin{axis}[
      xlabel={Block 2 noise variance $v_2$},
      xmin=1, xmax=20,
      xtick={1,2,...,20},
      ylabel={$r_1(\bmw,\bmY)$},
      ymin=0, ymax=1,
      ytick={0,0.25,0.5,0.75,1},
      yticklabels={0,1/4,1/2,3/4,1},
      width=0.9\linewidth, height=0.37\linewidth,
      font=\footnotesize,
      ylabel style={rotate=-90},
      cycle list/Dark2,
      clip=false,
    ]
      \foreach \meth in {invvar,optlim,unweight,onlyY1,onlyY2} {
        \pgfplotstableread{figs/weight,comp,recsweep,c-1-10,v1-1.0,v2range-1.0-20.0,lambda-2,d-1000,ntrials-400,\meth.dat}\loadedtable
        \addplot+[forget plot,draw=none,name path=q1] table [x=v21,y=q1] from \loadedtable ;  % Ribbon q1
        \addplot+[forget plot,draw=none,name path=q3] table [x=v21,y=q3] from \loadedtable ;  % Ribbon q3
        \addplot+[forget plot,opacity=0.25] fill between [of=q1 and q3];                      % Ribbon
        \addplot+[ultra thick] table [x=v21,y=avg] from \loadedtable ;                        % Mean
        \addplot[forget plot,very thick,dashed,black] table [x=v21,y=lim] from \loadedtable ; % Asymptotic behavior
      }

      \begin{scope}[font=\small]
        \node[anchor=south west,Dark2-A] at (axis cs:14.0,0.15) {inverse noise variance};
        \node[anchor=south west,Dark2-B] at (axis cs: 6.0,0.62) {asymptotic optimal};
        \node[anchor=south west,Dark2-C] at (axis cs: 7.0,0.05) {unweighted};
        \node[anchor=north east,Dark2-D] at (axis cs:20.0,0.48) {only use $\bmY_1$ (cleaner block)};
        \node[anchor=north east,Dark2-E] at (axis cs: 5.0,0.30) {\shortstack[r]{only use $\bmY_2$\\(larger block)}};
      \end{scope}
    \end{axis}
  \end{tikzpicture}
  \caption{Performance comparison of weighting schemes
    \cref{eq:opt:weight:asymp,eq:inv:weight,eq:only1:weight,eq:only2:weight,eq:unweight:weight}
    for an illustrative example with two blocks of data
    $\bmY_1 \in \bbR^{d \times n_1}$
    and
    $\bmY_2 \in \bbR^{d \times n_2}$
    with $d = 10^3$, $n_1 = 10^3$, $n_2 = 10^4$,
    and signal component variance $\lambda_1 = 2$.
    The first block has noise variance $v_1 = 1$,
    while the second block noise variance ranges from $v_2=1$ to $v_2=20$.
    For each weighting scheme,
    the solid colored curve is the average from 400 trials,
    the ribbon indicates the corresponding interquartile interval,
    and the dashed black curve is the asymptotic performance
    from \cref{thm:opt:weight:asymp,prop:existing:rec:asymp}.}
    \label{fig:weight:comp:rec}
\end{figure}

We make the following observations
from \cref{fig:weight:comp:rec}:
\begin{itemize}
  \item Across the entire sweep,
    the asymptotic optimal weights are generally best.
  \item The performance of the asymptotic optimal weights
    is well predicted by the asymptotic performance
    from \cref{thm:opt:weight:asymp}.
  \item When $v_2$ is small, there is a lot of clean data coming from $\bmY_2$
    since $n_2$ is fairly large.
    All weighting schemes do well except the scheme of using only $\bmY_1$.
  \item As $v_2$ grows, $\bmY_2$ becomes noisier and the methods that use $\bmY_2$
    degrade in performance.
    Asymptotic optimal weighting degrades the most slowly / gracefully.
  \item As $v_2$ continues to grow, all methods that use $\bmY_2$
    eventually do worse than using only $\bmY_1$
    and hit zero asymptotic performance,
    \emph{except for the asymptotic optimal weighting}.
  \item When $v_2$ is large, asymptotic optimal weighting
    performs similarly to using only $\bmY_1$.
  \item Asymptotic optimal weights naturally transition from using $\bmY_2$
    to largely ignoring $\bmY_2$ without any tuning or manual choice of a ``cutoff''.
  \item Unweighted PCA and inverse noise variance weighted PCA
    sometimes perform worse when given more data;
    in some cases using only $\bmY_1$ or $\bmY_2$ was better.
    In contrast, asymptotic optimal weighting always performs
    better than using only $\bmY_1$ or only $\bmY_2$.
    Moreover, it always uses all data blocks, i.e., all weights are nonzero.
    With optimal weighting, more data can only help, it never hurts.
\end{itemize}

\Cref{fig:weight:comp:rec} overlaid the asymptotic performance for each weighting scheme.
For optimally weighted PCA, this limit was given in \cref{eq:opt:rec:asymp}.
The following \lcnamecref{prop:existing:rec:asymp}
provides an analogous result for the existing weighting schemes.

\begin{proposition}[Asymptotic performance of the existing weighting schemes]%
  \label{prop:existing:rec:asymp}%
  The weighting schemes \cref{eq:inv:weight,eq:only1:weight,eq:only2:weight,eq:unweight:weight}
  have corresponding performance converging as
  \begin{alignat}{4}
    \label{eq:inv:rec:asymp}
    & \textnormal{inverse noise variance:}
    &\qquad
    r_i(\bmw,\bmY)
    &\quad&
    &\asto
    &\quad&
    \max\bigg(0,\frac{c_{\phantom{0}} - \brv_{\phantom{0}}^2/\lambda_i^2}{c_{\phantom{0}} + \brv_{\phantom{0}}/\lambda_i}\bigg)
    , \\
    \label{eq:only1:rec:asymp}
    & \textnormal{only use $\bmY_1$:}
    &
    r_i(\bmw,\bmY)
    &&
    &\asto
    &&
    \max\bigg(0,\frac{c_1 - v_1^2/\lambda_i^2}{c_1 + v_1/\lambda_i}\bigg)
    , \\
    \label{eq:only2:rec:asymp}
    & \textnormal{only use $\bmY_2$:}
    &
    r_i(\bmw,\bmY)
    &&
    &\asto
    &&
    \max\bigg(0,\frac{c_2 - v_2^2/\lambda_i^2}{c_2 + v_2/\lambda_i}\bigg)
    , \\
    \label{eq:unweight:rec:asymp}
    & \textnormal{unweighted:}
    &
    r_i(\bmw,\bmY)
    &&
    &\asto
    &&
    \max\bigg(0,\frac{A(\beta_i)}{\beta_i B_i'(\beta_i)}\bigg)
    ,
  \end{alignat}
  where
  $c \coloneqq c_1 + \cdots + c_L$,
  $\brv \coloneqq (p_1/v_1 + \cdots + p_L/v_L)^{-1}$,
  $p_\ell \coloneqq c_\ell / c$,
  \begin{align*}
    A(x)   &\coloneqq 1 - \sum_{\ell=1}^L \frac{c_\ell v_\ell^2}{(x - v_\ell)^2}
    , &
    B_i(x) &\coloneqq 1 - \lambda_i \sum_{\ell=1}^L \frac{c_\ell}{x - v_\ell}
    ,
  \end{align*}
  and $\beta_i$ is the largest real root of the rational function $B_i$.
\end{proposition}

\Cref{prop:existing:rec:asymp}
is a by-product of \cref{lem:rec:lim}
in our proof of \cref{thm:opt:weight:asymp},
with some parts shown previously and some shown in this paper.
Specifically, \cref{eq:only1:rec:asymp,eq:only2:rec:asymp}
are exactly the well-studied homoscedastic case
\cite{johnstone2009oca,johnstone2018pih,nadler2008fsa,paul2007aos},
since the noise is homoscedastic when using only $\bmY_1$ or $\bmY_2$.
For unweighted PCA \cref{eq:unweight:rec:asymp},
\cite{hong2018apo} derived the performance for the case where $A(\beta_i) > 0$
and conjectured the behavior for $A(\beta_i) \leq 0$.
Finally, for the performance of inverse noise variance weights \cref{eq:inv:rec:asymp},
closely related results were contemporaneously derived in the recent work
\cite{leeb2021oss}.

\subsection{Comparison of phase transitions}
\label{sec:weight:comp:phase}%
The asymptotic performances
\cref{eq:opt:rec:asymp,eq:inv:rec:asymp,eq:only1:rec:asymp,eq:only2:rec:asymp,eq:unweight:rec:asymp}
of the various weighting schemes
exhibit phase transitions
between settings with zero asymptotic performance
and those with nonzero asymptotic performance.
Namely, each scheme has nonzero asymptotic performance
for data parameters $\bmc$, $\bmv$ and $\lambda_i$
in an associated regime.
\Cref{fig:weight:comp:phase} compares these regimes.

\begin{figure}[h] \centering
  \subfloat[With respect to combined aspect ratio $c$ for entire dataset~$\bmY$
    and signal component variance $\lambda_i$,
    for blocks with associated
    aspect ratios $\bmc = c \cdot (1/11,10/11)$
    and variances $\bmv=(1,5)$.
    \label{fig:weight:comp:phase:c:lam}]{%
    \begin{tikzpicture}
      \begin{axis}[
        xlabel={Combined aspect ratio $c$},
        xmin=0.8, xmax=6,
        ylabel={Signal component variance $\lambda_i$},
        ymin=0.9, ymax=3.2,
        width=0.6\linewidth, height=0.4\linewidth,
        font=\footnotesize,
        cycle list/Dark2,
      ]
        \tikzmath{ \rat=10; \v1=1; \v2=5;
          \p1=1/(1+\rat); \p2=\rat/(1+\rat); }
        \pgfplotsset{
          nzregion/.style = {forget plot,draw=none,fill,fill opacity=0.05},
          boundary/.style = {ultra thick}
        }
  
        \node [anchor=south west] at (axis cs:0.9,1) {\shortstack[l]{region of zero\\asymptotic performance}};

        % Inverse noise variance
        \tikzmath{ \constinv=(\p1/\v1 + \p2/\v2)^2;
          \cmininv=1/(\constinv*(\pgfkeysvalueof{/pgfplots/ymax})^2); \cmaxinv=\pgfkeysvalueof{/pgfplots/xmax}; }
        \addplot+ [nzregion,domain=\cmininv:\cmaxinv,samples=40] {1/sqrt(x*\constinv)} -- (\pgfkeysvalueof{/pgfplots/xmax},\pgfkeysvalueof{/pgfplots/ymax}) -- cycle;
        \addplot+ [boundary,domain=\cmininv:\cmaxinv,samples=40] {1/sqrt(x*\constinv)}
          node [pos=0.45,rotate=-24,anchor=south,yshift=0.1em,font=\scriptsize] {inverse noise variance};
  
        % Optimal
        \tikzmath{ \constopt=(\p1/\v1^2 + \p2/\v2^2);
          \cminopt=1/(\constopt*(\pgfkeysvalueof{/pgfplots/ymax})^2); \cmaxopt=\pgfkeysvalueof{/pgfplots/xmax}; }
        \addplot+ [nzregion,domain=\cminopt:\cmaxopt,samples=40] {1/sqrt(x*\constopt)} -- (\pgfkeysvalueof{/pgfplots/xmax},\pgfkeysvalueof{/pgfplots/ymax}) -- cycle;
        \addplot+ [boundary,domain=\cminopt:\cmaxopt,samples=40] {1/sqrt(x*\constopt)}
          node [pos=0.45,rotate=-24,anchor=north,yshift=0.1em,font=\scriptsize] {asymptotic optimal};
  
        % Uniform
        \pgfplotstableread{figs/weight,comp,phase,clam,r21-\rat,v-\v1-\v2,unweight.dat}\loadedtable
        \addplot+ [nzregion] table [x=c,y=lam] from \loadedtable -- (\pgfkeysvalueof{/pgfplots/xmax},\pgfkeysvalueof{/pgfplots/ymax}) -- cycle;
        \addplot+ [boundary,smooth] table [x=c,y=lam] from \loadedtable
          node [pos=0.36,rotate=-24,anchor=north,yshift=0.1em,font=\scriptsize] {unweighted};
  
        % Only 1
        \tikzmath{ \constonly1=\p1/\v1^2;
          \cminonly1=1/(\constonly1*(\pgfkeysvalueof{/pgfplots/ymax})^2); \cmaxonly1=\pgfkeysvalueof{/pgfplots/xmax}; }
        \addplot+ [nzregion,domain=\cminonly1:\cmaxonly1,samples=40] {1/sqrt(x*\constonly1)} -- (\pgfkeysvalueof{/pgfplots/xmax},\pgfkeysvalueof{/pgfplots/ymax}) -- cycle;
        \addplot+ [boundary,domain=\cminonly1:\cmaxonly1,samples=40] {1/sqrt(x*\constonly1)}
          node [pos=0.45,rotate=-24,anchor=north,yshift=0.2em,font=\scriptsize] {only use $\bmY_1$ (cleaner block)};
  
        % Only 2
        \tikzmath{ \constonly2=\p2/\v2^2;
          \cminonly2=1/(\constonly2*(\pgfkeysvalueof{/pgfplots/ymax})^2); \cmaxonly2=\pgfkeysvalueof{/pgfplots/xmax}; }
        \addplot+ [nzregion,domain=\cminonly2:\cmaxonly2,samples=40] {1/sqrt(x*\constonly2)} -- (\pgfkeysvalueof{/pgfplots/xmax},\pgfkeysvalueof{/pgfplots/ymax}) -- cycle;
        \addplot+ [boundary,domain=\cminonly2:\cmaxonly2,samples=40] {1/sqrt(x*\constonly2)}
          node [pos=0.40,rotate=-24,anchor=south,yshift=0.1em,font=\scriptsize] {only use $\bmY_2$ (larger block)};
      \end{axis}
    \end{tikzpicture}%
  }
  \hfill
  \subfloat[With respect to signal-to-noise ratios $\lambda_i/v_1$ and $\lambda_i/v_2$,
    for blocks with associated aspect ratios $\bmc=(1/2,1/2)$.
    \label{fig:weight:comp:phase:snr}]{%
    \begin{tikzpicture}
      \begin{axis}[
        xlabel={Signal-to-noise ratio $\lambda_i/v_1$ for $\bmY_1$},
        xmin=0, xmax=2.2,
        ylabel={Signal-to-noise ratio $\lambda_i/v_2$ for $\bmY_2$},
        ymin=0, ymax=2.2,
        axis equal,
        width=0.4\linewidth, height=0.4\linewidth,
        font=\footnotesize,
        cycle list/Dark2,
      ]
        \tikzmath{ \c=1; \rat=1;
          \p1=1/(1+\rat); \p2=\rat/(1+\rat); }
        \pgfplotsset{
          nzregion/.style = {forget plot,draw=none,fill,fill opacity=0.05},
          boundary/.style = {ultra thick}
        }
  
        \addplot[forget plot,dashed,opacity=0.2] coordinates {(0,0) (\pgfkeysvalueof{/pgfplots/xmax},\pgfkeysvalueof{/pgfplots/xmax})};
        \node [anchor=south west] at (axis cs:0.05,0.05) {\shortstack[l]{region of zero\\asymptotic\\performance}};
  
        % Inverse noise variance
        \tikzmath{ \invthresh1=1/(sqrt(\c)*\p1); \invthresh2=1/(sqrt(\c)*\p2); }
        \addplot+ [nzregion] coordinates {(0,\invthresh2) (\invthresh1,0)} -- (\pgfkeysvalueof{/pgfplots/xmax},0) -- (\pgfkeysvalueof{/pgfplots/xmax},\pgfkeysvalueof{/pgfplots/ymax}) -- (0,\pgfkeysvalueof{/pgfplots/ymax}) -- cycle;
        \addplot+ [boundary] coordinates {(0,\invthresh2) (\invthresh1,0)}
          node [pos=0.0,sloped,anchor=south west,shift={(-0.2em,-0.2em)},font=\scriptsize] {inv noise var};
  
        % Optimal
        \tikzmath{ \only1=1/sqrt(\c*\p1); \only2=1/sqrt(\c*\p2); }
        \addplot+ [nzregion,domain=90:0,samples=20] ({\only1*cos(x)},{\only2*sin(x)}) -- (\pgfkeysvalueof{/pgfplots/xmax},0) -- (\pgfkeysvalueof{/pgfplots/xmax},\pgfkeysvalueof{/pgfplots/ymax}) -- (0,\pgfkeysvalueof{/pgfplots/ymax}) -- cycle;
        \addplot+ [boundary,domain=90:0,samples=20] ({\only1*cos(x)},{\only2*sin(x)})
          node [pos=0.25,rotate=-24,anchor=north,yshift=-0.2em,font=\scriptsize] {asymp optimal};
  
        % Uniform
        \pgfplotstableread{figs/weight,comp,phase,snr,c-\c,r21-\rat,unweight.dat}\loadedtable
        \addplot+ [nzregion] table [x=lv1,y=lv2] from \loadedtable -- (\pgfkeysvalueof{/pgfplots/xmax},\pgfkeysvalueof{/pgfplots/ymax}) -- cycle;
        \addplot+ [boundary] table [x=lv1,y=lv2] from \loadedtable
          node [pos=0.0,anchor=south east,shift={(0.15em,-0.15em)},font=\scriptsize] {unweighted};
  
        % Only 1
        \addplot+ [nzregion] coordinates {(\only1,\pgfkeysvalueof{/pgfplots/ymax}) (\only1,0)} -- (\pgfkeysvalueof{/pgfplots/xmax},0) -- (\pgfkeysvalueof{/pgfplots/xmax},\pgfkeysvalueof{/pgfplots/ymax}) -- cycle;
        \addplot+ [boundary] coordinates {(\only1,\pgfkeysvalueof{/pgfplots/ymax}) (\only1,0)}
          node [pos=0.0,sloped,anchor=south west,yshift=-0.2em,font=\scriptsize] {only $\bmY_1$};
  
        % Only 2
        \addplot+ [nzregion] coordinates {(0,\only2) (\pgfkeysvalueof{/pgfplots/xmax},\only2)} -- (\pgfkeysvalueof{/pgfplots/xmax},\pgfkeysvalueof{/pgfplots/ymax}) -- (0,\pgfkeysvalueof{/pgfplots/ymax}) -- cycle;
        \addplot+ [boundary] coordinates {(0,\only2) (\pgfkeysvalueof{/pgfplots/xmax},\only2)}
          node [pos=1.0,sloped,anchor=south east,yshift=-0.2em,font=\scriptsize] {only $\bmY_2$};
      \end{axis}
    \end{tikzpicture}%
  }
  \caption{Comparison of phase transitions
    from \cref{thm:opt:weight:asymp,prop:existing:phase}
    for an illustrative example of two blocks of data
    $\bmY_1 \in \bbR^{d \times n_1}$
    and
    $\bmY_2 \in \bbR^{d \times n_2}$
    with variances $v_1$ and $v_2$
    and signal component variance $\lambda_i$.
    For each weighting scheme,
    asymptotic performance is zero below the phase transition
    and nonzero above it.
    In all cases, the shading goes up and to the right.
  }
  \label{fig:weight:comp:phase}
\end{figure}

We make the following observations
from \cref{fig:weight:comp:phase}:
\begin{itemize}
  \item None of the existing schemes dominates the rest
    with respect to nonzero asymptotic performance.
    In some cases one is better than another.
  \item Asymptotic optimal weighting dominates all of them.
    Whenever nonzero asymptotic performance is possible for one of the existing schemes,
    it is also possible for asymptotic optimal weighting.
  \item Asymptotic optimal weighting also achieves nonzero asymptotic performance
    in settings where all of the existing schemes have zero asymptotic performance.
\end{itemize}

For optimally weighted PCA,
the condition defining the regime
is given in \cref{thm:opt:weight:asymp}.
The following \lcnamecref{prop:existing:phase} gives analogous conditions
for the existing weighting schemes
and follows straightforwardly from \cref{prop:existing:rec:asymp}.

\begin{proposition}[Phase transitions of existing weighting schemes]%
  \label{prop:existing:phase}%
  The asymptotic performance
  of the weighting schemes are nonzero
  if and only if,
  respectively,
  \begin{alignat*}{2}
    & \textnormal{inverse noise variance:}
    &\qquad
    & c_{\phantom{0}} \cdot ( \lambda_i / \brv_{\phantom{0}} )^2 > 1
    , \\
    & \textnormal{only use $\bmY_1$:}
    &
    & c_1 \cdot ( \lambda_i / v_1 )^2 > 1
    , \\
    & \textnormal{only use $\bmY_2$:}
    &
    & c_2 \cdot ( \lambda_i / v_2 )^2 > 1
    , \\
    & \textnormal{unweighted:}
    &
    & A(\beta_i) > 0
    ,
  \end{alignat*}
  where $c$, $\brv$, $A$, and $\beta_i$ are as in \cref{prop:existing:rec:asymp}.
\end{proposition}

%% file: content/ext,comp,more,methods.tex
%!TEX root = ../wpca.tex

\revone{%
\section{Comparison with additional methods}
\label{sec:comp:more:methods}%

This section compares the performance of optimally weighted PCA
with additional PCA methods designed for some form of heteroscedastic noise.
Specifically,
we consider the following iterative methods:
\begin{itemize}
  \item HePPCAT \cite{hong2021hpp}
  is a probabilistic PCA approach
  that accounts for noise with samplewise heteroscedasticity
  by modeling the heteroscedasticity
  in the statistical likelihood.
  We used $1000$ iterations.

  \item HeteroPCA \cite{zhang2022hpa}
  addresses the bias
  in the diagonal of the covariance matrix
  caused by
  noise with featurewise heteroscedasticity.
  It does so by iteratively replacing the biased entries
  using low-rank approximation.
  We used $100$ iterations.
\end{itemize}
\Cref{fig:comp:more:methods} shows the nonasymptotic distribution of performance
for these methods
for the setup considered in \cref{sec:weight:comp:rec}:
$\bmY_1 \in \bbR^{d \times n_1}$ and $\bmY_2 \in \bbR^{d \times n_2}$
are generated from the model \cref{eq:model:blocks}
with Gaussian coefficients and noise,
and a single signal component $\bmu_1 \in \bbR^d$ drawn uniformly from the unit sphere.
The signal component variance is $\lambda_1 = 2$,
and the data sizes are $d = 10^3$, $n_1 = 10^3$ and $n_2 = 10^4$.
The first noise variance is $v_1 = 1$
and the second noise variance ranges from $v_2=1$ to $v_2=20$.
}%

\begin{figure}[h] \centering
  \begin{tikzpicture}
    \begin{axis}[
      xlabel={Block 2 noise variance $v_2$},
      xmin=1, xmax=20,
      xtick={1,2,...,20},
      ylabel={$|\bmu_1^\CT \bhtu_1(\bmY)|^2$},
      ymin=0, ymax=1,
      ytick={0,0.25,0.5,0.75,1},
      yticklabels={0,1/4,1/2,3/4,1},
      width=0.9\linewidth, height=0.37\linewidth,
      font=\footnotesize,
      ylabel style={rotate=-90},
      cycle list={{Dark2-A,solid},{Dark2-F,solid},{Dark2-B,dashed},{Dark2-C,dashed}},
      clip=false,
    ]
      \foreach \meth in {heppcat,heteropca,optwpca,unweight} {
        \pgfplotstableread{figs/comp,more,methods,recsweep,c-1-10,v1-1.0,v2range-1.0-20.0,lambda-2,d-1000,ntrials-400,\meth.dat}\loadedtable
        \addplot+[forget plot,draw=none,name path=q1] table [x=v21,y=q1] from \loadedtable ;  % Ribbon q1
        \addplot+[forget plot,draw=none,name path=q3] table [x=v21,y=q3] from \loadedtable ;  % Ribbon q3
        \addplot+[forget plot,opacity=0.25] fill between [of=q1 and q3];                      % Ribbon
        \addplot+[ultra thick] table [x=v21,y=avg] from \loadedtable ;                        % Mean
      }

      \begin{scope}[font=\small]
        \node[anchor=north east,Dark2-A] at (axis cs:19.0,0.50) {HePPCAT \cite{hong2021hpp}};
        \node[anchor=north east,Dark2-F] at (axis cs: 5.3,0.35) {HeteroPCA \cite{zhang2022hpa}};
        \node[anchor=south west,Dark2-B] at (axis cs: 6.0,0.62) {Optimally Weighted PCA};
        \node[anchor=south west,Dark2-C] at (axis cs: 6.8,0.08) {Unweighted PCA};
      \end{scope}
    \end{axis}
  \end{tikzpicture}
  \caption{\revone{Performance comparison with additional methods
    for the illustrative example in \cref{fig:weight:comp:rec}:
    two blocks of data
    $\bmY_1 \in \bbR^{d \times n_1}$
    and
    $\bmY_2 \in \bbR^{d \times n_2}$
    with $d = 10^3$, $n_1 = 10^3$, $n_2 = 10^4$,
    and signal component variance $\lambda_1 = 2$.
    The first block has noise variance $v_1 = 1$,
    while the second block noise variance ranges from $v_2=1$ to $v_2=20$.
    For each method,
    the colored curve is the average from 400 trials,
    and
    the ribbon indicates the corresponding interquartile interval.}}
    \label{fig:comp:more:methods}
\end{figure}

\revone{%
We make the following observations
from \cref{fig:comp:more:methods}:
\begin{itemize}
  \item HePPCAT accounts for the heterogeneous quality of the data blocks,
  and
  it performs very similarly to optimally weighted PCA in this case
  (the curves overlap).
  Note that
  HePPCAT involves solving a nonconvex optimization problem,
  and it currently lacks a guarantee of convergence
  to a global optimizer.
  In contrast,
  optimally weighted PCA can be computed
  simply and reliably
  via the well-studied singular value decomposition.

  \item HeteroPCA is also designed to account for heteroscedastic noise,
  but does so primarily for featurewise heteroscedasticity.
  It treats the samples uniformly
  and performs very similarly to unweighted PCA in this case
  (the curves overlap);
  it performs worse than optimally weighted PCA
  and is even eventually worse than using only $\bmY_1$.
  This example highlights how
  samplewise and featurewise
  heteroscedasticity in the noise differ.
  They have qualitatively different impacts
  that seem to call for distinct approaches.
\end{itemize}
}%

%% file: content/derivation.tex
%!TEX root = ../wpca.tex

\section{Proof of main result} \label{proof:main:result}%
\Cref{thm:opt:weight:asymp}
states that
unless $\sum_{\ell=1}^L c_\ell (\lambda_i /v_\ell)^2 \leq 1$,
the optimal weights $\bmw^{\opt}_i(\bmY)$
and corresponding optimal performance $r^{\opt}_i(\bmY)$
for the \ith component
converge almost surely to
$\bbrw^{\opt}_i$ (up to scale)
and
$\brr^{\opt}_i$
in
the right hand sides of \cref{eq:opt:weight:asymp,eq:opt:rec:asymp}.

Namely,
$\bbrw^{\opt}_i$ and $\brr^{\opt}_i$ are the result of
first optimizing (with respect to the weights $\bmw$)
then taking the limit (as the data $\bmY$ grows in size).
Unfortunately,
$\bmw^{\opt}_i(\bmY)$ and $r^{\opt}_i(\bmY)$
are complicated functions,
making it challenging to directly analyze their limit.
So, we instead first take the limit then optimize.
\revone{%
More precisely,
using $\aslim$ to denote almost sure limits
and
writing $\bmY^{(d)}$ to make the limits more explicit,
we first derive
\begin{equation*}
  \argmax_{\bmw} \brr_i(\bmw)
  \quad \text{where} \quad
  \brr_i(\bmw) = \aslim_{d \to \infty} r_i(\bmw,\bmY^{(d)})
  ,
\end{equation*}
then use that result to obtain
the result we want,
i.e.,
\begin{equation*}
  \aslim_{d \to \infty} \bmw^{\opt}_i(\bmY^{(d)})
  \quad \text{where} \quad
  \bmw^{\opt}_i(\bmY^{(d)}) = \argmax_{\bmw} r_i(\bmw,\bmY^{(d)})
  .
\end{equation*}
}%
The following diagram illustrates the approach:
\begin{equation} \label{eq:conv:argmax:diagram}
  \begin{tikzcd}[row sep=4.8em,column sep=10em]
    r_i(\bmw,\bmY)
    \arrow[r, "\text{optimize w.r.t. $\bmw$}"' Blue, "\substack{\text{Definition} \\ \text{(\Cref{eq:opt:weight:oracle})}}"]
    \arrow[d, "\text{a.s. limit}" Blue, "\substack{\text{\Cref{sct:rec:lim}} \\ \text{(\cref{lem:rec:lim})}}"']
    &
    \bmw^{\opt}_i(\bmY), r^{\opt}_i(\bmY)
    \arrow[d, "\text{a.s. limit}"' Blue, "\substack{\text{Main result} \\ \text{(\cref{thm:opt:weight:asymp})}}"]
    \\
    \brr_i(\bmw)
    \arrow[r, "\text{optimize w.r.t. $\bmw$}" Blue, "\substack{\text{\Cref{sct:rec:lim:opt}} \\ \text{(\cref{lem:rec:lim:opt})}}"']
    &
    \bbrw^{\opt}_i, \brr^{\opt}_i
  \end{tikzcd}
\end{equation}
For this approach to work, the diagram must commute,
i.e., the optimizer of the almost sure limit (which we derive)
must match the almost sure limit of the optimizer (which we want).
The following \lcnamecref{lem:conv:argmax} states a suitable sufficient condition under which this happens,
i.e., the maximizer of the limit is the limit of the maximizer.
This \lcnamecref{lem:conv:argmax}
may be proved using techniques and results from variational analysis,
e.g., \cite[Chapter 7]{rockafellar1998va}.
For convenience,
\insupp{calc:proof:lem:conv:argmax}
also provides an elementary, self-contained, and concise proof.

\begin{lemma}[Diagram commutes]%
  \label{lem:conv:argmax}%
  Let $\clX \subseteq \bbR^d$ be compact,
  $f_n : \clX \to \bbR$ be a sequence of functions,
  and $f : \clX \to \bbR$ be a function,
  such that on $\clX$
  \begin{enumerate}
    \item each $f_n$ has a maximum,
    \item $f_n$ converges uniformly to $f$,
    \item $f$ is continuous,
    and
    \item $f$ has a unique maximizer.
  \end{enumerate}
  Then the maximum of $f_n$
  and the set of maximizers of $f_n$
  both converge,
  i.e.,
  \begin{align*}
    \max_{\bmx \in \clX} f_n(\bmx)
    \; \; &\to \; \;
    \max_{\bmx \in \clX} f(\bmx)
    , &
    \argmax_{\bmx \in \clX} f_n(\bmx)
    \; \; &\to \; \;
    \argmax_{\bmx \in \clX} f(\bmx)
    ,
  \end{align*}
  where the set convergence is with respect to the Hausdorff distance.
\end{lemma}

With \cref{lem:conv:argmax} in hand,
it now remains to:
a) derive the almost sure limit $\brr_i(\bmw)$
of $r_i(\bmw,\bmY)$,
and show that the convergence is uniform in $\bmw$ (\cref{lem:rec:lim} in
\cref{sct:rec:lim});
and b) optimize $\brr_i(\bmw)$
and show that the optimizer is unique
(up to scaling)
except when $\sum_{\ell=1}^L c_\ell (\lambda_i /v_\ell)^2 \leq 1$
(\cref{lem:rec:lim:opt} in \cref{sct:rec:lim:opt}).

\subsection{Almost sure limit of the performance}
\label{sct:rec:lim}%
This section derives the almost sure limit
of the performance $r_i(\bmw,\bmY)$,
where the convergence is uniform in the weights $\bmw$.

\begin{lemma}[Almost sure limit of the performance]%
  \label{lem:rec:lim}%
  For $i = 1,\dots,k$,
  \begin{equation} \label{eq:rec:lim}
    r_i(\bmw, \bmY)
    \quad
    \asto
    \quad
    \brr_i(\bmw)
    \coloneqq
    \max\bigg(
      0,
      \frac{1}{\beta_{i,\bmw}} \frac{A_{\bmw}(\beta_{i,\bmw})}{B_{i,\bmw}'(\beta_{i,\bmw})}
    \bigg)
    ,
  \end{equation}
  where the convergence is uniform with respect to $\bmw$ on $\bbR_{\geq 0}^L \setminus \{\bm0_L\}$,
  \begin{align} \label{eq:rec:lim:AB}
    A_{\bmw}(x) &\coloneqq
    1 - \sum_{\ell=1}^L \frac{c_\ell w_\ell^2 v_\ell^2}{(x - w_\ell v_\ell)^2}
    , &
    B_{i,\bmw}(x) &\coloneqq
    1 - \lambda_i \sum_{\ell=1}^L \frac{c_\ell w_\ell}{x - w_\ell v_\ell}
    ,
  \end{align}
  and $\beta_{i,\bmw}$ is the largest real root of $B_{i,\bmw}$.
\end{lemma}

The remainder of this subsection proves \cref{lem:rec:lim}.
After defining some notations,
we derive the limit of the singular values,
then derive the limit of $r_i(\bmw,\bmY)$
in two regimes (above and below the phase transition),
and finally derive the above algebraic form.
There are several other ways to structure
these derivations; see, e.g., \cite{bloemendal2014ill}.
The approach we take here
carefully combines the general perturbation approach of \cite{benaychgeorges2012tsv}
with celebrated random matrix theoretic results on local laws \cite{ding2018ana, knowles2016all, xi2020coe}
(reviewed in \insupp{sec:local:law})
to obtain uniform convergence for the singular values and vectors.
The algebraic form is derived following the approach in \cite{hong2018apo}.
Throughout the proof,
we postpone
some detailed calculations
to the \supporapp.

\subsubsection{Notation and preliminaries}
\label{proof:rec:lim:prelim}%
Let $\httheta_{i,\bmw}$, $\bhtu_{i,\bmw}$, and $\bhtq_{i,\bmw}$
denote, respectively,
the \ith singular value, \ith left singular vector and \ith right singular vector
of the normalized and weighted data matrix
\begin{equation} \label{p:eq:model:wmatrix}
  \btlY_{\bmw}
  \coloneqq
  \frac{1}{\sqrt{n}} \, [\sqrt{w_1}\,\bmY_1,\dots,\sqrt{w_L}\,\bmY_L]
  = \bmU \bmTheta \btlQ_{\bmw}^\CT + \btlE_{\bmw}
  ,
\end{equation}
where
$\bmU \coloneqq [\bmu_1,\dots,\bmu_k]$,
$\bmTheta \coloneqq \diag(\theta_1,\dots,\theta_k)$
has diagonal entries $\theta_i \coloneqq \sqrt{\lambda_i}$,
and
\begin{align*}
  \btlQ_{\bmw} &\coloneqq \frac{1}{\sqrt{n}} \, [\sqrt{w_1}\,\bmZ_1,\dots,\sqrt{w_L}\,\bmZ_L]^\CT
  , &
  \btlE_{\bmw} &\coloneqq \frac{1}{\sqrt{n}} \, [\sqrt{w_1}\,\bmE_1,\dots,\sqrt{w_L}\,\bmE_L]
  ,
\end{align*}
are normalized and weighted coefficients and noise.
We indicate that these are functions of the weights via the subscript $\bmw$,
and omit the dependence
of the singular values and vectors on the data (for brevity).
We also omit the index for $d$;
all limits are as $d \to \infty$
unless otherwise specified.

Note that the left singular vectors of $\btlY_{\bmw}$
are exactly the eigenvectors of
the weighted sample covariance,
so the performance can be written as
$r_i(\bmw,\bmY) = |\bmu_i^\CT \bhtu_{i,\bmw}|^2$.
The noise matrix $\btlE_{\bmw}$ satisfies
the usual random matrix theoretic conditions,
and is well known to have
a singular value distribution
that converges weakly almost surely
to a nonrandom compactly supported measure $\mu_{\bmw}$
whose Stieltjes transform
\begin{equation}
    \label{eq:stieltjes:def}
    \mplaw(\bmw,\zeta)
    \coloneqq
    \int \frac{1}{\zeta^2-t^2} \rmd \mu_{\bmw}(t)
\end{equation}
is the unique solution
to the generalized Marchenko-Pastur equation
\begin{equation}
  \label{eq_gen_MP}
  \frac{1}{\mplaw(\bmw,\zeta)}
  =
  \zeta^2
  -
  \sum_{\ell=1}^L \frac{p_\ell w_\ell v_\ell}{1 - w_\ell v_\ell \mplaw(\bmw,\zeta) / c}
  ,
\end{equation}
for which $\Im \mplaw(\bmw,\zeta) < 0$ for $\zeta^2$ in the upper half complex plane \cite{marcenko1967doe}.

Moreover,
the operator norm of the noise matrix converges
to the upper-edge of $\mu_{\bmw}$
(see \insupp{calc:conv:noise} for detailed derivation),
i.e.,
\begin{equation}
  \label{eq:conv:noise}
  \|\btlE_{\bmw}\|_{\opnorm}
  \asunifto{\bmw \in \Delta_L}
  b_{\bmw}
  \coloneqq
  \sup\{\text{support of $\mu_{\bmw}$}\}
  ,
\end{equation}
where $\Delta_L \coloneqq \{\bmw \in \bbR_{\geq 0}^L : w_1 + \cdots + w_L = 1\}$
is the probability simplex,
and $\asunifto{}$ denotes almost sure uniform convergence.

\subsubsection{Limits of the singular values}
\label{proof:analysis:singvals}%
Using a similar argument as \cite[Section 4]{benaychgeorges2012tsv},
one can show that
any singular value of $\btlY_{\bmw}$
that does not tend to a limit greater than $b_{\bmw}$
tends to $b_{\bmw}$,
so we focus on singular values with limits greater than $b_{\bmw}$.
Following the approach of \cite[Section 4]{benaychgeorges2012tsv},
this section studies those singular values
through the matrix-valued function
\begin{equation}
  \label{p:eq:singval:kernel}
  \bmM(\bmw,\zeta)
  \coloneqq
  \begin{bmatrix} \bmU & \\ & \btlQ_{\bmw} \end{bmatrix}^\CT
  \bmG(\bmw,\zeta)
  \begin{bmatrix} \bmU & \\ & \btlQ_{\bmw} \end{bmatrix}
  -
  \begin{bmatrix} & \bmTheta^{-1} \\ \bmTheta^{-1} & \end{bmatrix}
  ,
\end{equation}
where $\bmG(\bmw,\zeta)$
is the resolvent (or Green's function) defined as
\begin{equation}
  \label{eq_resolvent}
  \bmG(\bmw,\zeta)
  \coloneqq
  \Bigg(
    \zeta \bmI_{d+n}
    -
    \begin{bmatrix} & \btlE_{\bmw}   \\ \btlE_{\bmw}^\CT   & \end{bmatrix}
  \Bigg)^{-1}
  .
\end{equation}
The link to singular values
is made through an extension of \cite[Lemma 4.1]{benaychgeorges2012tsv}
that incorporates the weights
(see \insupp{calc:proof:lem:singval:kernel} for detailed derivation).
It states that
\begin{equation} \label{eq:singval:kernel}
  \forall_{\zeta > \|\btlE_{\bmw}\|_{\opnorm}}
  \quad
  \text{$\zeta$ is a singular value of $\btlY_{\bmw}$}
  \iff
  \det\bmM(\bmw,\zeta) = 0
  ,
\end{equation}
so we instead study $\bmM$ in the limit.
A careful application of anisotropic local laws \cite{ding2018ana, knowles2016all, xi2020coe}
(see \insupp{calc:lem:kernel:lims} for detailed calculations)
yields
that for any $\tau > 0$,
\begin{equation} \label{wtMas}
  \bmM(\bmw,\zeta)
  \quad
  \asunifto{(\bmw,\zeta) \in \Omega(\tau)}
  \quad
  \bbrM(\bmw,\zeta)
  \coloneqq
  \begin{bmatrix} \varphi_{1,\bmw}(\zeta) \bmI_k &  \\  & \varphi_{2,\bmw}(\zeta) \bmI_k \end{bmatrix}
  -
  \begin{bmatrix} & \bmTheta^{-1} \\ \bmTheta^{-1} & \end{bmatrix}
  ,
\end{equation}
where
$\Omega(\tau) \coloneqq \{(\bmw,\zeta) \in \Delta_L \times \bbC : (\Re\zeta,\Im\zeta) \in [b_{\bmw} + \tau,\infty) \times [-1 , 1]\}$
and
\begin{align}
  \label{p:eq:vp1}
  \varphi_{1,\bmw}(\zeta)
  &\coloneqq \zeta \mplaw(\bmw,\zeta)
  = \int \frac{\zeta}{\zeta^2-t^2} \; \rmd \mu_{\bmw}(t)
  , &
  \varphi_{2,\bmw}(\zeta) &\coloneqq \sum_{\ell=1}^L \frac{p_\ell w_\ell}{\zeta - w_\ell v_\ell \varphi_{1,\bmw}(\zeta)/c}
  .
\end{align}

Finally, we apply \cite[Lemma A.1]{benaychgeorges2012tsv}
with \cref{eq:singval:kernel}
and \cref{wtMas}
in the same way as~\cite[Section 4]{benaychgeorges2012tsv};
straightforward calculations (see \insupp{calc:traceconv:props})
verify that $\varphi_{1,\bmw}$ and $\varphi_{2,\bmw}$ satisfy
the conditions of \cite[Lemma A.1]{benaychgeorges2012tsv}.
Moreover, noting that these arguments extend
to uniform convergence in $\bmw \in \Delta_L$
yields that
\begin{equation} \label{p:eq:singval:limit}
  \httheta_{i,\bmw}
  \asunifto{\bmw \in \Delta_L}
  \brtheta_{i,\bmw}
  ,
  \quad \text{where} \quad
  \brtheta_{i,\bmw}
  \coloneqq
  \begin{cases}
    \rho_{i,\bmw} & \text{if } \theta_i > \tltheta_{\bmw} , \\
    b_{\bmw} & \text{otherwise} ,
  \end{cases}
\end{equation}
where
$D_{\bmw}(\zeta) \coloneqq \varphi_{1,\bmw}(\zeta) \varphi_{2,\bmw}(\zeta)$ for $\zeta > b_{\bmw}$,
$\rho_{i,\bmw} \coloneqq D_{\bmw}^{-1}(1/\theta_i^2)$, and
$\tltheta_{\bmw}^2 \coloneqq 1/D_{\bmw}(b_{\bmw}^+)$.
Here $D_{\bmw}^{-1}$ denotes the inverse function of $D_{\bmw}$, and $f(b^+) \coloneqq \lim_{\zeta \to b^+} f(\zeta)$
is the limit from above.

\subsubsection{Performance above the phase transition}
\label{proof:analysis:singvecs}%
This section derives the limit of
the performance $r_i(\bmw,\bmY)$
above the phase transition.
Namely,
for any $\nu > 0$,
we prove uniform convergence
with respect to $\bmw$ over the domain
$\clW_>(\nu)
\coloneqq
\{\bmw \in \Delta_L : \brtheta_{i,\bmw} > b_{\bmw} + \nu\}$.

Following the approach of \cite[Section 5]{benaychgeorges2012tsv},
we study $r_i(\bmw,\bmY) = |\bmu_i^\CT \bhtu_{i,\bmw}|^2$
through the following extension of
\cite[Lemma 5.1]{benaychgeorges2012tsv}
(see \insupp{calc:lem:singvec:kernel} for detailed derivation):
\begin{subequations}
\label{lem:singvec:kernel}
\begin{alignat}{6}
  \label{p:eq:singvec:kernel}
  &
  \forall_{\bmw \in \Delta_L}
  &\quad&
  \httheta_{i,\bmw} > \|\btlE_{\bmw}\|_{\opnorm}
  &\quad&
  \implies
  &\quad&
  &\bigg\{\;\;&
  0 =
  \bmM(\bmw,\httheta_{i,\bmw})
  \begin{bmatrix}
    \bmTheta \btlQ_{\bmw}^\CT \bhtq_{i,\bmw} \\
    \bmTheta \bmU^\CT \bhtu_{i,\bmw}
  \end{bmatrix}
  \text{ and}
  &\;&
  \\
  \label{p:eq:singvec:sum}
  &
  &&
  &&
  &&
  &&
  1 = \chi_1(\bmw) + \chi_2(\bmw) + 2 \Re \chi_3(\bmw)
  &&\;\;\bigg\}
  ,
\end{alignat}
\end{subequations}
where $\bmGamma_{\bmw} \coloneqq (\httheta_{i,\bmw}^2 \bmI_d - \btlE_{\bmw}\btlE_{\bmw}^\CT)^{-1}$
and
\begin{alignat}{13}
  \label{p:eq:singvec:sum:parts}
  \chi_1(\bmw)
  &\;\;&
  &\coloneqq
  &\;\;&
  \sum_{j_1,j_2 = 1}^k
  &\;&
  \theta_{j_1} \theta_{j_2}
  &\;&
  \cdot
  &\;&
  & (\btlq_{j_1,\bmw}^\CT \bhtq_{i,\bmw})
  &\;&
  &\cdot
  &\;&
  & (\btlq_{j_2,\bmw}^\CT \bhtq_{i,\bmw})^*
  &\;&
  &\cdot
  &\;&
  &
  \bmu_{j_2}^\CT \httheta_{i,\bmw}^2 &\bmGamma_{\bmw}^2 \bmu_{j_1}
  ,
  \\
  \nonumber
  \chi_2(\bmw)
  &&
  &\coloneqq
  &&
  \sum_{j_1,j_2 = 1}^k
  &&
  \theta_{j_1} \theta_{j_2}
  &&
  \cdot
  &&
  & (\bmu_{j_1}^\CT \bhtu_{i,\bmw})
  &&
  &\cdot
  &&
  & (\bmu_{j_2}^\CT \bhtu_{i,\bmw})^*
  &&
  &\cdot
  &&
  &
  \btlq_{j_2,\bmw}^\CT \btlE_{\bmw}^\CT &\bmGamma_{\bmw}^2 \btlE_{\bmw} \btlq_{j_1,\bmw}
  ,
  \\
  \nonumber
  \chi_3(\bmw)
  &&
  &\coloneqq
  &&
  \sum_{j_1,j_2 = 1}^k
  &&
  \theta_{j_1} \theta_{j_2}
  &&
  \cdot
  &&
  & (\bmu_{j_1}^\CT \bhtu_{i,\bmw})
  &&
  &\cdot
  &&
  & (\btlq_{j_2,\bmw}^\CT \bhtq_{i,\bmw})^*
  &&
  &\cdot
  &&
  &
  \bmu_{j_2}^\CT \httheta_{i,\bmw} &\bmGamma_{\bmw}^2 \btlE_{\bmw} \btlq_{j_1,\bmw}
  .
\end{alignat}
It follows from \cref{p:eq:singval:limit}
that almost surely, eventually,
for all $\bmw \in \clW_>(\nu)$,
$\httheta_{i,\bmw} > \|\btlE_{\bmw}\|_{\opnorm}$
and hence the condition of \cref{lem:singvec:kernel} holds.
It remains to study the limits of $\chi_1$, $\chi_2$, and $\chi_3$.

Carefully applying \cref{wtMas,p:eq:singval:limit}
to \cref{p:eq:singvec:kernel}
in a similar way as \cite[Section 5]{benaychgeorges2012tsv}
yields that (see \insupp{calc:lem:singvec:sum:parts:lim} for detailed calculations),
for any $\nu > 0$,
\begin{subequations}
\label{p:eq:singvec:sum:parts:lim:all}
\begin{alignat}{10}
  \label{p:eq:singvec:sum:parts:lim}
  \chi_1(\bmw)
  &\quad&
  -
  &\quad&
  &
  \theta_i^2
  \bigg[ +\frac{\varphi_{1,\bmw}(\rho_{i,\bmw})}{2 \rho_{i,\bmw}} - \frac{\varphi_{1,\bmw}'(\rho_{i,\bmw})}{2} \bigg]
  |\bmu_i^\CT \bhtu_{i,\bmw}|^2
  \frac{\varphi_{2,\bmw}(\rho_{i,\bmw})}{\varphi_{1,\bmw}(\rho_{i,\bmw})}
  &
  &\quad&
  &\asunifto{\bmw \in \clW_>(\nu)}
  &\quad&
  0
  , \\
  \label{p:eq:singvec:sum:parts:lim:chi2}
  \chi_2(\bmw)
  &&
  -
  &&
  &
  \theta_i^2
  \bigg[ -\frac{\varphi_{2,\bmw}(\rho_{i,\bmw})}{2 \rho_{i,\bmw}} - \frac{\varphi_{2,\bmw}'(\rho_{i,\bmw})}{2} \bigg]
  |\bmu_i^\CT \bhtu_{i,\bmw}|^2
  &
  &&
  &\asunifto{\bmw \in \clW_>(\nu)}
  &&
  0
  , \\
  \label{p:eq:singvec:sum:parts:lim:chi3}
  \chi_3(\bmw)
  &&
  &&
  &&
  &&
  &\asunifto{\bmw \in \clW_>(\nu)}
  &&
  0
  .
\end{alignat}
\end{subequations}
Finally,
applying \cref{p:eq:singvec:sum:parts:lim:all}
to \cref{p:eq:singvec:sum}
yields that (see \insupp{calc:unif:conv:after:lemma} for detailed calculations),
for any $\nu > 0$,
\begin{equation} \label{eq:singvec:lim:above}
  r_i(\bmw,\bmY)
  =
  |\bmu_i^\CT \bhtu_{i,\bmw}|^2
  \asunifto{\bmw \in \clW_>(\nu)}
  - \frac{2 \varphi_{1,\bmw}(\rho_{i,\bmw})}{\theta_i^2 D_{\bmw}'(\rho_{i,\bmw})}
  .
\end{equation}

\subsubsection{Performance below the phase transition}%
This section bounds the limit of
the performance $r_i(\bmw,\bmY)$
below the phase transition.
Namely, for any $\nu > 0$,
we derive a uniform bound with respect to $\bmw$
over the domain
$\clW_\leq(\nu)
\coloneqq \Delta_L \setminus \clW_>(\nu)
=
\{\bmw \in \Delta_L : \brtheta_{i,\bmw} \leq b_{\bmw} + \nu\}$.

Following the approach of \cite{bloemendal2014ill},
we study $r_i(\bmw,\bmY) = |\bmu_i^\CT \bhtu_{i,\bmw}|^2$
by first obtaining the following deterministic bound
(see \insupp{calc:below:phase:det:bound} for detailed calculations):
\begin{equation} \label{bdd_spectal}
  |\bmu_i^\CT \bhtu_{i,\bmw}|^2
  \leq
  -\nu \cdot
  \Im\Big\{
    \zeta_{i,\bmw}^{-1}
    \,
    \big[
      \btlM(\bmw,\zeta_{i,\bmw})
      - \btlM(\bmw,\zeta_{i,\bmw}) \left[\bmM(\bmw,\zeta_{i,\bmw})\right]^{-1} \btlM(\bmw,\zeta_{i,\bmw})
    \big]_{ii}
  \Big\}
\end{equation}
where $\zeta_{i,\bmw}^2 \coloneqq \httheta_{i,\bmw}^2 + \imath \nu$
and
\begin{equation}
  \btlM(\bmw,\zeta)
  \coloneqq
  \begin{bmatrix} \bmU & \\ & \btlQ_{\bmw} \end{bmatrix}^\CT
  \bmG(\bmw,\zeta)
  \begin{bmatrix} \bmU & \\ & \btlQ_{\bmw} \end{bmatrix}
  .
\end{equation}

Next, by standard calculations (see \insupp{calc:below:phase:bounds}),
we have that almost surely, eventually,
the following bounds hold for all $\bmw \in \clW_{\leq}(\nu)$:
\begin{align} \label{eq:singvec:below:bound}
  |\zeta_{i,\bmw}^{-1}| &\leq \tlC_1
  , &
  \big\| \btlM(\bmw,\zeta_{i,\bmw}) \big\|_{\opnorm} &\leq \tlC_2
  , &
  \|\bmM(\bmw,\zeta_{i,\bmw})^{-1}\|_{\opnorm} &\leq \tlC_3 \nu^{-1/2}
  ,
\end{align}
where $\tlC_1$, $\tlC_2$, and $\tlC_3$ do not depend on $\nu$ or $\bmw$.

Finally,
applying \cref{eq:singvec:below:bound}
to \cref{bdd_spectal}
yields that for any $\nu > 0$,
almost surely, eventually,
\begin{equation} \label{eq:singvec:bound:below}
  \sup_{\bmw \in \clW_{\leq}(\nu)}
  r_i(\bmw,\bmY)
  =
  \sup_{\bmw \in \clW_{\leq}(\nu)}
  |\bmu_i^\CT \bhtu_{i,\bmw}|^2
  \leq
  \tlC_4 (\nu + \nu^{1/2})
  ,
\end{equation}
where $\tlC_4 \coloneqq \tlC_1 \max(\tlC_2,\tlC_2^2\tlC_3)$
does not depend on $\nu$.

\subsubsection{Uniform convergence and algebraic form of performance}
\label{proof:rec:lim:unif:alg}%
Noting that $\nu > 0$ can be arbitrarily small
in \cref{eq:singvec:lim:above,eq:singvec:bound:below}
yields uniform convergence across $\bmw \in \Delta_L$,
i.e.,
\begin{equation} \label{eq:rec:lim:nonalg}
  r_i(\bmw,\bmY)
  \quad
  \asunifto{\bmw \in \Delta_L}
  \quad
  \brr_i(\bmw)
  \coloneqq
  \begin{cases}
    - \frac{2 \varphi_{1,\bmw}(\rho_{i,\bmw})}{\theta_i^2 D_{\bmw}'(\rho_{i,\bmw})}
    , & \text{if } \theta_i > \tltheta_{\bmw}
    , \\
    0
    , & \text{otherwise}
    .
  \end{cases}
\end{equation}
Since $r_i(\bmw,\bmY)$ is scale invariant,
i.e., $r_i(\gamma\bmw,\bmY) = r_i(\bmw,\bmY)$ for any $\gamma > 0$,
it immediately follows that the convergence
is also uniform over $\bbR_{\geq 0}^L \setminus \{\bm0_L\}$ as well.

The proof concludes by deriving the algebraic description \cref{eq:rec:lim}
of $\brr_i(\bmw)$.
Following the approach of \cite[Section 5.2]{hong2018apo},
we change variables to
\begin{equation}
  \psi_{\bmw}(\zeta) \coloneqq \frac{c\zeta}{\varphi_{1,\bmw}(\zeta)}
  ,
\end{equation}
and observe that,
analogously to \cite[Section 5.3]{hong2018apo},
\begin{align}
  \label{p:eq:Dpsi}
  D_{\bmw}(\zeta)
  &
  = \frac{1-B_{i,\bmw}(\psi_{\bmw}(\zeta))}{\theta_i^2}
  , &
  \frac{D_{\bmw}'(\zeta)}{\zeta}
  &
  =
  -\frac{2c}{\theta_i^2}
  \frac{B_{i,\bmw}'(\psi_{\bmw}(\zeta))}{A_{\bmw}(\psi_{\bmw}(\zeta))}
  ,
\end{align}
$\psi_{\bmw}(b_{\bmw}^+) = \alpha_{\bmw}$
and $\psi_{\bmw}(\rho_{i,\bmw}) = \beta_{i,\bmw}$ when $\theta_i^2 > \tltheta_{\bmw}^2$,
where $\alpha_{\bmw}$ and $\beta_{i,\bmw}$ are
the largest real roots of $A_{\bmw}$ and $B_{i,\bmw}$, respectively.
See \insupp{calc:alg:desc:rewrite} for the detailed derivations.

Even though $\psi_{\bmw}(\rho_i)$ is defined only when $\theta_i^2 > \tltheta_{\bmw}^2$,
the largest real roots $\alpha_{\bmw}$ and $\beta_{i,\bmw}$ are always defined
and always larger than $\max_\ell(w_\ell v_\ell)$.
Thus
\begin{equation} \label{p:eq:phasetrans}
\theta_i^2 > \tltheta_{\bmw}^2
= \frac{\theta_i^2}{1-B_{i,\bmw}(\psi_{\bmw}(b_{\bmw}^+))}
\; \Leftrightarrow \;
B_{i,\bmw}(\alpha_{\bmw}) < 0
\Leftrightarrow \;
\alpha_{\bmw} < \beta_{i,\bmw}
\; \Leftrightarrow \;
A_{\bmw}(\beta_{i,\bmw}) > 0
,
\end{equation}
where the final two equivalences hold because
$A_{\bmw}(x)$ and $B_{i,\bmw}(x)$ are both strictly increasing functions
for $x > \max_\ell(w_\ell v_\ell)$
and $A_{\bmw}(\alpha_{\bmw}) = B_{i,\bmw}(\beta_{i,\bmw}) = 0$.

Finally,
using \cref{p:eq:Dpsi,p:eq:phasetrans}
to rewrite \cref{eq:rec:lim:nonalg}
concludes the proof
of \cref{lem:rec:lim}.

\subsection{Optimization of the almost sure limit}
\label{sct:rec:lim:opt}%
This section optimizes $\brr_i(\bmw)$
and shows that the maximizer is unique
(up to scaling).

\begin{lemma}[Optimization of the almost sure limit]%
  \label{lem:rec:lim:opt}%
  The asymptotic performance $\brr_i(\bmw)$
  is continuous
  and is maximized as:
  \begin{align}
    \{\gamma \bbrw^{\opt}_i : \gamma > 0\}
    &=
    \argmax_{\bmw \in \bbR_{\geq 0}^L \setminus \{\bm0_L\}}
    \; \brr_i(\bmw)
    , &
    \brr^{\opt}_i
    &=
    \max_{\bmw \in \bbR_{\geq 0}^L \setminus \{\bm0_L\}}
    \; \brr_i(\bmw)
    ,
  \end{align}
  except when
  $\sum_{\ell=1}^L c_\ell (\lambda_i /v_\ell)^2 \leq 1$,
  in which case
  $\brr_i(\bmw) = 0$ for all weights $\bmw \in \bbR_{\geq 0}^L \setminus \{\bm0_L\}$.
\end{lemma}

The remainder of this subsection proves \cref{lem:rec:lim:opt}.
A major challenge for the proof is that
$\brr_i(\bmw)$ is defined implicitly via the root $\beta_{i,\bmw}$,
and setting the gradient equal to zero
to obtain local maxima
yields a complicated nonlinear system of equations to solve.
Surprisingly, the system turns out to have a simple solution
that we derive by carefully exploiting the structure of the system.
Moreover, we show that the solution is globally optimal,
obtaining the optimal weights and their corresponding performance.

Before deriving the gradient,
note that
$\brr_i(\bmw)$
is not always differentiable everywhere
due to its truncation at zero.
However, it can be rewritten as
\begin{equation} \label{eq:rec:lim:untruncated}
  \brr_i(\bmw) = \max\big(0,\tlr_i(\bmw)\big)
  \quad \text{where} \quad
  \tlr_i(\bmw)
  \coloneqq
  \frac{1}{\beta_{i,\bmw}} \frac{A_{\bmw}(\beta_{i,\bmw})}{B_{i,\bmw}'(\beta_{i,\bmw})}
  ,
\end{equation}
so the problem reduces to maximizing
the differentiable function $\tlr_i(\bmw)$
then checking whether it is positive.
Furthermore, $\tlr_i(\bmw)$ achieves
its maximum over the feasible region $\bbR_{\geq 0}^L \setminus \{\bm0_L\}$.
To see why,
note that $\tlr_i(\bmw)$ is scale-invariant,
i.e.,
\begin{equation*}
  \forall_{\gamma > 0} \quad
  \tlr_i(\gamma \bmw)
  =
  \frac{1}{\beta_{i,\gamma \bmw}}
  \frac{A_{\gamma \bmw}(\beta_{i,\gamma \bmw})}{B_{i,\gamma \bmw}'(\beta_{i,\gamma \bmw})}
  =
  \frac{1}{\gamma \beta_{i,\bmw}}
  \frac{A_{\bmw}(\beta_{i,\bmw})}{(1/\gamma)B_{i,\bmw}'(\beta_{i,\bmw})}
  =
  \tlr_i(\bmw)
\end{equation*}
since
$A_{\gamma \bmw}(x) = A_{\bmw}(x/\gamma)$,
$B_{i,\gamma \bmw}(x) = B_{i,\bmw}(x/\gamma)$,
and
$B_{i,\gamma \bmw}'(x) = (1/\gamma) B_{i,\bmw}'(x/\gamma)$,
resulting additionally in
$\beta_{i,\gamma\bmw} = \gamma \beta_{i,\bmw}$.
Thus,
$\tlr_i(\bmw)$ can equivalently be maximized over the compact set
$\Delta_L \coloneqq \{\bmw \in \bbR_{\geq 0}^L : w_1 + \cdots + w_L = 1\}$
and hence achieves its maximum.

Next, note that
the feasible region $\bbR_{\geq 0}^L \setminus \{\bm0_L\}$ is not open,
so we partition it into $2^L-1$ sets
according to which weights are zero.
Namely, consider partitions of the form
\begin{equation*}
  \clP_{\clL}
  \coloneqq
  \{\bmw \in \bbR_{\geq 0}^L :
    \forall_{\ell \in    \clL} \ w_\ell = 0,
    \forall_{\ell \notin \clL} \ w_\ell > 0
  \}
  ,
  \quad
  \text{for $\clL \subset \{1,\dots,L\}$ a proper subset}
  .
\end{equation*}
Since $\tlr_i(\bmw)$ achieves its maximum,
a maximizer exists within at least one of the partitions.
Moreover, since $\tlr_i(\bmw)$ is differentiable,
$\tlr_i(\bmw)$ is maximized at a critical point
of a partition.
It remains to identify and compare
the critical points of all the partitions $\clP_{\clL}$.

First consider $\clP_\emptyset$,
i.e., the set of positive weights $w_1,\dots,w_L > 0$,
and let $\tlw_j \coloneqq 1/w_j$.
This parameterization ends up simplifying the manipulations.
Differentiating $\tlr_i(\bmw)$ and \cref{eq:rec:lim:AB}
with respect to $\tlw_j$ yields
\begin{subequations}
\label{eq:opt:weight:system}
\begin{align}
  \label{eq:rdiff}
  \diff{\tlr_i(\bmw)}{\tlw_j}
  =
  \tlr_i(\bmw)
  &
  \bigg[
    -\frac{1}{            \beta_{i,\bmw} } \diff{            \beta_{i,\bmw} }{\tlw_j}
    +\frac{1}{   A_{\bmw}(\beta_{i,\bmw})} \diff{   A_{\bmw}(\beta_{i,\bmw})}{\tlw_j}
    -\frac{1}{B_{i,\bmw}'(\beta_{i,\bmw})} \diff{B_{i,\bmw}'(\beta_{i,\bmw})}{\tlw_j}
  \bigg]
  , \\
  \label{eq:Abdiff}
  \diff{A_{\bmw}(\beta_{i,\bmw})}{\tlw_j}
  &=
  A_{\bmw}'(\beta_{i,\bmw}) \diff{\beta_{i,\bmw}}{\tlw_j}
  + 2 \frac{c_j v_j^2}{(\beta_{i,\bmw}\tlw_j - v_j)^3} \beta_{i,\bmw}
  , \\
  \label{eq:Bpbdiff}
  \diff{B_{i,\bmw}'(\beta_{i,\bmw})}{\tlw_j}
  &=
  B_{i,\bmw}''(\beta_{i,\bmw}) \diff{\beta_{i,\bmw}}{\tlw_j}
  - 2 \lambda_i \frac{c_j \tlw_j}{(\beta_{i,\bmw}\tlw_j - v_j)^3} \beta_{i,\bmw}
  +   \lambda_i \frac{c_j       }{(\beta_{i,\bmw}\tlw_j - v_j)^2}
  , \\
  \label{eq:Bbdiff}
  0 &= \diff{B_{i,\bmw}(\beta_{i,\bmw})}{\tlw_j}
  = B_{i,\bmw}'(\beta_{i,\bmw}) \diff{\beta_{i,\bmw}}{\tlw_j}
  + \lambda_i \frac{c_j}{(\beta_{i,\bmw}\tlw_j - v_j)^2} \beta_{i,\bmw}
  ,
\end{align}
\end{subequations}
where one must carefully account for the fact that
$A_{\bmw}(x)$, $B_{i,\bmw}(x)$ and $\beta_{i,\bmw}$ are functions of $\tlw_j$.
The problem now is to set $\idiff{\tlr_i(\bmw)}{\tlw_j}$ equal to zero
and carefully exploit the structure of the system
\cref{eq:opt:weight:system}
to solve it.

In particular,
rewriting \cref{eq:Abdiff,eq:Bpbdiff}
in terms of $\idiff{\beta_{i,\bmw}}{\tlw_j}$
using \cref{eq:Bbdiff} yields
\begin{subequations}
\label{eq:Abdiff1:Bpbdiff1}
\begin{align}
  \label{eq:Abdiff1}
  \diff{A_{\bmw}(\beta_{i,\bmw})}{\tlw_j}
  &=
  \bigg[
    A_{\bmw}'(\beta_{i,\bmw})
    -
    \frac{2 B_{i,\bmw}'(\beta_{i,\bmw})}{\lambda_i}
    \frac{v_j^2}{\beta_{i,\bmw} \tlw_j - v_j}
  \bigg]
  \diff{\beta_{i,\bmw}}{\tlw_j}
  , \\
  \label{eq:Bpbdiff1}
  \diff{B_{i,\bmw}'(\beta_{i,\bmw})}{\tlw_j}
  &=
  \bigg[
    B_{i,\bmw}''(\beta_{i,\bmw})
    +
    2 B_{i,\bmw}'(\beta_{i,\bmw}) \frac{\tlw_j}{\beta_{i,\bmw} \tlw_j - v_j}
  \bigg]
  \diff{\beta_{i,\bmw}}{\tlw_j}
  -
  \frac{B_{i,\bmw}'(\beta_{i,\bmw})}{\beta_{i,\bmw}} \diff{\beta_{i,\bmw}}{\tlw_j}
  .
\end{align}
\end{subequations}
Substituting \cref{eq:Abdiff1,eq:Bpbdiff1}
into \cref{eq:rdiff} then rearranging yields
\begin{equation} \label{eq:rdiff1}
  \diff{\tlr_i(\bmw)}{\tlw_j}
  =
  \frac{2}{\lambda_i \beta_{i,\bmw}}
  \diff{\beta_{i,\bmw}}{\tlw_j}
  \bigg[
    \lambda_i \Delta_{i,\bmw}
    -
    \frac{
      \lambda_i \beta_{i,\bmw} \tlr_i(\bmw) \tlw_j + v_j^2
    }{
      \beta_{i,\bmw} \tlw_j - v_j
    }
  \bigg]
  ,
\end{equation}
where the following term
is independent of $j$:
\begin{equation*}
  \Delta_{i,\bmw}
  \coloneqq
  \frac{1}{2}
  \frac{A_{\bmw}(\beta_{i,\bmw})}{B_{i,\bmw}'(\beta_{i,\bmw})}
  \bigg[
    \frac{A_{\bmw}'(\beta_{i,\bmw})}{A_{\bmw}(\beta_{i,\bmw})}
    -
    \frac{B_{i,\bmw}''(\beta_{i,\bmw})}{B_{i,\bmw}'(\beta_{i,\bmw})}
  \bigg]
  .
\end{equation*}
Since $\beta_{i,\bmw} > \max_\ell(w_\ell v_\ell) > 0$,
it follows from \cref{eq:Bbdiff} that $\idiff{\beta_{i,\bmw}}{\tlw_j} \neq 0$,
so it follows from \cref{eq:rdiff1} that $\idiff{\tlr_i(\bmw)}{\tlw_j}$ is zero
exactly when
\begin{equation} \label{eq:rsopt1}
  \lambda_i \Delta_{i,\bmw}
  =
  \frac{
    \lambda_i \beta_{i,\bmw} \tlr_i(\bmw) \tlw_j + v_j^2
  }{
    \beta_{i,\bmw} \tlw_j - v_j
  }
  .
\end{equation}
Rearranging \cref{eq:rec:lim:untruncated} and substituting \cref{eq:rsopt1} yields
\begin{align*}
  0 &=
  A_{\bmw}(\beta_{i,\bmw}) - \tlr_i(\bmw) \beta_{i,\bmw} B_{i,\bmw}'(\beta_{i,\bmw})
  = 1 -
  \sum_{\ell=1}^L
  \frac{
    c_\ell (v_\ell^2 + \lambda_i \beta_{i,\bmw} \tlr_i(\bmw) \tlw_\ell)
  }{
    (\beta_{i,\bmw} \tlw_\ell - v_\ell)^2
  }
  \\&
  = 1 -
  \Delta_{i,\bmw} \lambda_i
  \sum_{\ell=1}^L \frac{c_\ell}{\beta_{i,\bmw} \tlw_\ell - v_\ell}
  =
  1 - \Delta_{i,\bmw}(1-B_{i,\bmw}(\beta_{i,\bmw}))
  = 1 - \Delta_{i,\bmw}
  ,
\end{align*}
so $\Delta_{i,\bmw} = 1$.
Substituting into \cref{eq:rsopt1}
and solving for $\tlw_j$ yields
\begin{equation} \label{eq:design:weights:scaling}
  w_j
  = \frac{1}{\tlw_j}
  = \frac{(1 - \tlr_i(\bmw)) \beta_{i,\bmw}}{v_j (1 + v_j/\lambda_i)}
  = \frac{\gamma_{i,\bmw}}{v_j (1 + v_j/\lambda_i)}
  ,
\end{equation}
where the constant
$\gamma_{i,\bmw} \coloneqq (1 - \tlr_i(\bmw)) \beta_{i,\bmw}$
is: a) independent of $j$,
b) parameterizes the ray of critical points in $\clP_\emptyset$,
and c) can be chosen freely,
e.g., as unity yielding $\bbrw^{\opt}_i$
from \cref{eq:opt:weight:asymp}.

Solving \cref{eq:design:weights:scaling} for $\beta_{i,\bmw}$,
substituting into $0 = B_{i,\bmw}(\beta_{i,\bmw})$,
and rearranging yields that the corresponding $\tlr_i(\bbrw^{\opt}_i)$
is a root of
\begin{equation*}
  R_{i,\emptyset}(x)
  \coloneqq
  1 -
  \sum_{\ell=1}^L
    \frac{c_\ell}{v_\ell/\lambda_i}
    \,\frac{1-x}{v_\ell/\lambda_i+x}
  .
\end{equation*}
Since $R_{i,\emptyset}(x)$
increases from negative infinity to one
as $x$ increases from $-\min_\ell(v_\ell)/\lambda_i$ to one,
it has exactly one real root in that domain.
In particular, this root is the largest real root
since $R_{i,\emptyset}(x) \geq 1$ for $x \geq 1$.
Furthermore, $\tlr_i(\bbrw^{\opt}_i)$
increases continuously to one
as $c \coloneqq c_1 + \cdots + c_L$ increases to infinity,
so $\tlr_i(\bbrw^{\opt}_i)$ must be the largest real root.

Likewise,
the critical points of other partitions $\clP_{\clL}$
are given by setting
the positive weights proportional to $\bbrw^{\opt}_i$
with the corresponding $\tlr_i(\bbrw^{\opt}_i)$
given by the largest real root of
\begin{equation*}
  R_{i,\clL}(x)
  \coloneqq
  1 -
  \sum_{\ell \notin \clL}
    \frac{c_\ell}{v_\ell/\lambda_i}
    \,\frac{1-x}{v_\ell/\lambda_i+x}
  .
\end{equation*}
For $\clL_1 \subset \clL_2$ a proper subset,
the largest real root of $R_{i,\clL_1}$
is greater than that of $R_{i,\clL_2}$
since
$R_{i,\clL_1}(x) < R_{i,\clL_2}(x)$
for any $x \in (-\min_\ell(v_\ell)/\lambda_i,1)$.
As a result, $\tlr_i(\bmw)$ is maximized in $\clP_\emptyset$.

Finally, we check when $\tlr_i(\bbrw^{\opt}_i)$ is positive.
Recalling that $R_{i,\emptyset}(x)$ is an increasing function
on $x \in (0,1)$
and noting that
$R_{i,\emptyset}(0) = 1 - \sum_{\ell=1}^L c_\ell (\lambda_i/v_\ell)^2$
yields that $\tlr_i(\bbrw^{\opt}_i)$ is positive
if and only if $\sum_{\ell=1}^L c_\ell (\lambda_i/v_\ell)^2 > 1$.
When it is positive, the maximizers are given by the critical points above;
otherwise $\brr_i(\bmw) = 0$ for all $\bmw$.
This concludes the proof of \cref{lem:rec:lim:opt}.

%% file: content/ext,signal,var,est.tex
%!TEX root = ../wpca.tex

\section{When signal and noise variances are unknown}
\label{sec:signal:var:est}%
The asymptotic optimal weights
from the main result (\cref{thm:opt:weight:asymp})
depend on both the signal component variance $\lambda_i$
and the noise variances $\bmv$,
but one or both of these variances may be unknown in some settings.
Of course,
one could estimate these variances
using existing ideas and approaches,
then plug them into \cref{eq:opt:weight:asymp}
in \cref{thm:opt:weight:asymp}.
The question then is,
are these resulting estimated weights also asymptotically optimal?
Fortunately, it is straightforward to see
that the answer is yes
under natural conditions on the variance estimates.
The following \lcnamecref{prop:est:var:consistency}
makes this statement precise;
it follows straightforwardly from \cref{lem:rec:lim,lem:rec:lim:opt}.

\begin{proposition}[Asymptotic optimality with estimated variances]%
  \label{prop:est:var:consistency}%
  Suppose
  $\htlambda_i(\bmY)$ and $\bhtv(\bmY)$
  are consistent estimates of
  $\lambda_i$ and $\bmv$,
  i.e.,
  \begin{align}
    \label{eq:est:var:consistency}
    \htlambda_i(\bmY) &\asto \lambda_i
    , &
    \bhtv(\bmY) &\asto \bmv
    .
  \end{align}
  Let $\bhtw^{\opt}_i(\bmY)$
  be the estimated weights
  obtained by plugging $\bhtv(\bmY)$ and $\htlambda_i(\bmY)$
  into \cref{eq:opt:weight:asymp}
  in \cref{thm:opt:weight:asymp},
  i.e.,
  \begin{equation}
    \label{eq:est:weight}
    \bhtw^{\opt}_i(\bmY)
    \coloneqq
    \bigg(
      \frac{1}{\htv_1(\bmY)}\frac{1}{1 + \htv_1(\bmY)/\htlambda_i(\bmY)},
      \dots,
      \frac{1}{\htv_L(\bmY)}\frac{1}{1 + \htv_L(\bmY)/\htlambda_i(\bmY)}
    \bigg)
    .
  \end{equation}
  Then
  the estimated weights and their corresponding performance converge
  to the asymptotic optimal weights and performance, i.e.,
  \begin{align}
    \bhtw^{\opt}_i(\bmY)
    &\asto
    \bbrw^{\opt}_i
    , &
    r_i(\bhtw^{\opt}_i(\bmY),\bmY)
    &\asto
    \brr^{\opt}_i
    .
  \end{align}
  Namely, the estimated weights are asymptotically optimal.
\end{proposition}

Thus,
optimal weighting
only needs consistent estimates of $\lambda_i$ and $\bmv$
that may be obtained using any one of various existing approaches;
which one is most appropriate will depend on the specific application.
Here we consider a simple pair of estimators as an illustrative example.

\begin{example}[Variance estimators]%
  \label{ex:var:est}%
  As an illustrative example,
  consider the following simple estimators
  for the signal and noise variances
  \begin{align}
    \label{eq:ex:var:est}
    \htlambda_i(\bmY)
    &\coloneqq
    \Xi\Bigg(
      \htlambda^{\inv}_i\big(\bmY; \bhtv(\bmY)\big);
      \Bigg[
        \sum_{\ell=1}^L \frac{p_\ell}{\htv_\ell(\bmY)}
      \Bigg]^{-1}
    \Bigg)
    , &
    \bhtv(\bmY)
    &\coloneqq \bigg(
      \frac{\|\bmY_1\|_{\frob}^2}{d \, n_1},
      \dots,
      \frac{\|\bmY_L\|_{\frob}^2}{d \, n_L}
    \bigg)
    ,
  \end{align}
  where
  \begin{align}
    \htlambda^{\inv}_i(\bmY; \bmv)
    &\coloneqq
    \text{\ith leading eigenvalue of }
    \sum_{\ell = 1}^L
    \bigg(\frac{1/v_\ell}{n_1 / v_1 + \cdots + n_L / v_L} \bigg) \,
    \bmY_\ell \bmY_\ell^\CT
    , \\
    \Xi(\lambda; v)
    &\coloneqq
    \text{the larger root of the quadratic polynomial }
    (x + v/c)(x + v) - \lambda x
    .
  \end{align}
  It is straightforward to verify
  with standard techniques
  (see \insupp{sec:supp:ext:signal:var} for details)
  that
  these estimators are consistent
  as long as
  $c > ( \brv / \lambda_i)^2$,
  i.e., when inverse noise variance weighting is above the phase transition.
  Thus, by \cref{prop:est:var:consistency},
  the resulting estimated weights and their corresponding performance converge
  to the asymptotic optimal weights and performance.

  \Cref{fig:est:weight:sim} illustrates the nonasymptotic behavior
  of these estimated weights
  in numerical simulations.
  The data is generated as in \cref{sec:opt:weight:sim},
  with dimensionality $d = 1000$
  and block sizes $\bmn = (4000,8000)$.
  The estimated weights and their performance
  generally concentrate around
  the asymptotic optimal weights and performance.
  Moreover, the estimated weights achieve performance closely matching
  the nonasymptotic empirically optimized weights
  in \cref{fig:opt:weight:sim:1000}.
\end{example}

\begin{figure}[t] \centering
  \begin{tikzpicture}
    \tikzmath{
      \vrat=3; \lmb=1;
      \wlim=(1*(1+\lmb))/(\vrat*(\vrat+\lmb));
      \rlim=0.6782867492531689;
    };
    \begin{axis}[
      xlabel={Estimated relative weight $\htw^{\opt}_{1,2}(\bmY)/\htw^{\opt}_{1,1}(\bmY)$},
      xmin=0, xmax=0.55,
      xtick={0,0.16667,0.33333,0.5},
      xticklabels={0,1/6,1/3,1/2},
      ymin=0, ymax=100,
      ytick={0,20,...,100},
      width=0.5\linewidth, height=0.25\linewidth,
      font=\footnotesize,
      clip=false,
    ]
      \addplot[ybar interval,fill opacity=0.5,fill=RoyalBlue]
        table {figs/est,var,sim,weights,c-4-8,v-1-3,lambda-1,d-1000,ntrials-500,ratiores-0.01.dat};
      \addplot[ultra thick,Dark2-B,dashed] coordinates {(\wlim,0) (\wlim,104)}
        node[anchor=south] {\shortstack{asymptotic optimal\\relative weight $\brw^{\opt}_{1,2}/\brw^{\opt}_{1,1}$}};
    \end{axis}
  \end{tikzpicture}
  \qquad
  \begin{tikzpicture}
    \tikzmath{
      \vrat=3; \lmb=1;
      \wlim=(1*(1+\lmb))/(\vrat*(\vrat+\lmb));
      \rlim=0.6782867492531689;
    };
    \begin{axis}[
      xlabel={Performance of estimated weights $r_1(\bhtw^{\opt}_1(\bmY),\bmY)$},
      xmin=0, xmax=1,
      xtick={0,0.25,0.5,0.75,1},
      xticklabels={0,1/4,1/2,3/4,1},
      ymin=0, ymax=25,
      ytick={0,5,...,25},
      width=0.5\linewidth, height=0.25\linewidth,
      font=\footnotesize,
      clip=false,
    ]
      \addplot[ybar interval,fill opacity=0.5,fill=RoyalBlue]
        table {figs/est,var,sim,rec,c-4-8,v-1-3,lambda-1,d-1000,ntrials-500,ratiores-0.01.dat};
      \addplot[ultra thick,Dark2-B,dashed] coordinates {(\rlim,0) (\rlim,26)}
        node[anchor=south] {\shortstack{asymptotic optimal\\performance $\brr^{\opt}_1$}};
    \end{axis}
  \end{tikzpicture}%
  \caption{Nonasymptotic empirical distributions of estimated weights $\bhtw^{\opt}_1(\bmY)$
    from \cref{eq:est:weight}
    and corresponding performance $r_1(\bhtw^{\opt}_1(\bmY),\bmY)$
    for an illustrative example with two blocks of data
    $\bmY_1 \in \bbR^{d \times n_1}$
    and
    $\bmY_2 \in \bbR^{d \times n_2}$,
    where the estimated weights are computed
    using the signal and noise variance estimators
    \cref{eq:ex:var:est}
    from \cref{ex:var:est}.
    The data blocks are generated
    with noise variances $v_1 = 1$ and $v_2 = 3$,
    one underlying component having variance $\lambda_1 = 1$,
    and dimensions $d = 1000$, $n_1 = 4000$ and $n_2 = 8000$.}
  \label{fig:est:weight:sim}
\end{figure}

%% file: content/experiment.tex
%!TEX root = ../wpca.tex

\section{Illustration on real data from astronomy}
\label{sec:exp:astro}%
This section illustrates optimally weighted PCA
on real data from astronomy.
In particular, we consider quasar spectra
from the Sloan Digital Sky Survey (SDSS) Data Release 16 \cite{ahumada2020t1d}
using the associated DR16Q quasar catalog \cite{lyke2020tsd}.
Each spectrum is a vector of flux measurements across wavelengths
for a particular quasar,
and comes with associated noise variance estimates
across wavelengths.
The noise is heteroscedastic across both quasars and wavelengths,
but here we focus on a subset
that is somewhat homoscedastic across wavelengths.
\Insupp{exp:preprocess}
describes the details of
the subset selected
and the preprocessing performed
(filtering, interpolation, centering, normalization).

\begin{figure} \centering
  \includegraphics[width=\linewidth]{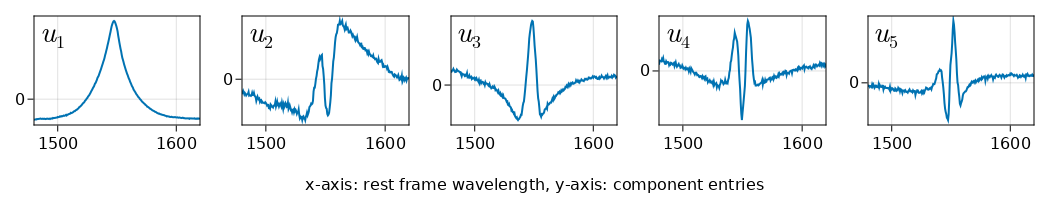}
  \caption{Ground truth components computed from 5000 samples with smallest variances.}
  \label{fig:exp:ground}
\end{figure}

The resulting data has $d=281$ wavelengths measured for $n^{\full}=10459$ spectra,
yielding a data matrix
$\bmY^{\full} \in \bbR^{d \times n^{\full}}$
with a vector of associated variances
$\bmv^{\full} \in \bbR_{\geq 0}^{n^{\full}}$.
\Cref{fig:exp:ground} shows plots of components $\bmu_1,\dots,\bmu_5 \in \bbR^d$
computed via an unweighted PCA on the 5000 samples
from $\bmY^{\full}$
with smallest variances.
We regard these components as ``ground truth.''

\begin{figure} \centering
  \subfloat[Data matrix $\bmY \in \bbR^{d \times n}$ (top)
    and associated noise variances $\bmv \in \bbR_{\geq 0}^n$ (bottom).
    \label{fig:exp:test:matrix}]{%
    \includegraphics[width=0.5\linewidth]{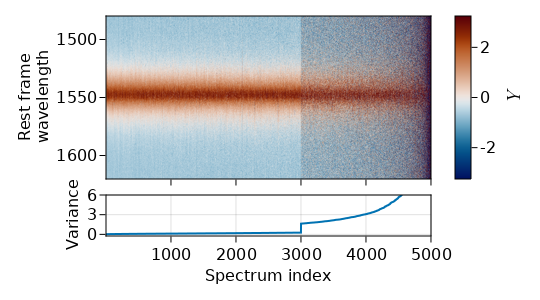}%
  }
  \hfill
  \subfloat[Example spectra
    with associated high, medium, and low noise variances.
    \label{fig:exp:test:example:spectra}]{%
    \includegraphics[width=0.45\linewidth]{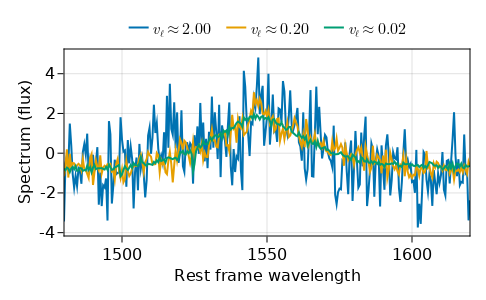}
  }
  \caption{Quasar spectra dataset and example spectra.}
\end{figure}

To evaluate the various weighting schemes,
we formed a test dataset $\bmY \in \bbR^{d \times n}$
containing $n=5000$ samples
by combining
the 3000 samples with the smallest variances
(a subset of the 5000 samples used to produce ground truth)
and the 2000 samples with the largest variances.
This provides a dataset with noise heteroscedasticity across samples,
shown as a heatmap in \cref{fig:exp:test:matrix}.
\Cref{fig:exp:test:example:spectra} shows a few sample spectra
from the dataset;
note that they have a common shape but have varying levels of noise.

We computed the leading singular vectors $\bhtu_1(\bmw\revone{_1},\bmY),\dots,\bhtu_5(\bmw\revone{_5},\bmY)$
via unweighted PCA,
inverse variance weighted PCA,
and optimally weighted PCA\revone{;
the optimally weighted PCA
has component-specific weights,
so its weights are not the same across the components.%
\footnote{\revone{While the components
$\bhtu_1(\bhtw^{\opt}_1(\bmY),\bmY),\dots,\bhtu_5(\bhtw^{\opt}_1(\bmY),\bmY)$
from optimally weighted PCA
are not guaranteed to be orthogonal in general
(due to the component-specific weighting),
they were approximately orthogonal here.
In particular,
$\max_{i \neq j}
 |\bhtu_i(\bhtw^{\opt}_i(\bmY),\bmY)^\CT \bhtu_j(\bhtw^{\opt}_j(\bmY),\bmY)|^2
 \approx 0.00013$.}}
}%
We used the provided noise variances,
and estimated the signal variances
using the estimator from \cref{ex:var:est}
(with the provided noise variances substituted).
\Cref{tbl:exp:astro} shows the recovery of the ground truth singular vectors
calculated from the ``clean'' samples above.

\begin{table}[h] \centering
  \caption{Recoveries $r_i(\bmw,\bmY)$
    for unweighted PCA ($w_\ell = 1$),
    inverse variance weighted PCA ($w_\ell = 1/v_\ell$),
    and optimally weighted PCA.
    Higher is better, best value (up to rounding) is shown in bold.}
  \label{tbl:exp:astro}
  \begin{tabular}{|l|ccccc|}
    \hline
    Component                      &            1   &            2   &            3   &            4   &            5   \\ \hline
    Unweighted PCA                 &  $    0.003 $  &  $    0.307 $  &  $    0.009 $  &  $    0.004 $  &  $    0.018 $  \\
    Inverse variance weighted PCA  &  $\bm{1.000}$  &  $    0.903 $  &  $    0.915 $  &  $    0.817 $  &  $    0.811 $  \\
    Optimally weighted PCA         &  $\bm{1.000}$  &  $\bm{0.920}$  &  $\bm{0.934}$  &  $\bm{0.884}$  &  $\bm{0.880}$  \\ \hline
  \end{tabular}
\end{table}

The following observations
apply to this data
and also summarize some of the main themes
of this paper:
\begin{itemize}
  \item unweighted PCA performs poorly for heteroscedastic data,
  \item inverse variance weighted PCA performs much better than unweighted PCA,
  and
  \item optimally weighted PCA performs even better than inverse variance weighted PCA\@.
\end{itemize}
Similar comparisons for these leading components
occurred for many of the other test datasets we tried.
The comparisons were less consistent
for components 6 and on;
optimal weights were sometimes better
and inverse noise variance weights were sometimes better.
This inconsistent behavior
is potentially due to
the data not matching the model closely enough
or perhaps the ground truth needing to be chosen differently.
A detailed investigation
of this phenomenon is beyond our present scope.

Overall, this example with quasar spectra coming from astronomy
illustrates the potential for optimally weighted PCA
to improve the recovery of underlying
principal components
from real data that has heteroscedastic noise.

%% file: content/conclusion.tex
%!TEX root = ../wpca.tex

\revone{%
\section{Conclusion}
\label{sec:conc}%
This paper derived asymptotic optimal weights
for weighted PCA
when the data is high-dimensional
with noise that is heteroscedastic across samples.
The optimal weights are a simple function
of the signal and noise variances,
and are not the inverse noise variance weights
commonly used in practice.
Numerical simulations
illustrated that
the asymptotic optimal weights
are often close to the optimal weights in finite dimensions
when the dimensions are large enough.
Comparisons of the asymptotic optimal weights
with existing weighting schemes
illustrated
that
the asymptotic optimal weights:
a) are generally closer to the optimal weights in finite dimensions,
b) appropriately combine all the data to achieve the best performance,
and
c) achieve nonzero asymptotic performance in the widest range of settings.
Additional simulations
illustrated that optimally weighted PCA
compares favorably with
other PCA methods
designed for some form of heteroscedastic noise.
Finally,
we briefly discussed how
one can use
estimated signal and noise variances
when the true variances are unknown,
and illustrated the optimal weights on real data from astronomy.

Overall,
optimally weighted PCA
is a simple, principled, and promising method
for estimating underlying principal components
from high-dimensional data
with noise that is heteroscedastic across samples.
However, many open questions remain.
While the asymptotic optimal performance
\cref{eq:opt:rec:asymp}
is often close
to the performance of optimally weighted PCA
in finite dimensions,
it would be useful to also derive
higher-order asymptotics for the performance
as well as nonasymptotic characterizations.
These results would provide refined estimates of the performance.
Another interesting variation of the asymptotic regime
is to allow the number of blocks $L$ to grow with $d$,
potentially with $n_1,\dots,n_L = \BigO(1)$.
This regime may better capture datasets
where each sample has its own associated noise variance
(e.g., like the astronomy data in \cref{sec:exp:astro})
or where the block structure is unknown.
While the proof of \cref{thm:opt:weight:asymp}
does not seem to readily generalize
to such cases,
we conjecture that
the optimal weights will still have the same form
and that estimates of the signal and noise variances
can still be used when the true variances are unknown.
Another interesting direction is to study
whether optimally weighted PCA
is optimal
across not just weighted PCA
but also across more general classes of methods,
e.g.,
by deriving fundamental bounds on the achievable performance.
Finally,
it would be interesting
to study how optimally weighted PCA
might be combined
with other techniques
to handle settings
where the noise is not just heteroscedastic
across samples
but also across features.
}%

%% file: acknowledgement.tex
%!TEX root = wpca.tex

\section*{Acknowledgements}
The authors thank
Raj Rao Nadakuditi,
Romain Couillet
and Edgar Dobriban
for helpful discussions
regarding the singular values and vectors
of low rank perturbations of large random matrices.
The authors also thank Dejiao Zhang for
helpful comments about the form
of the optimal weights.

The example quasar spectra were provided by the Sloan Digital Sky Survey
\cite{ahumada2020t1d,lyke2020tsd}.
Funding for the Sloan Digital Sky 
Survey IV has been provided by the 
Alfred P. Sloan Foundation, the U.S. 
Department of Energy Office of 
Science, and the Participating 
Institutions. 

SDSS-IV acknowledges support and 
resources from the Center for High 
Performance Computing at the 
University of Utah. The SDSS 
website is www.sdss.org.

SDSS-IV is managed by the 
Astrophysical Research Consortium 
for the Participating Institutions 
of the SDSS Collaboration including 
the Brazilian Participation Group, 
the Carnegie Institution for Science, 
Carnegie Mellon University, Center for 
Astrophysics | Harvard \& 
Smithsonian, the Chilean Participation 
Group, the French Participation Group, 
Instituto de Astrof\'isica de 
Canarias, The Johns Hopkins 
University, Kavli Institute for the 
Physics and Mathematics of the 
Universe (IPMU) / University of 
Tokyo, the Korean Participation Group, 
Lawrence Berkeley National Laboratory, 
Leibniz Institut f\"ur Astrophysik 
Potsdam (AIP), Max-Planck-Institut 
f\"ur Astronomie (MPIA Heidelberg), 
Max-Planck-Institut f\"ur 
Astrophysik (MPA Garching), 
Max-Planck-Institut f\"ur 
Extraterrestrische Physik (MPE), 
National Astronomical Observatories of 
China, New Mexico State University, 
New York University, University of 
Notre Dame, Observat\'ario 
Nacional / MCTI, The Ohio State 
University, Pennsylvania State 
University, Shanghai 
Astronomical Observatory, United 
Kingdom Participation Group, 
Universidad Nacional Aut\'onoma 
de M\'exico, University of Arizona, 
University of Colorado Boulder, 
University of Oxford, University of 
Portsmouth, University of Utah, 
University of Virginia, University 
of Washington, University of 
Wisconsin, Vanderbilt University, 
and Yale University.

%% file: content/ext,het,signal.tex
%!TEX root = ../wpca,supp.tex

\section{Extension to heterogeneous signal strengths}
\label{sec:het:signal}%
The model stated in \cref{sec:opt:weight}
has heteroscedastic noise
with a common signal covariance $\bmF \bmF^\CT$.
This model may at first appear to limit the scope of \cref{thm:opt:weight:asymp}
to only datasets with homogeneous signal strength.
However, as mentioned in \cref{rem:het:signal},
the result immediately generalizes
to cases with heterogeneous signal strengths $\tau_1,\dots,\tau_L > 0$.
Namely, suppose
$\bmY_1,\dots,\bmY_L$
are now generated as:
\begin{equation} \label{eq:model:blocks:hetsignal}
  \bmY_\ell = \sqrt{\tau_\ell} \bmF \bmZ_\ell + \bmE_\ell \in \bbC^{d \times n_\ell}
  , \quad \text{for } \ell = 1,\dots,L
  ,
\end{equation}
where
$\bmF \in \bbC^{d \times k}$,
$\bmZ_\ell \in \bbC^{k \times n_\ell}$,
and $\bmE_\ell \in \bbC^{d \times n_\ell}$
are as in \cref{eq:model:blocks}.
This model arises, e.g., in low-rank clutter estimation for RADAR
\cite{breloy2015cse,breloy2016rcm,sun2016lca}.

Under this model,
we have the following straightforward \lcnamecref{cor:opt:weight:asymp:hetsignal}.

\begin{corollary}[Extension to heterogeneous signal strengths]%
  \label{cor:opt:weight:asymp:hetsignal}%
  Let $\lambda_{i,\ell} \coloneqq \lambda_i \tau_\ell$.
  The optimal weights $\bmw^{\opt}_i(\bmY)$
  and corresponding optimal performance
  $r^{\opt}_i(\bmY)$
  converge as
  \begin{align} \label{eq:opt:weight:asymp:hetsignal}
    \bmw^{\opt}_i(\bmY)
    &\;\; \asto \;\;
    \bigg(
      \frac{1}{v_1}\frac{1}{1 + v_1/\lambda_{i,1}},
      \dots,
      \frac{1}{v_L}\frac{1}{1 + v_L/\lambda_{i,L}}
    \bigg)
    \qquad
    \text{up to scaling}
    , \\[0.5em]
    \label{eq:opt:rec:asymp:hetsignal}
    r^{\opt}_i(\bmY)
    &\;\; \asto \;\;
    \text{the unique solution $x \in (0,1)$ of }
    \sum_{\ell=1}^L
      \frac{c_\ell}{v_\ell/\lambda_{i,\ell}}
      \,\frac{1-x}{v_\ell/\lambda_{i,\ell}+x}
    = 1
    ,
  \end{align}
  except when
  $\sum_{\ell=1}^L c_\ell (\lambda_{i,\ell} /v_\ell)^2 \leq 1$,
  in which case
  $\bmu_i$ is asymptotically unrecoverable
  by any weighted PCA,
  i.e., $r_i(\bmw,\bmY) \asto 0$ for all weights $\bmw$.
\end{corollary}

\begin{proof}[Proof of \cref{cor:opt:weight:asymp:hetsignal}]
  Let $\btlY_\ell \coloneqq \bmY_\ell / \sqrt{\tau_\ell}$
  and $\tlv_\ell \coloneqq v_\ell/\tau_\ell$
  for $\ell = 1,\dots,L$.
  Note that
  $\bhtu_i(\bmw,\bmY) = \bhtu_i(\bmw \odot \bmtau,\btlY)$,
  so
  \begin{align*}
    \bmw^{\opt}_i(\btlY) \oslash \bmtau
    &\in
    \argmax_{\bmw}
    \; r_i(\bmw \odot \bmtau,\btlY)
    =
    \argmax_{\bmw}
    \; r_i(\bmw,\bmY)
    , \\
    r^{\opt}_i(\btlY)
    &=
    \max_{\bmw}
    \; r_i(\bmw \odot \bmtau,\btlY)
    =
    \max_{\bmw}
    \; r_i(\bmw,\bmY)
    ,
  \end{align*}
  where $\odot$ and $\oslash$ denote entrywise multiplication and division.
  Since $\btlY_1,\dots,\btlY_L$ obeys \cref{eq:model:blocks}
  with noise variances $\tlv_1,\dots,\tlv_L$,
  applying \cref{thm:opt:weight:asymp}
  and simplifying yields \cref{eq:opt:weight:asymp:hetsignal,eq:opt:rec:asymp:hetsignal}.
\end{proof}

Note that the form of the weights is the same as before:
inverse noise variance weights times
an SNR-dependent term that downweights low-SNR blocks.
Moreover, note that data with homoscedastic noise
and heterogeneous signal strengths
are uniformly weighted by inverse noise variance weights
but optimal weights downweight the blocks with weaker signals.
\revone{%
When the
signal strengths $\bmtau$,
signal variance $\lambda_i$,
and noise variances $\bmv$
are unknown,
estimates of these parameters
may be obtained
using any one of
various existing approaches,
and which one is most appropriate
will depend on the specific application.
Notably,
the simple estimators in \cref{ex:var:est}
could still be used in this case
to estimate the signal and noise variances,
and the same ideas may be adapted to also estimate the signal strengths.
For example,
one could apply the same signal variance estimator on each block individually,
then divide the resulting per-block signal variance estimates
by the signal variance estimated from the data overall.
}%

%% file: content/supp,derivation,calc.tex
%!TEX root = ../wpca,supp.tex

\section{Proof of \cref{lem:conv:argmax}}
\label{calc:proof:lem:conv:argmax}%
First, we prove convergence of the maximum, i.e., that
\begin{equation*}
  \forall_{\varepsilon > 0} \; \exists_{N} \; \forall_{n > N}
  \quad \Big| \max_{\bmx \in \clX} f_n(\bmx) - \max_{\bmx \in \clX} f(\bmx) \Big| < \varepsilon
  .
\end{equation*}
Let $\varepsilon > 0$ be arbitrary,
and let $\bmx^{\opt}$ be the unique maximizer of $f$.
Since $f_n \to f$ uniformly,
there exists $N$ such that
\begin{equation*}
  \forall_{n > N} \; \forall_{\bmx \in \clX} \quad |f_n(\bmx) - f(\bmx)| < \varepsilon
  .
\end{equation*}
Now, let $n > N$ be arbitrary.
Since $f_n$ has a maximum,
there exists $\btlx \in \argmax_{\bmx \in \clX} f_n(\bmx)$.
Noting that
\begin{align*}
  \max_{\bmx \in \clX} f_n(\bmx)
  &= f_n(\btlx)
  < f(\btlx) + \varepsilon
  \leq f(\bmx^{\opt}) + \varepsilon
  = \max_{\bmx \in \clX} f(\bmx) + \varepsilon
  , \\
  \max_{\bmx \in \clX} f_n(\bmx)
  &\geq f_n(\bmx^{\opt})
  > f(\bmx^{\opt}) - \varepsilon
  = \max_{\bmx \in \clX} f(\bmx) - \varepsilon
  ,
\end{align*}
yields $|\max_{\bmx \in \clX} f_n(\bmx) - \max_{\bmx \in \clX} f(\bmx)| < \varepsilon$,
establishing convergence for the maximum.

Next, we prove convergence of the set of maximizers.
Let $\bmx^{\opt}$ be the unique maximizer of $f$.
We want to show that the set of maximizers of $f_n$
(nonempty by assumption)
converges to $\{\bmx^{\opt}\}$
with respect to the Hausdorff distance,
i.e.,
that
$\sup\{ \|\btlx - \bmx^{\opt}\|_2 : \btlx \in \argmax_{\bmx \in \clX} f_n(\bmx) \} \to 0$,
or equivalently
\begin{equation*}
  \forall_{\delta > 0} \; \exists_{N} \; \forall_{n > N}
  \quad \forall_{\btlx \in \argmax_{\bmx \in \clX} f_n(\bmx)} \; \|\btlx - \bmx^{\opt}\|_2 < \delta
  .
\end{equation*}
Let $\delta > 0$ be arbitrary.
Note that $\clK \coloneqq \{\bmx \in \clX : \|\bmx - \bmx^{\opt}\|_2 \geq \delta\}$ is compact
since $\clX$ is compact.
Since $f$ is continuous,
it then follows from the extreme value theorem that
$f$ is maximized over $\clK$ at some point $\bhtx$.
Next, since $f_n \to f$ uniformly and the maximizer $\bmx^{\opt}$ is unique,
there exists $N$ such that
\begin{equation*}
  \forall_{n > N} \; \forall_{\bmx \in \clX} \quad |f_n(\bmx) - f(\bmx)| < \frac{f(\bmx^{\opt}) - f(\bhtx)}{2}
  .
\end{equation*}
Now, let $n > N$ be arbitrary.
Let $\btlx \in \argmax_{\bmx \in \clX} f_n(\bmx)$.
Then
\begin{equation*}
  f(\bmx^{\opt}) - f(\btlx)
  = \underbrace{f(\bmx^{\opt}) - f_n(\bmx^{\opt})}_{< \; (f(\bmx^{\opt}) - f(\bhtx))/2}
  + \underbrace{f_n(\bmx^{\opt}) - f_n(\btlx)}_{\leq \; 0}
  + \underbrace{f_n(\btlx) - f(\btlx)}_{< \; (f(\bmx^{\opt}) - f(\bhtx))/2}
  < f(\bmx^{\opt}) - f(\bhtx)
  ,
\end{equation*}
so $f(\btlx) > f(\bhtx) = \max_{\bmx \in \clK} f(\bmx)$.
Thus $\btlx \notin \clK$,
and so $\|\btlx - \bmx^{\opt}\|_2 < \delta$,
concluding the proof.

\section{Detailed calculations from the proof of \cref{lem:rec:lim}}%
This section provides detailed calculations
from the proof of \cref{lem:rec:lim}
in \cref{sct:rec:lim}.

\subsection{Anisotropic local law} \label{sec:local:law}%
This section states an \emph{anisotropic local law} for $\bmG(\bmw,\zeta)$
in a form that will be convenient for the rest of the proof.
Roughly speaking, it says that
the random resolvent matrix $\bmG(\bmw,\zeta)$ is well approximated
in the limit
by the deterministic matrix
\begin{equation}
  \label{defn_pi}
  \bmPi(\bmw,\zeta)
  \coloneqq
  \begin{bmatrix}
    \zeta \mplaw(\bmw,\zeta) \bmI_p & \\
    & \big( \zeta \bmI_n - c^{-1} \zeta \mplaw(\bmw,\zeta) \cdot \bmV \bmW \big)^{-1}
  \end{bmatrix}
  ,
\end{equation}
where
$\bmV \coloneqq \blockdiag(v_1 \bmI_{n_1},\dots,v_L \bmI_{n_L})$,
$\bmW \coloneqq \blockdiag(w_1 \bmI_{n_1},\dots,w_L \bmI_{n_L})$,
and $\mplaw(\bmw,\zeta)$
is the Stieltjes transform defined in \cref{eq:stieltjes:def}.
The local law is stated over two domains:
$\Sedge$ (around the edge $b_{\bmw}$ of the support of $\mu_{\bmw}$)
and $\Sout$ (outside the support of $\mu_{\bmw}$).
It follows straightforwardly from \cite{ding2021ssc,xi2020coe}.
Similar local laws have also been proved in \cite{bloemendal2014ill,knowles2016all,yang2019euo}
under various different assumptions.

\begin{lemma}[Anisotropic local laws]
  \label{lem_localin}
  Under the model assumptions,
  there exists a constant $\kappa>0$
  such that for any constants $\tau, \delta > 0$
  and any deterministic sequences of unit vectors $\bmxi_n, \bmnu_n \in \bbC^{d+n}$,
  \begin{alignat}{9}
    &
     \label{aniso_in}
      \sup_{\bmw \in \Delta_L} &\;\;& \sup_{\zeta \in \Sedge(\bmw,\tau,\delta;\kappa)}
      &\quad&
      \big| &\,& \bmxi_n^\CT &\,& [\bmG(\bmw,\zeta) - \bmPi(\bmw,\zeta)] &\,& \bmnu_n &\,& \big|
      &\qquad&
      \asto
      &\qquad&
      0
    ,
    \\
    &
     \label{aniso_outstrong}
      \sup_{\bmw \in \Delta_L} && \sup_{\zeta \in \Sout(\bmw,\tau)}
      &&
      \big| && \bmxi_n^\CT && [\bmG(\bmw,\zeta) - \bmPi(\bmw,\zeta)] && \bmnu_n && \big|
      &&
      \asto
      &&
      0
    ,
  \end{alignat}
  where
  \begin{alignat}{4}
    \label{eq_parain}
    \Sedge(\bmw&,\tau,\delta;\kappa&&)
    \coloneqq \big\{ \zeta \in \bbC :
      (\Re\zeta,\Im\zeta) \in [b_{\bmw} - \kappa, b_{\bmw} + \tau&&) \times [n^{-1+\delta} &&, \tau]
    \big\}
    , \\
    \label{eq_paraout}
    \Sout(\bmw&,\tau&&)
    \coloneqq \big\{ \zeta \in \bbC :
      (\Re\zeta,\Im\zeta) \in [b_{\bmw} + \tau,\infty&&) \times [-1 &&, 1]
    \big\}
    .
  \end{alignat}
\end{lemma}

\begin{proof}[Proof of \cref{lem_localin}]
Let
$\bbrE \coloneqq \bmE \bmV^{-1/2}$,
i.e., $\bbrE$ is whitened noise.
It follows immediately from the moment condition on $\bmE$ that
there exist constants $a>4$ and $C_0>0$
such that
\begin{equation}
  \label{condition_4e}
  \max_{i,j} \bbE |\bbrE_{ij}|^a \leq C_0
  .
\end{equation}
The local laws \cite{ding2021ssc, knowles2016all,xi2020coe}
are stated for matrices with truncated entries,
so we define
\begin{equation} \label{eq:truncated:noise}
  \bhtE_{\bmw} \coloneqq \bhtE \bmV^{1/2} \bmW^{1/2}
  \quad \text{where} \quad
  \bhtE_{ij} \coloneqq
  \begin{cases}
    \bbrE_{ij} & \text{if } |\bbrE_{ij}| \leq \phi_n n^{\varepsilon_0}
    , \\
    0 & \text{otherwise}
    ,
  \end{cases}
\end{equation}
where $\phi_n \coloneqq n^{2/a}$ and $\varepsilon_0 \in (0,1/2-2/a)$ is a constant to be chosen later.
Likewise, let $\bhtG(\bmw,\zeta)$ be the analogue of $\bmG(\bmw,\zeta)$
with $\bhtE_{\bmw}$ used in place of $\btlE_{\bmw}$.

Note that it follows from \cref{condition_4e} and integration by parts that
\begin{align}
  \label{error_3rd4th}
  | \bbE \bhtE_{ij} | &\leq n^{-3/2-\varepsilon_0}
  , &
  \bbE | \bhtE_{ij} |^2 &\leq n^{-1-\varepsilon_0}
  , &
  \bbE|\bhtE_{ij}|^3 &= \BigO(1)
  , &
  \bbE|\bhtE_{ij}|^4 &= \BigO(1)
  .
\end{align}
Moreover,
the empirical spectral distribution of $\bmV^{1/2}\bmW^{1/2}$ satisfies \cite[Example 2.8]{knowles2016all},
which leads to the edge regularity condition
\begin{equation}
  \label{eq:edge:regularity}
  \forall_{\bmw \in \Delta_L}
  \quad
  \min_{\ell} \big[ 1 - c^{-1} \mplaw(\bmw,b_{\bmw}) \cdot v_\ell w_\ell \big]
  > 0
  .
\end{equation}
Under
\cref{error_3rd4th},
$\max_{i,j}|\widehat{\bmE}_{ij}| \leq \phi_n n^{\varepsilon_0}$,
and the edge regularity condition \cref{eq:edge:regularity},
the following local laws follow immediately from \cite{ding2021ssc,xi2020coe}:
for any $\bmw \in \Delta_L$,
there exists a constant $\kappa>0$
such that for any constants $\tau, \delta > 0$,
any deterministic unit vectors $\bmxi, \bmnu \in \bbC^{d+n}$,
and any $\varepsilon_1, C_1 > 0$,
eventually (i.e., for $d$ large enough),
with probability at least $1 - n^{-C_1}$
we have
\begin{alignat}{8}
  \label{aniso_in2}
  &
  \sup_{\zeta \in \Sedge(\bmw,\tau,\delta;\kappa)}
  &\;\;&
  \big| &\,& \bmxi^\CT &\,& [\bhtG(\bmw,\zeta) - \bmPi(\bmw,\zeta)] &\,& \bmnu &\,& \big|
  &\quad&
  \leq
  &\quad&
  \phi_n n^{-1/2+\varepsilon_0+\varepsilon_1} + \frac{n^{\varepsilon_1}}{\sqrt{n \Im\zeta}}
  ,
  \\
  \label{aniso_outstrong2}
  &
  \sup_{\zeta \in \Sout(\bmw,\tau)}
  &&
  \big| && \bmxi^\CT && [\bhtG(\bmw,\zeta) - \bmPi(\bmw,\zeta)] && \bmnu && \big|
  &&
  \leq
  &&
  \phi_n n^{-1/2+\varepsilon_0+\varepsilon_1}
  .
\end{alignat}
Namely, \cref{aniso_in2} follows from \cite[Theorem 3.4]{xi2020coe}
and \cref{aniso_outstrong2} follows from \cite[Theorem S.3.12]{ding2021ssc}.

Using a standard $\varepsilon$-net argument
and taking a union bound with respect to $\bmw$ 
yields \cref{aniso_in2,aniso_outstrong2} uniformly in $\bmw \in \Delta_L$.
Moreover, choosing $\varepsilon_0$ and $\varepsilon_1$ small enough
yields $\phi_n n^{-1/2+\varepsilon_0+\varepsilon_1} \to 0$
and $n^{\varepsilon_1}/\sqrt{n \Im\zeta} \leq n^{-\delta/2+\varepsilon_1} \to 0$
for all $\zeta \in \Sedge(\bmw,\tau,\delta;\kappa)$.
Thus, applying the Borel-Cantelli lemma
yields analogues of \cref{aniso_in,aniso_outstrong}
with $\bhtE$ in place of $\bbrE$.
The proof concludes by observing that%
\footnote{In the derivation,
we regard $n \equiv n(d)$ as a sequence indexed by $d$,
which has the same order as $d$.}
\begin{align}
  \label{XneXio}
  \Pr( \bhtE \neq \bbrE \; \text{i.o.} )
  &
  = \lim_{k \to \infty}
  \Pr\big(
    \cup_{d=2^k}^\infty \cup_{i=1}^d \cup_{j=1}^{n(d)}
    \{ |\bbrE_{ij}| \geq \phi_n n^{\varepsilon_0} \}
  \big)
  \\&
  \nonumber
  = \lim_{k \to \infty}
  \Pr\big(
    \cup_{t=k}^\infty \cup_{d \in [2^t,2^{t+1})} \cup_{i=1}^d \cup_{j=1}^{n(d)}
    \{ |\bbrE_{ij}| \geq \phi_n n^{\varepsilon_0} \}
  \big)
  \\&
  \nonumber
  \leq C_0 \lim_{k \to \infty}
  \sum_{t=k}^\infty (2^{t+1})^2 (2^{t(2/a+\varepsilon_0)})^{-a}
  \leq 4 C_0 \lim_{k \to \infty} \sum_{t=k}^\infty 2^{-t a \varepsilon_0}
  = 0
  ,
\end{align}
i.e. $\bhtE = \bbrE$ eventually, almost surely.
\end{proof}

\subsection{Detailed derivation of \cref{eq:conv:noise}}
\label{calc:conv:noise}%
It follows immediately from \cite[Theorem 3.7]{xi2020coe} that
\begin{equation}
\label{noiseopnorm}
\Pr\bigg\{
  \sup_{\bmw \in \Delta_L} \big| \|\bhtE_{\bmw}\|_{\opnorm} - b_{\bmw} \big|
  \to 0
\bigg\}
= 1
,
\end{equation}
where $\bhtE_{\bmw}$ is the truncated weighted noise
defined in \cref{eq:truncated:noise}.
Combining \cref{noiseopnorm}
with \cref{XneXio}
(i.e., the fact that $\bhtE_{\bmw} = \btlE_{\bmw}$ eventually, almost surely)
yields \cref{eq:conv:noise}.

\subsection{Detailed derivation of \cref{eq:singval:kernel}}
\label{calc:proof:lem:singval:kernel}%
\Cref{eq:singval:kernel} modifies \cite[Lemma 4.1]{benaychgeorges2012tsv}
to account for the weights
and is proved in essentially the same way.
Namely,
recall that
by \cite[Theorem 7.3.3]{horn2013ma},
$\zeta$ is a singular value of $\btlY_{\bmw}$
if and only if it is a root of the characteristic polynomial
\begin{align}
  0
  &
  =
  \det \Bigg(
    \zeta \bmI_{d+n}
    -
    \begin{bmatrix} & \btlY_{\bmw} \\ \btlY_{\bmw}^\CT & \end{bmatrix}
  \Bigg)
  \\&
  \label{p:eq:singval:kernel:det1}
  =
  \det \Bigg(
    \underbrace{
    \zeta \bmI_{d+n}
    -
    \begin{bmatrix} & \btlE_{\bmw} \\ \btlE_{\bmw}^\CT & \end{bmatrix}
    }_{\bmA}
    -
    \underbrace{
    \begin{bmatrix} \bmU & \\ & \btlQ_{\bmw} \end{bmatrix}
    }_{\bmB}
    \underbrace{
    \begin{bmatrix} & \bmTheta \\ \bmTheta & \end{bmatrix}
    }_{\bmD}
    \underbrace{
    \begin{bmatrix} \bmU & \\ & \btlQ_{\bmw} \end{bmatrix}^\CT
    }_{\bmC}
  \Bigg)
  \\&
  \label{p:eq:singval:kernel:det2}
  =
  \det(\bmA)
  \det(\bmD)
  \det ( -\bmM(\bmw,\zeta) )
  ,
\end{align}
where \cref{p:eq:singval:kernel:det1}
is a convenient form of the matrix,
and \cref{p:eq:singval:kernel:det2}
follows from the identity
\begin{equation}
  \det(\bmA-\bmB\bmD\bmC)
  =
  \det(\bmA)
  \det(\bmD)
  \det(\bmD^{-1}-\bmC\bmA^{-1}\bmB)
  ,
\end{equation}
for invertible matrices $\bmA$ and $\bmD$.
Note that
$\bmA$ above
is invertible
because $\zeta > \|\btlE_{\bmw}\|_{\opnorm}$
(indeed, one only needs that $\zeta$ is not a singular value of $\btlE_{\bmw}$).
As a result, \cref{p:eq:singval:kernel:det2}
is zero if and only if $\det \bmM(\bmw,\zeta) = 0$,
yielding \cref{eq:singval:kernel}.

\subsection{Detailed derivation of \cref{wtMas}}
\label{calc:lem:kernel:lims}%
Let $\tau > 0$ be arbitrary,
and decompose $\bmM$ as
\begin{align}
    \label{p:lim:diff:0}
    \bmM(\bmw,\zeta)
    &
    =
    \begin{bmatrix} \bmU & \\ & \btlQ_{\bmw} \end{bmatrix}^\CT
    [\bmG(\bmw,\zeta) - \bmPi(\bmw,\zeta)]
    \begin{bmatrix} \bmU & \\ & \btlQ_{\bmw} \end{bmatrix}
    \\&\qquad
    \nonumber
    +
    \begin{bmatrix} \bmU & \\ & \btlQ_{\bmw} \end{bmatrix}^\CT
    \bmPi(\bmw,\zeta)
    \begin{bmatrix} \bmU & \\ & \btlQ_{\bmw} \end{bmatrix}
    -
    \begin{bmatrix} & \bmTheta^{-1} \\ \bmTheta^{-1} & \end{bmatrix}
    .
\end{align}
We will show that the first term vanishes
and the final two terms converge to $\bbrM(\bmw,\zeta)$.

For the first term of \cref{p:lim:diff:0},
applying \cref{aniso_outstrong} from \cref{lem_localin}
yields
\begin{equation}
  \label{p:lim:diff:1}
  \Bigg\|
    \begin{bmatrix} \bmU & \\ & \btlQ_{\bmw} \end{bmatrix}^\CT
    [\bmG(\bmw,\zeta) - \bmPi(\bmw,\zeta)]
    \begin{bmatrix} \bmU & \\ & \btlQ_{\bmw} \end{bmatrix}
  \Bigg\|_{\opnorm}
  \qquad
  \asunifto{(\bmw,\zeta) \in \Omega(\tau)}
  \qquad
  0
  ,
\end{equation}
where we use the observation that
$\Omega(\tau) = \{(\bmw,\zeta) \in \Delta_L \times \bbC : \zeta \in \Sout(\bmw,\tau)\}$.

For the second term of \cref{p:lim:diff:0}, note that
\begin{align}
  \label{p:lim:diff:2}
  &
  \begin{bmatrix} \bmU & \\ & \btlQ_{\bmw} \end{bmatrix}^\CT
  \bmPi(\bmw,\zeta)
  \begin{bmatrix} \bmU & \\ & \btlQ_{\bmw} \end{bmatrix}
  \\&\qquad
  \nonumber
  =
  \begin{bmatrix}
    \zeta \mplaw(\bmw,\zeta) \bmI_k & \\
    & \btlQ_{\bmw}^\CT \big( \zeta \bmI_n - c^{-1} \zeta \mplaw(\bmw,\zeta) \cdot \bmV \bmW \big)^{-1} \btlQ_{\bmw}
  \end{bmatrix}
  .
\end{align}
With the regularity condition \cref{eq:edge:regularity},
it follows immediately from \cite[Lemma A.6]{ding2018ana} that
\begin{equation}
  \label{eq:strong:regularity}
  \inf_{\bmw \in \Delta_L} \;\;
  \inf_{\zeta \in \Sedge(\bmw,\tau,\delta;\kappa) \cup \Sout(\bmw,\tau)} \;\;
  \min_{i=1,\dots,n} \bigg[1 - c^{-1} \mplaw(\bmw,\zeta) \cdot (\bmV \bmW)_{ii} \bigg]
  \quad
  >
  \quad
  0
  .
\end{equation}
Note that \cite[Lemma A.6]{ding2018ana} states the result for fixed $\bmw$;
taking the minimum over $\bmw$ in the compact set $\Delta_L$
yields \cref{eq:strong:regularity}.
Thus, $( \zeta \bmI_n - c^{-1} \zeta \mplaw(\bmw,\zeta) \cdot \bmV \bmW )^{-1}$
is a diagonal matrix with bounded entries.
Then using a standard $\varepsilon$-net argument
together with the law of large numbers
yields that
\begin{align}
  \label{p:lim:diff:3}
  &
  \btlQ_{\bmw}^\CT \big( \zeta \bmI_n - c^{-1} \zeta \mplaw(\bmw,\zeta) \cdot \bmV \bmW \big)^{-1} \btlQ_{\bmw}
  \\&\quad
  \nonumber
  -
  \frac{1}{n}
  \tr\bigg[
    \bmW^{1/2} \big( \zeta \bmI_n - c^{-1} \zeta \mplaw(\bmw,\zeta) \cdot \bmV \bmW \big)^{-1} \bmW^{1/2}
  \bigg]
  \cdot
  \bmI_k
  \quad
  \asunifto{(\bmw,\zeta) \in \Omega(\tau)}
  \quad
  \bm0_{k \times k}
  .
\end{align}
Note also that
\begin{equation}
  \label{p:lim:diff:4}
  \frac{1}{n}
  \tr\bigg[
    \bmW^{1/2} \big( \zeta \bmI_n - c^{-1} \zeta \mplaw(\bmw,\zeta) \cdot \bmV \bmW \big)^{-1} \bmW^{1/2}
  \bigg]
  =
  \sum_{\ell=1}^L \frac{p_\ell w_\ell}{\zeta - w_\ell v_\ell \zeta \mplaw(\bmw,\zeta) /c}
  .
\end{equation}
Combining \cref{p:lim:diff:0,p:lim:diff:1,p:lim:diff:2,p:lim:diff:3,p:lim:diff:4}
with $\zeta\mplaw(\bmw,\zeta) = \varphi_{1,\bmw}(\zeta)$
yields \cref{wtMas}.

\subsection{Verification of properties
for \texorpdfstring{$\varphi_{1,\bmw}$}{phi1w}
and \texorpdfstring{$\varphi_{2,\bmw}$}{phi2w}}
\label{calc:traceconv:props}%
We verify that,
for any $\bmw \in \Delta_L$,
the functions
$\varphi_{1,\bmw}$ and $\varphi_{2,\bmw}$ have the following properties:
\begin{align}
  \label{eq:traceconv:props:1}
  &\forall_{\zeta > b_{\bmw}} \; \varphi_{1,\bmw}(\zeta) > 0
  , &
  &\varphi_{1,\bmw}(\zeta) \to 0 \text{ as } |\zeta| \to \infty
  , &
  &\forall_{\zeta \notin \supp(\mu_{\bmw})} \;
  \varphi_{1,\bmw}(\zeta) \in \bbR \Leftrightarrow \zeta \in \bbR
  , \\
  \label{eq:traceconv:props:2}
  &\forall_{\zeta > b_{\bmw}} \; \varphi_{2,\bmw}(\zeta) > 0
  , &
  &\varphi_{2,\bmw}(\zeta) \to 0 \text{ as } |\zeta| \to \infty
  , &
  &\forall_{\zeta \notin \supp(\mu_{\bmw})} \;
  \varphi_{2,\bmw}(\zeta) \in \bbR \Leftrightarrow \zeta \in \bbR
  ,
\end{align}
where $\supp(\mu_{\bmw})$ denotes the support of $\mu_{\bmw}$.
First, we verify the properties in \cref{eq:traceconv:props:1} for $\varphi_{1,\bmw}$:
\begin{enumerate}[label=(\alph*),parsep=1em,topsep=1em]
  \item
  For any $\zeta > b_{\bmw}$,
  the integrand in \cref{p:eq:vp1}
  is positive and bounded away from zero since
  the support of $\mu_{\bmw}$
  lies between zero and $b_{\bmw}$.

  Thus, $\forall_{\zeta > b_{\bmw}} \; \varphi_{1,\bmw}(\zeta) > 0$.

  \item
  As $|\zeta| \to \infty$,
  the integrand in \cref{p:eq:vp1} goes to zero uniformly in $t$.

  Thus, $\varphi_{1,\bmw}(\zeta) \to 0$ as $|\zeta| \to \infty$.

  \item
  For $\zeta \notin \supp(\mu_{\bmw})$,
  the imaginary part of $\varphi_{1,\bmw}(\zeta)$ is
  \begin{equation*}
    \Im \varphi_{1,\bmw}(\zeta)
    =
    \int \Im\bigg( \frac{\zeta}{\zeta^2-t^2} \bigg) \; \rmd \mu_{\bmw}(t)
    =
    - \Im\zeta
    \underbrace{
      \int \frac{|\zeta|^2+t^2}{|\zeta^2-t^2|^2} \; \rmd \mu_{\bmw}(t)
    }_{> 0}
    .
  \end{equation*}
  Thus,
  $\forall_{\zeta \notin \supp(\mu_{\bmw})} \;
  \varphi_{1,\bmw}(\zeta) \in \bbR \Leftrightarrow \zeta \in \bbR$.
\end{enumerate}
Next we verify the properties in \cref{eq:traceconv:props:2} for $\varphi_{2,\bmw}$:
\begin{enumerate}[label=(\alph*),parsep=1em,topsep=1em]
  \item
  For any $\zeta > b_{\bmw}$,
  $\mplaw(\bmw,\zeta) \leq \mplaw(\bmw,b_{\bmw})$
  since the integrand in \cref{eq:stieltjes:def} is decreasing in $\zeta$.
  So, it follows from \cref{eq:edge:regularity} that
  for $\ell = 1,\dots,L$,
  \begin{equation*}
    1 - c^{-1} \mplaw(\bmw,\zeta) \cdot v_\ell w_\ell
    \geq
    1 - c^{-1} \mplaw(\bmw,b_{\bmw}) \cdot v_\ell w_\ell
    > 0
    ,
  \end{equation*}
  so the denominator of each summand
  in the definition of $\varphi_{2,\bmw}$ in \cref{p:eq:vp1}
  is also positive,
  i.e.,
  $\zeta - w_\ell v_\ell \varphi_{1,\bmw}(\zeta)/c
  = \zeta [ 1 - c^{-1} \mplaw(\bmw,\zeta) \cdot v_\ell w_\ell ]
  > 0$.

  Thus, $\forall_{\zeta > b_{\bmw}} \; \varphi_{2,\bmw}(\zeta) > 0$.

  \item
  As $|\zeta| \to \infty$,
  $|\zeta - w_\ell v_\ell \varphi_{1,\bmw}(\zeta) / c| \to \infty$
  for each $\ell \in \{1,\dots,L\}$
  since $\varphi_{1,\bmw}(\zeta) \to 0$
  as shown above.

  Thus, $\varphi_{2,\bmw}(\zeta) \to 0$ as $|\zeta| \to \infty$.

  \item
  As shown above,
  for any $\zeta \notin \supp(\mu_{\bmw})$,
  $\Im \varphi_{1,\bmw}(\zeta)$ is zero
  if $\Im\zeta$ is zero
  and has the opposite sign of $\Im\zeta$ otherwise.
  As a result,
  \begin{equation*}
    \Im\{ \zeta - w_\ell v_\ell \varphi_{1,\bmw}(\zeta) / c \}
    =
    \Im\zeta - (w_\ell v_\ell / c) \Im \varphi_{1,\bmw}(\zeta)
  \end{equation*}
  is zero if $\Im\zeta$ is zero
  and has the same sign as $\Im\zeta$ otherwise.
  
  Thus,
  $\forall_{\zeta \notin \supp(\mu_{\bmw})} \;
  \varphi_{2,\bmw}(\zeta) \in \bbR \Leftrightarrow \zeta \in \bbR$.
\end{enumerate}
As a result,
$\varphi_{1,\bmw}$ and $\varphi_{2,\bmw}$ satisfy
the needed conditions of \cite[Lemma A.1]{benaychgeorges2012tsv}.

\subsection{Detailed derivation of \cref{lem:singvec:kernel}}
\label{calc:lem:singvec:kernel}%
\Cref{lem:singvec:kernel} modifies \cite[Lemma 5.1]{benaychgeorges2012tsv}
to account for the weights
and is proved in essentially the same way.
Namely, let $\bmw$ be such that $\httheta_{i,\bmw} > \|\btlE_{\bmw}\|_{\opnorm}$
and let $\btlX_{\bmw} \coloneqq \bmU\bmTheta\btlQ_{\bmw}^\CT$
be the weighted and normalized signal.

To derive \cref{p:eq:singvec:kernel},
we first use the block matrix inverse
\cite[Equation (0.7.3.1)]{horn2013ma}
to write the resolvent $\bmG(\bmw,\zeta)$ from \cref{eq_resolvent}
in the following block matrix form:
\begin{equation}
\label{eq_resolvent:block}
\bmG(\bmw,\zeta)
=
\begin{bmatrix}
  \zeta (\zeta^2 \bmI_d - \btlE_{\bmw}\btlE_{\bmw}^\CT)^{-1} &
  (\zeta^2 \bmI_d - \btlE_{\bmw}\btlE_{\bmw}^\CT)^{-1} \btlE_{\bmw}
  \\
  \btlE_{\bmw}^\CT (\zeta^2 \bmI_d - \btlE_{\bmw}\btlE_{\bmw}^\CT)^{-1} &
  \zeta (\zeta^2 \bmI_n - \btlE_{\bmw}^\CT\btlE_{\bmw})^{-1}
\end{bmatrix}
.
\end{equation}
Then \cref{p:eq:singvec:kernel} follows by
substituting \cref{p:eq:singval:kernel,eq_resolvent:block}
then factoring.
Namely,
\begin{align}
  &
  \bmM(\bmw,\httheta_{i,\bmw})
  \begin{bmatrix}
  \bmTheta \btlQ_{\bmw}^\CT \bhtq_{i,\bmw} \\
  \bmTheta \bmU^\CT \bhtu_{i,\bmw}
  \end{bmatrix}
  \\&
  \label{p:eq:singvec:kernel:prod1}
  =
  \begin{bmatrix}
    \bmU^\CT
    \big(
      (\httheta_{i,\bmw}^2 \bmI_d - \btlE_{\bmw}\btlE_{\bmw}^\CT)^{-1}
      (\httheta_{i,\bmw} \btlX_{\bmw} \bhtq_{i,\bmw}
        + \btlE_{\bmw} \btlX_{\bmw}^\CT \bhtu_{i,\bmw})
      - \bhtu_{i,\bmw}
    \big)
    \\
    \btlQ_{\bmw}^\CT
    \big(
      (\httheta_{i,\bmw}^2 \bmI_n - \btlE_{\bmw}^\CT\btlE_{\bmw})^{-1}
      (\btlE_{\bmw}^\CT \btlX_{\bmw} \bhtq_{i,\bmw}
        + \httheta_{i,\bmw} \btlX_{\bmw}^\CT \bhtu_{i,\bmw})
      - \bhtq_{i,\bmw}
    \big)
  \end{bmatrix}
  \\&
  \label{p:eq:singvec:kernel:prod2}
  =
  \begin{bmatrix}
  \bmU^\CT (\bhtu_{i,\bmw} - \bhtu_{i,\bmw}) \\
  \btlQ_{\bmw}^\CT (\bhtq_{i,\bmw} - \bhtq_{i,\bmw})
  \end{bmatrix}
  = 0
  ,
\end{align}
where \cref{p:eq:singvec:kernel:prod1} uses
the identity
$\btlE_{\bmw}^\CT (\httheta_{i,\bmw}^2 \bmI_d - \btlE_{\bmw}\btlE_{\bmw}^\CT)^{-1}
=(\httheta_{i,\bmw}^2 \bmI_n - \btlE_{\bmw}^\CT\btlE_{\bmw})^{-1} \btlE_{\bmw}^\CT$,
and \cref{p:eq:singvec:kernel:prod2}
follows by substituting $\btlX_{\bmw} = \btlY_{\bmw} - \btlE_{\bmw}$
and using the singular vector identities
\begin{align*}
  \btlY_{\bmw} \bhtq_{i,\bmw} &= \httheta_{i,\bmw} \bhtu_{i,\bmw}
  , &
  \btlY_{\bmw}^\CT \bhtu_{i,\bmw} &= \httheta_{i,\bmw} \bhtq_{i,\bmw}
  .
\end{align*}

To derive \cref{p:eq:singvec:sum},
combine the identity
$\bhtu_{i,\bmw}
=
\bmGamma_{\bmw}
(
  \httheta_{i,\bmw} \btlX_{\bmw} \bhtq_{i,\bmw}
  + \btlE_{\bmw} \btlX_{\bmw}^\CT \bhtu_{i,\bmw}
)$
used to obtain \cref{p:eq:singvec:kernel:prod2}
with the fact that $\|\bhtu_{i,\bmw}\|_2 = 1$
and expand as
\begin{align*}
  1
  &= \bhtu_{i,\bmw}^\CT \bhtu_{i,\bmw}
  \\&
  =
  (
    \httheta_{i,\bmw} \btlX_{\bmw} \bhtq_{i,\bmw}
    + \btlE_{\bmw} \btlX_{\bmw}^\CT \bhtu_{i,\bmw}
  )^\CT
  \;
  \bmGamma_{\bmw}^2
  \;
  (
    \httheta_{i,\bmw} \btlX_{\bmw} \bhtq_{i,\bmw}
    + \btlE_{\bmw} \btlX_{\bmw}^\CT \bhtu_{i,\bmw}
  )
  \\&
  = \chi_1(\bmw) + \chi_2(\bmw) + 2 \Re \chi_3(\bmw)
  ,
\end{align*}
where the outer terms are
\begin{align}
  \label{p:eq:singvec:sum:parts:1}
  \chi_1(\bmw)
  &\coloneqq
  \bhtq_{i,\bmw}^\CT \btlX_{\bmw}^\CT
  \httheta_{i,\bmw}^2 \bmGamma_{\bmw}^2
  \btlX_{\bmw} \bhtq_{i,\bmw}
  , &
  \chi_2(\bmw)
  &\coloneqq
  \bhtu_{i,\bmw}^\CT \btlX_{\bmw} \btlE_{\bmw}^\CT
  \bmGamma_{\bmw}^2
  \btlE_{\bmw} \btlX_{\bmw}^\CT \bhtu_{i,\bmw}
  ,
\end{align}
and the cross term is
\begin{equation}
  \label{p:eq:singvec:sum:parts:2}
  \chi_3(\bmw)
  \coloneqq
  \bhtq_{i,\bmw}^\CT \btlX_{\bmw}^\CT
  \httheta_{i,\bmw} \bmGamma_{\bmw}^2
  \btlE_{\bmw} \btlX_{\bmw}^\CT \bhtu_{i,\bmw}
  .
\end{equation}
Expanding
$\btlX_{\bmw}
= \bmU\bmTheta\btlQ_{\bmw}^\CT
= \theta_1 \bmu_1 \btlq_{1,\bmw}^\CT
+ \cdots
+ \theta_k \bmu_k \btlq_{k,\bmw}^\CT$
in \cref{p:eq:singvec:sum:parts:1,p:eq:singvec:sum:parts:2}
then simplifying yields \cref{p:eq:singvec:sum:parts}.

\subsection{Detailed derivation of \cref{p:eq:singvec:sum:parts:lim:all}}
\label{calc:lem:singvec:sum:parts:lim}%
To derive \cref{p:eq:singvec:sum:parts:lim:all},
it is helpful to first establish
that several terms appearing in \cref{p:eq:singvec:sum:parts}
are bounded.
Namely,
almost surely,
eventually,
the following upper bounds hold
for all $\bmw \in \clW_>(\nu)$:
\begin{subequations}
\label{p:eq:bounded}
\begin{align}
  \label{p:eq:bounded:httheta}
  \httheta_{i,\bmw}
  &\leq \|\btlY_{\bmw}\|_{\opnorm}
  = \| \btlY_{\bm1_L} \bmW^{1/2} \|_{\opnorm}
  \leq \| \btlY_{\bm1_L} \|_{\opnorm} \| \bmW^{1/2} \|_{\opnorm}
  \leq \| \btlY_{\bm1_L} \|_{\opnorm}
  < 2 \cdot \brtheta_{1,\bm1_L}
  , \\
  \label{p:eq:bounded:bmGamma}
  \| \bmGamma_{\bmw} \|_{\opnorm}
  &= \| (\httheta_{i,\bmw}^2 \bmI_d - \btlE_{\bmw}\btlE_{\bmw}^\CT)^{-1} \|_{\opnorm}
  < 2 \cdot \frac{1}{\brtheta_{i,\bmw}^2 - b_{\bmw}^2}
  < 2 \cdot \frac{1}{\nu^2}
  , \\
  \label{p:eq:bounded:btlE}
  \| \btlE_{\bmw} \|_{\opnorm}
  &= \| \btlE_{\bm1_L} \bmW^{1/2} \|_{\opnorm}
  \leq \| \btlE_{\bm1_L} \|_{\opnorm} \| \bmW^{1/2} \|_{\opnorm}
  \leq \| \btlE_{\bm1_L} \|_{\opnorm}
  < 2 \cdot b_{\bm1_L}
  , \\
  \label{p:eq:bounded:btlqW}
  \| \btlq_{j,\bmw} \|_2
  &= \| \bmW^{1/2} \btlq_{j,\bm1_L} \|_2
  \leq \|\bmW^{1/2}\|_{\opnorm} \|\btlq_{j,\bm1_L}\|_2
  < 2
  , \\
  \label{p:eq:bounded:btlqWbtlq}
  |\btlq_{j,\bmw}^\CT \bhtq_{i,\bmw}|
  &\leq \|\btlq_{j,\bmw}\|_2 \|\bhtq_{i,\bmw}\|_2
  < 2
  ,\\
  \label{p:eq:bounded:bmubmu}
  |\bmu_j^\CT \bhtu_{i,\bmw}|
  &\leq \|\bmu_j\|_2 \|\bhtu_{i,\bmw}\|_2 = 1
\end{align}
\end{subequations}
where
\begin{itemize}
  \item \cref{p:eq:bounded:httheta} follows from
  the identity $\btlY_{\bmw} = \btlY_{\bm1_L} \bmW^{1/2}$,
  submultiplicativity of the operator norm,
  the fact that $\|\bmW^{1/2}\|_{\opnorm} = \|\bmw\|_\infty^{1/2} \leq 1$ for $\bmw \in \Delta_L$,
  and \cref{p:eq:singval:limit};
  \item \cref{p:eq:bounded:bmGamma} follows from \cref{p:eq:singval:limit,eq:conv:noise},
  and the fact that $\brtheta_{i,\bmw}^2 > (b_{\bmw} + \nu)^2 > b_{\bmw}^2 + \nu^2$
  for $\bmw \in \clW_>(\nu)$;
  \item \cref{p:eq:bounded:btlE} follows from
  the identity $\btlE_{\bmw} = \btlE_{\bm1_L} \bmW^{1/2}$,
  submultiplicativity of the operator norm,
  the fact that $\|\bmW^{1/2}\|_{\opnorm} = \|\bmw\|_\infty^{1/2} \leq 1$ for $\bmw \in \Delta_L$,
  and \cref{eq:conv:noise};
  \item \cref{p:eq:bounded:btlqW} follows from
  the identity $\btlq_{j,\bmw} = \bmW^{1/2} \btlq_{j,\bm1_L}$,
  the operator norm inequality,
  the fact that $\|\bmW^{1/2}\|_{\opnorm} = \|\bmw\|_\infty^{1/2} \leq 1$ for $\bmw \in \Delta_L$,
  and the fact that $\|\btlq_{j,\bm1_L}\|_2^2 \asto 1$ by the law of large numbers;
  \item \cref{p:eq:bounded:btlqWbtlq} follows from
  the Cauchy-Schwarz inequality
  and \cref{p:eq:bounded:btlqW};
  \item \cref{p:eq:bounded:bmubmu} follows from the Cauchy-Schwarz inequality.
\end{itemize}
In other words,
these terms are almost surely eventually uniformly bounded.

We now derive \cref{p:eq:singvec:sum:parts:lim,p:eq:singvec:sum:parts:lim:chi2}.
Note first that
almost surely,
eventually,
for any $\bmw \in \clW_>(\nu)$,
$\bmM(\bmw,\cdot)$ and $\bbrM(\bmw,\cdot)$
are both Lipschitz functions for $\zeta > b_{\bmw} + \nu/2$
with Lipschitz constant $\BigO(1/\nu^2)$.
Moreover,
almost surely,
eventually,
for all $\bmw \in \clW_>(\nu)$,
$\httheta_{i,\bmw} > b_{\bmw} + \nu/2$.
Hence,
it follows from the limits \cref{wtMas,p:eq:singval:limit} that
\begin{equation}
  \bmM(\bmw,\httheta_{i,\bmw})
  \quad
  \asunifto{\bmw \in \clW_>(\nu)}
  \quad
  \bbrM(\bmw,\rho_{i,\bmw})
  .
\end{equation}
Applying this to \cref{p:eq:singvec:kernel}
yields
\begin{equation} \label{p:eq:singvec:proj:lim}
  \begin{bmatrix} \xi(\bmw) \\ \eta(\bmw) \end{bmatrix}
  \coloneqq \proj_{(\ker \bbrM(\bmw,\rho_{i,\bmw}))^\perp}
  \begin{bmatrix}
    \bmTheta \btlQ_{\bmw}^\CT \bhtq_{i,\bmw} \\
    \bmTheta \bmU^\CT \bhtu_{i,\bmw}
  \end{bmatrix}
  \quad
  \asunifto{\bmw \in \clW_>(\nu)}
  \quad
  0
  .
\end{equation}
Observe next that, similar to~\cite[Section 5]{benaychgeorges2012tsv},
\begin{equation}
  \forall_{\bmw \in \clW_>(\nu)} \ \
  \ker \bbrM(\bmw,\rho_{i,\bmw}) =
  \Bigg\{
    \begin{bmatrix} s \\ t \end{bmatrix} \in \bbC^{2k} :
    \begin{aligned}
      t_j &= \theta_i \varphi_{1,\bmw}(\rho_{i,\bmw}) s_j &\text{for } j \text{ s.t. } \theta_j =    \theta_i \\
      t_j &= s_j = 0                                      &\text{for } j \text{ s.t. } \theta_j \neq \theta_i
    \end{aligned}
  \Bigg\}
  ,
\end{equation}
so the projection entries are
\begin{align}
  \label{p:eq:singvec:proj}
  \begin{bmatrix} \xi_i(\bmw) \\ \eta_i(\bmw) \end{bmatrix}
  &=
  \big(
  \theta_i \varphi_{1,\bmw}(\rho_{i,\bmw})
  \btlq_{i,\bmw}^\CT \bhtq_{i,\bmw}
  - \bmu_i^\CT \bhtu_{i,\bmw}
  \big)
  \frac{\theta_i
  }{\theta_i^2 \varphi_{1,\bmw}^2(\rho_{i,\bmw}) + 1}
  \begin{bmatrix} \theta_i \varphi_{1,\bmw}(\rho_{i,\bmw}) \\ -1 \end{bmatrix}
  , \\
  \nonumber
  \forall_{j : j \neq i} \quad
  \begin{bmatrix} \xi_j(\bmw) \\ \eta_j(\bmw) \end{bmatrix}
  &=
  \theta_j
  \begin{bmatrix}
  \btlq_{j,\bmw}^\CT \bhtq_{i,\bmw} \\
  \bmu_j^\CT \bhtu_{i,\bmw}
  \end{bmatrix}
  ,
\end{align}
Applying \cref{p:eq:singvec:proj:lim}
to \cref{p:eq:singvec:proj} yields
\begin{alignat}{3} \label{p:eq:singvec:lim:1}
  \sum_{j : j \neq i}
  |\bmu_j^\CT \bhtu_{i,\bmw}|^2
  +
  |\btlq_{j,\bmw}^\CT \bhtq_{i,\bmw}|^2
  &
  &\qquad&
  \asunifto{\bmw \in \clW_>(\nu)}
  &\qquad&
  0
  , \\
  \label{p:eq:singvec:lim:2}
  \Bigg|
    \sqrt{\frac{\varphi_{1,\bmw}(\rho_{i,\bmw})}{\varphi_{2,\bmw}(\rho_{i,\bmw})}}
    \btlq_{i,\bmw}^\CT \bhtq_{i,\bmw}
    - \bmu_i^\CT \bhtu_{i,\bmw}
  \Bigg|^2
  &
  &&
  \asunifto{\bmw \in \clW_>(\nu)}
  &&
  0
  ,
\end{alignat}
where we recall that
$D_{\bmw}(\rho_{i,\bmw}) = \varphi_{1,\bmw}(\rho_{i,\bmw})\varphi_{2,\bmw}(\rho_{i,\bmw}) = 1/\theta_i^2$
and note that
$\varphi_{1,\bmw}(\rho_{i,\bmw})$
and $\varphi_{2,\bmw}(\rho_{i,\bmw})$
are uniformly upper and lower bounded
with respect to $\bmw \in \Delta_L$.

Now, by \cref{p:eq:singvec:lim:1}
and the bounds \cref{p:eq:bounded},
it follows that
\begin{subequations}
\label{eq:chi1:chi2:offdiag}
\begin{alignat}{3}
  \label{eq:chi1:offdiag}
  &
  \sum_{j_1,j_2 : (j_1, j_2) \neq (i,i)}
  \theta_{j_1} \theta_{j_2}
  (\btlq_{j_1,\bmw}^\CT \bhtq_{i,\bmw})
  (\btlq_{j_2,\bmw}^\CT \bhtq_{i,\bmw})^*
  \bmu_{j_2}^\CT \httheta_{i,\bmw}^2 \bmGamma_{\bmw}^2 \bmu_{j_1}
  &\;\;&
  \asunifto{\bmw \in \clW_>(\nu)}
  &\;\;&
  0
  , \\
  \label{eq:chi2:offdiag}
  &
  \sum_{j_1,j_2 : (j_1, j_2) \neq (i,i)}
  \theta_{j_1} \theta_{j_2}
  (\bmu_{j_1}^\CT \bhtu_{i,\bmw})
  (\bmu_{j_2}^\CT \bhtu_{i,\bmw})^*
  \btlq_{j_2,\bmw}^\CT \btlE_{\bmw}^\CT \bmGamma_{\bmw}^2 \btlE_{\bmw} \btlq_{j_1,\bmw}
  &&
  \asunifto{\bmw \in \clW_>(\nu)}
  &&
  0
  ,
\end{alignat}
\end{subequations}
so it remains to analyze the summands
of $\chi_1$ and $\chi_2$
with $(j_1, j_2) = (i,i)$.
For this, note that
\begin{subequations}
\label{eq:vp:expand}
\begin{align}
  \bmu_i^\CT \httheta_{i,\bmw}^2 \bmGamma_{\bmw}^2 \bmu_i
  &
  =
  \bigg(+\frac{1}{2 \zeta} - \frac{1}{2} \diff{}{\zeta}\bigg)
  [\bmM(\bmw,\zeta)]_{i,i}
  \Bigg|_{\zeta=\httheta_{i,\bmw}}
  ,
  \\
  \btlq_{i,\bmw}^\CT \btlE_{\bmw}^\CT \bmGamma_{\bmw}^2 \btlE_{\bmw} \btlq_{i,\bmw}
  &
  =
  \bigg(-\frac{1}{2 \zeta} - \frac{1}{2} \diff{}{\zeta}\bigg)
  [\bmM(\bmw,\zeta)]_{k+i,k+i}
  \Bigg|_{\zeta=\httheta_{i,\bmw}}
  ,
\end{align}
\end{subequations}
which can be verified by substituting \cref{p:eq:singval:kernel,eq_resolvent:block},
expanding and simplifying.

Now recall that \cref{wtMas} yields, for $i=1,\dots,k$,
\begin{align}
  \label{eq:conv:vp}
  [\bmM(\bmw,\zeta)]_{i,i}
  &
  \;
  \asunifto{(\bmw,\zeta) \in \Omega(\tau)}
  \;
  \varphi_{1,\bmw}(\zeta)
  , &
  [\bmM(\bmw,\zeta)]_{k+i,k+i}
  &
  \;
  \asunifto{(\bmw,\zeta) \in \Omega(\tau)}
  \;
  \varphi_{2,\bmw}(\zeta)
  .
\end{align}
Moreover,
almost surely,
eventually,
for all $(\bmw,\zeta) \in \Omega(\tau)$
these are both holomorphic functions
with respect to $\zeta$,
because eventually,
for all $\bmw \in \Delta_L$,
$\|\btlE_{\bmw}\|_{\opnorm} < b_{\bmw} + \tau$.
Thus, by a standard application of Cauchy's integral formula,
their derivatives converge uniformly
on any compact subset.
Namely,
for any compact subset $\clC \subset \Omega(\tau)$,
for all $i = 1,\dots,k$,
\begin{align}
  \label{eq:conv:vp:diff}
  \diff{}{\zeta} [\bmM(\bmw,\zeta)]_{i,i}
  &
  \;
  \asunifto{(\bmw,\zeta) \in \clC}
  \;
  \varphi_{1,\bmw}'(\zeta)
  , &
  \diff{}{\zeta} [\bmM(\bmw,\zeta)]_{k+i,k+i}
  &
  \;
  \asunifto{(\bmw,\zeta) \in \clC}
  \;
  \varphi_{2,\bmw}'(\zeta)
  .
\end{align}

Thus, applying \cref{eq:conv:vp,eq:conv:vp:diff} to \cref{eq:vp:expand} yields
\begin{subequations}
\label{eq:vp:expand:conv}
\begin{align}
  \label{eq:vp:expand:conv:chi1}
  \bmu_i^\CT \httheta_{i,\bmw}^2 \bmGamma_{\bmw}^2 \bmu_i
  &
  \quad
  \asunifto{\bmw \in \clW_>(\nu)}
  \quad
  +\frac{\varphi_{1,\bmw}(\rho_{i,\bmw})}{2 \rho_{i,\bmw}} - \frac{\varphi_{1,\bmw}'(\rho_{i,\bmw})}{2}
  ,
  \\
  \label{eq:vp:expand:conv:chi2}
  \btlq_{i,\bmw}^\CT \btlE_{\bmw}^\CT \bmGamma_{\bmw}^2 \btlE_{\bmw} \btlq_{i,\bmw}
  &
  \quad
  \asunifto{\bmw \in \clW_>(\nu)}
  \quad
  -\frac{\varphi_{2,\bmw}(\rho_{i,\bmw})}{2 \rho_{i,\bmw}} - \frac{\varphi_{2,\bmw}'(\rho_{i,\bmw})}{2}
  ,
\end{align}
\end{subequations}
where we also used that
for any $\bmw \in \clW_>(\nu)$,
almost surely,
eventually,
$\bmM(\bmw,\cdot)$, $\bbrM(\bmw,\cdot)$ and their derivatives are Lipschitz functions for $\zeta > b_{\bmw} + \nu/2$
with Lipschitz constants uniformly bounded across $\bmw \in \clW_>(\nu)$.

Finally,
combining \cref{p:eq:singvec:lim:2,eq:vp:expand:conv:chi1,eq:chi1:offdiag}
yields \cref{p:eq:singvec:sum:parts:lim}.
Similarly, combining
\cref{eq:vp:expand:conv:chi2,eq:chi2:offdiag}
yields \cref{p:eq:singvec:sum:parts:lim:chi2}.
\Cref{p:eq:singvec:sum:parts:lim:chi3} follows from a similar argument.
More precisely,
combining \cref{p:eq:singval:kernel,eq_resolvent:block,wtMas} yields that
\begin{equation*}
  \bmu_{j_2}^\CT \httheta_{i,\bmw} \bmGamma_{\bmw}^2 \btlE_{\bmw} \btlq_{j_1,\bmw}
  =
  - \frac{1}{2} \diff{}{\zeta}
  [\bmM(\bmw,\zeta)]_{j_2,k+j_1}
  \Bigg|_{\zeta=\httheta_{i,\bmw}}
  \quad
  \asunifto{\bmw \in \clW_>(\nu)}
  \quad
  0
  .
\end{equation*}
Combining this with the bounds \cref{p:eq:bounded}
yields \cref{p:eq:singvec:sum:parts:lim:chi3}.

\subsection{Detailed derivation of \cref{eq:singvec:lim:above}}
\label{calc:unif:conv:after:lemma}%
Note that
\begin{subequations}
\begin{align}
  \nonumber
  \bigg|
    1
    +
    \frac{\theta_i^2 D_{\bmw}'(\rho_{i,\bmw})}{2 \varphi_{1,\bmw}(\rho_{i,\bmw})}
    |\bmu_i^\CT \bhtu_{i,\bmw}|^2
  \bigg|
  &
  =
  \begin{alignedat}[t]{1}
  \bigg|
    1
    &
    -
    \theta_i^2
    \bigg[
      + \frac{\varphi_{1,\bmw}(\rho_{i,\bmw})}{2\rho_{i,\bmw}}
      - \frac{\varphi_{1,\bmw}'(\rho_{i,\bmw})}{2}
    \bigg]
    |\bmu_i^\CT \bhtu_{i,\bmw}|^2
    \frac{\varphi_{2,\bmw}(\rho_{i,\bmw})}{\varphi_{1,\bmw}(\rho_{i,\bmw})}
  \\&
    -
    \theta_i^2
    \bigg[
      - \frac{\varphi_{2,\bmw}(\rho_{i,\bmw})}{2\rho_{i,\bmw}}
      - \frac{\varphi_{2,\bmw}'(\rho_{i,\bmw})}{2}
    \bigg]
    |\bmu_i^\CT \bhtu_{i,\bmw}|^2
  \bigg|
  \end{alignedat}
  \\
  \label{eq:unif:conv:after:1}
  &
  \leq
  \big| 1 - \chi_1(\bmw) - \chi_2(\bmw) - 2 \Re \chi_3(\bmw) \big|
  \\&
  \label{eq:unif:conv:after:2}
  +
  \bigg|
    \chi_1(\bmw)
    -
    \theta_i^2
    \bigg[
      + \frac{\varphi_{1,\bmw} (\rho_{i,\bmw})}{2\rho_{i,\bmw}}
      - \frac{\varphi_{1,\bmw}'(\rho_{i,\bmw})}{2}
    \bigg]
    |\bmu_i^\CT \bhtu_{i,\bmw}|^2
    \frac{\varphi_{2,\bmw}(\rho_{i,\bmw})}{\varphi_{1,\bmw}(\rho_{i,\bmw})}
  \bigg|
  \\&
  \label{eq:unif:conv:after:3}
  +
  \bigg|
    \chi_2(\bmw)
    -
    \theta_i^2
    \bigg[
      - \frac{\varphi_{2,\bmw} (\rho_{i,\bmw})}{2\rho_{i,\bmw}}
      - \frac{\varphi_{2,\bmw}'(\rho_{i,\bmw})}{2}
    \bigg]
    |\bmu_i^\CT \bhtu_{i,\bmw}|^2
  \bigg|
  \\&
  \label{eq:unif:conv:after:4}
  +
  \big| 2 \Re \chi_3(\bmw) \big|
  ,
\end{align}
\end{subequations}
where almost surely
\begin{itemize}
  \item \cref{eq:unif:conv:after:1} is equal to zero for all $\bmw \in \clW_>(\nu)$ eventually
  by \cref{p:eq:singvec:sum},
  \item \cref{eq:unif:conv:after:2} converges to zero uniformly over $\bmw \in \clW_>(\nu)$
  by \cref{p:eq:singvec:sum:parts:lim},
  \item \cref{eq:unif:conv:after:3} converges to zero uniformly over $\bmw \in \clW_>(\nu)$
  by \cref{p:eq:singvec:sum:parts:lim:chi2},
  \item \cref{eq:unif:conv:after:4} converges to zero uniformly over $\bmw \in \clW_>(\nu)$
  by \cref{p:eq:singvec:sum:parts:lim:chi3},
\end{itemize}
and we used the fact that
$D_{\bmw}'(\zeta)
= \varphi_{1,\bmw}'(\zeta) \varphi_{2,\bmw}(\zeta)
+ \varphi_{1,\bmw}(\zeta) \varphi_{2,\bmw}'(\zeta)$.
Thus,
\begin{equation*}
  \frac{\theta_i^2 D_{\bmw}'(\rho_{i,\bmw})}{2 \varphi_{1,\bmw}(\rho_{i,\bmw})}
  |\bmu_i^\CT \bhtu_{i,\bmw}|^2
  \quad
  \asunifto{\bmw \in \clW_>(\nu)}
  \quad
  -1
  ,
\end{equation*}
and \cref{eq:singvec:lim:above}
follows
since $2 \varphi_{1,\bmw}(\rho_{i,\bmw}) / (\theta_i^2 D_{\bmw}'(\rho_{i,\bmw}))$
is bounded over $\bmw \in \clW_>(\nu)$.

\subsection{Detailed derivation of \cref{bdd_spectal}}
\label{calc:below:phase:det:bound}%
Consider
\begin{equation}
  \bmR(\bmw,\zeta)
  \coloneqq
  \Bigg(
    \zeta \bmI_{d+n}
    -
    \begin{bmatrix} & \btlY_{\bmw} \\ \btlY_{\bmw}^\CT & \end{bmatrix}
  \Bigg)^{-1}
  =
  \begin{bmatrix}
    \zeta \bmR_1(\bmw,\zeta) & \bmR_1(\bmw,\zeta) \btlY_{\bmw} \\
    \btlY_{\bmw} ^\CT \bmR_1(\bmw,\zeta) & \zeta \bmR_2(\bmw,\zeta)
  \end{bmatrix}
  ,
\end{equation}
where
$\bmR_1(\bmw,\zeta) \coloneqq (\zeta^2 \bmI_d - \btlY_{\bmw}\btlY_{\bmw}^\CT)^{-1}$
and
$\bmR_2(\bmw,\zeta) \coloneqq (\zeta^2 \bmI_n - \btlY_{\bmw}^\CT\btlY_{\bmw})^{-1}$.
With the spectral decomposition of $\btlY_{\bmw} \btlY_{\bmw}^\CT$,
it is easy to see that for $i = 1,\dots,k$,
\begin{equation}
  \label{bdd_spectal_first}
  | \bmu_i^\CT \bhtu_{i,\bmw} |^2 
  \leq
  - \nu \cdot \Im\big\{ \bmu_i^\CT \bmR_1(\bmw,\zeta_{i,\bmw}) \bmu_i \big\}
  =
  - \nu \cdot \Im\bigg\{
    \zeta_{i,\bmw}^{-1}
    \begin{bmatrix} \bmu_i \\ \bm0_n \end{bmatrix}^\CT
    \bmR(\bmw,\zeta_{i,\bmw})
    \begin{bmatrix} \bmu_i \\ \bm0_n \end{bmatrix}
  \bigg\}
  .
\end{equation}
Noting that
\begin{equation*}
  \bmR(\bmw,\zeta)
  =
  \Bigg(
    [\bmG(\bmw,\zeta)]^{-1}
    -
    \begin{bmatrix} \bmU & \\ & \btlQ_{\bmw} \end{bmatrix}
    \begin{bmatrix} & \bmTheta \\ \bmTheta & \end{bmatrix}
    \begin{bmatrix} \bmU & \\ & \btlQ_{\bmw} \end{bmatrix}^\CT
  \Bigg)^{-1}
  ,
\end{equation*}
and applying the Sherman-Morrison-Woodbury formula \cite[Equation (0.7.4.1)]{horn2013ma},
we get that 
\begin{equation}
  \label{woodR1}
  \begin{bmatrix} \bmu_i \\ \bm0_n \end{bmatrix}^\CT
  \bmR(\bmw,\zeta_{i,\bmw})
  \begin{bmatrix} \bmu_i \\ \bm0_n \end{bmatrix}
  =
  \big[
    \btlM(\bmw,\zeta_{i,\bmw})
    -
    \btlM(\bmw,\zeta_{i,\bmw}) [\bmM(\bmw,\zeta_{i,\bmw})]^{-1} \btlM(\bmw,\zeta_{i,\bmw})
  \big]_{ii}
  .
\end{equation}
Combining \cref{bdd_spectal_first,woodR1} yields \cref{bdd_spectal}.

\subsection{Detailed derivation of \cref{eq:singvec:below:bound}}
\label{calc:below:phase:bounds}%
We derive the three bounds of \cref{eq:singvec:below:bound}
one by one.

\begin{enumerate}
\item
The first bound is straightforward.
Note that for all $\bmw \in \Delta_L$,
$b_{\bmw} \geq \sqrt{\min_{\ell} p_\ell v_\ell}$
since
$\|\btlE_{\bmw}\|_{\opnorm} \asto b_{\bmw}$
and
\begin{equation*}
  \|\btlE_{\bmw}\|_{\opnorm}^2
  \geq \frac{1}{d} \|\btlE_{\bmw}\|_{\frob}^2
  = \sum_{\ell=1}^L \frac{w_\ell}{dn} \|\bmE_{\ell}\|_{\frob}^2
  \asto \sum_{\ell=1}^L p_\ell w_\ell v_\ell
  \geq
  C_1
  ,
\end{equation*}
where $C_1 \coloneqq \min_{\ell} p_\ell v_\ell > 0$,
and the convergence follows from the law of large numbers.
Thus,
it follows from \cref{p:eq:singval:limit} that
almost surely,
eventually,
for all $\bmw \in \Delta_L \supseteq \clW_{\leq}(\nu)$,
\begin{equation}
  \label{eq:singvec:below:bound:zeta}
  |\zeta_{i,\bmw}^{-1}|
  \leq |\httheta_{i,\bmw}^{-1}|
  < 2 \cdot \brtheta_{i,\bmw}^{-1}
  \leq 2 \cdot b_{\bmw}^{-1}
  \leq \tlC_1
  ,
\end{equation}
where $\tlC_1 \coloneqq 2 / C_1$
does not depend on $\nu$ or $\bmw$.

\item
For the second bound, observe that
for $\delta=1/2$
and $\tau$ sufficiently large,
it follows from \cref{p:eq:singval:limit} that
almost surely,
eventually,
for all $\bmw \in \clW_{\leq}(\nu)$,
$\zeta_{i,\bmw} \in \Sedge(\bmw,\tau,\delta;\kappa)$.
Thus, using \cref{eq:strong:regularity} yields
that almost surely, eventually,
\begin{equation*}
  \sup_{\bmw \in \Delta_L} \| \bmPi(\bmw,\zeta_{i,\bmw}) \|_{\opnorm} \leq C_2
\end{equation*}
where $C_2$ does not depend on $\nu$ or $\bmw$.
Next using \cref{aniso_in},
yields that almost surely, eventually,
\begin{align}
  \label{bdd_spec1}
  \|\btlM(\bmw,\zeta_{i,\bmw})\|_{\opnorm}
  &
  \leq
  \| \bmPi(\bmw,\zeta_{i,\bmw}) \|_{\opnorm}
  \bigg\| \begin{bmatrix} \bmU & \\ & \btlQ_{\bmw} \end{bmatrix}\bigg\|_{\opnorm}^2
  \\&\quad
  \nonumber
  +
  \bigg\|
    \begin{bmatrix} \bmU & \\ & \btlQ_{\bmw} \end{bmatrix}^\CT
    [\bmG(\bmw,\zeta_{i,\bmw}) - \bmPi(\bmw,\zeta_{i,\bmw})]
    \begin{bmatrix} \bmU & \\ & \btlQ_{\bmw} \end{bmatrix}
  \bigg\|_{\opnorm}
  \\&
  \nonumber
  \leq
  \tlC_2
  ,
\end{align}
where $\tlC_2$ does not depend on $\nu$ or $\bmw$.

\item
For the third bound,
observe that
\begin{equation}
  \label{eq:Mbar:inv:bound}
  \|\bbrM(\bmw,\zeta_{i,\bmw})^{-1}\|_{\opnorm}
  \leq C_3 (\Im \mplaw(\bmw,\zeta_{i,\bmw}))^{-1}
\end{equation}
where $C_3$ does not depend on $\nu$ or $\bmw$.
It has been shown in \cite[Lemma 3.6]{ding2018ana}
that
\begin{equation}
  \label{eq:Mbar:inv:bound:part2}
  \Im \mplaw(\bmw,\zeta_{i,\bmw})
  \geq
  \frac{C_4 \nu}{\sqrt{|\httheta_{i,\bmw} - b_{\bmw}| + \nu}}
\end{equation}
where $C_4$ does not depend on $\nu$ or $\bmw$.
Moreover, note that
the uniform convergence \cref{wtMas}
extends to the domain
$\Omega_{\leq}(\tau)
\coloneqq \{(\bmw,\zeta) \in \Delta_L \times \bbC : \zeta \in \Sedge(\bmw,\tau,\delta;\kappa)\}$.
Combining this with
\cref{eq:Mbar:inv:bound,eq:Mbar:inv:bound:part2}
yields that
almost surely,
eventually,
\begin{align}
  \label{bdd_spec2}
  \|\bmM(\bmw,\zeta_{i,\bmw})^{-1}\|_{\opnorm}
  &
  \leq
  2 \cdot \|\bbrM(\bmw,\zeta_{i,\bmw})^{-1}\|_{\opnorm}
  \leq
  2 \cdot \frac{C_3}{C_4 \nu} \sqrt{|\httheta_{i,\bmw} - b_{\bmw}| + \nu}
  \\&
  \nonumber
  \leq
  2 \cdot \frac{C_3}{C_4 \nu} \sqrt{2\nu + \nu}
  \leq
  \tlC_3 \nu^{-1/2}
  ,
\end{align}
where $\tlC_3 \coloneqq 4 \cdot C_3/C_4$ does not depend on $\nu$ or $\bmw$.
\end{enumerate}

\subsection{Detailed derivation of \cref{p:eq:Dpsi}}
\label{calc:alg:desc:rewrite}%
Analogous to \cite[Section 5.3]{hong2018apo},
observe that
$\psi_{\bmw}$ has the following properties
for all $\bmw \in \Delta_L$:
\begin{propenum}
  \item \label{alg:desc:prop:a}
  $0 = Q_{\bmw}(\psi_{\bmw}(\zeta),\zeta)$ for all $\zeta > b_{\bmw}$ where
  \begin{equation}
    \label{p:eq:Qdef}
    Q_{\bmw}(s,\zeta)
    \coloneqq \frac{c\zeta^2}{s^2}
    + \frac{c-1}{s}
    - \sum_{\ell=1}^L \frac{c_\ell}{s - w_\ell v_\ell}
    ,
  \end{equation}

  \item \label{alg:desc:prop:b}
  $\max_\ell(w_\ell v_\ell) < \psi_{\bmw}(\zeta) < c \zeta^2$,

  \item \label{alg:desc:prop:c}
  $0 < \psi_{\bmw}(b_{\bmw}^+) < \infty$ and $\psi_{\bmw}'(b_{\bmw}^+) = \infty$.
\end{propenum}
Expressing $D_{\bmw}$ in terms of $\psi_{\bmw}$ yields
\begin{align}
  \label{p:eq:Dpsi:calc}
  D_{\bmw}(\zeta)
  &
  = \varphi_{1,\bmw}(\zeta)
    \sum_{\ell=1}^L
      \frac{p_\ell w_\ell}{\zeta - w_\ell v_\ell \varphi_{1,\bmw}(\zeta)/c}
  =
    \sum_{\ell=1}^L
      \frac{c_\ell w_\ell}{\psi_{\bmw}(\zeta) - w_\ell v_\ell}
  = \frac{1-B_{i,\bmw}(\psi_{\bmw}(\zeta))}{\theta_i^2}
  , \\
  \label{p:eq:Dppsi:calc}
  \frac{D_{\bmw}'(\zeta)}{\zeta}
  &
  =
  -\frac{c \psi_{\bmw}'(\zeta)}{\zeta}
  \sum_{\ell=1}^L
  \frac{p_\ell w_\ell}{(\psi_{\bmw}(\zeta) - w_\ell v_\ell)^2}
  =
  -\frac{2c}{\theta_i^2}
  \frac{B_{i,\bmw}'(\psi_{\bmw}(\zeta))}{A_{\bmw}(\psi_{\bmw}(\zeta))}
  .
\end{align}
The second equality in \cref{p:eq:Dppsi:calc}
follows analogously to \cite[Section 5.4]{hong2018apo}
by deriving the identity
\begin{equation}
  \label{p:eq:psip}
  \psi_{\bmw}'(\zeta) = \frac{2c\zeta}{A_{\bmw}(\psi_{\bmw}(\zeta))}
  ,
\end{equation}
from \cref{alg:desc:prop:a} then simplifying.

Rearranging \cref{p:eq:psip}
then applying \cref{alg:desc:prop:c}
yields
\begin{equation}
  A_{\bmw}(\psi_{\bmw}(b_{\bmw}^+))
  = \frac{2cb_{\bmw}}{\psi_{\bmw}'(b_{\bmw}^+)}
  = 0
  ,
\end{equation}
so $\psi_{\bmw}(b_{\bmw}^+)$ is a root of $A_{\bmw}$.
If $\theta_i > \tltheta_{\bmw}$,
then $\rho_{i,\bmw} = D_{\bmw}^{-1}(1/\theta_i^2)$
and rearranging \cref{p:eq:Dpsi:calc} yields
\begin{equation}
  B_{i,\bmw}(\psi_{\bmw}(\rho_{i,\bmw}))
  = 1 - \theta_i^2 D_{\bmw}(\rho_{i,\bmw})
  = 0
  ,
\end{equation}
so $\psi_{\bmw}(\rho_{i,\bmw})$ is a root of $B_{i,\bmw}$.
Recall that
$\psi_{\bmw}(b_{\bmw}^+),\psi_{\bmw}(\rho_{i,\bmw}) \geq \max_\ell(w_\ell v_\ell)$
by \cref{alg:desc:prop:b},
and observe that both $A_{\bmw}(x)$ and $B_{i,\bmw}(x)$
monotonically increase
for $x > \max_\ell(w_\ell v_\ell)$
from negative infinity to one.
Thus, each has exactly
one real root larger than $\max_\ell(w_\ell v_\ell)$,
i.e., its largest real root,
and so $\psi_{\bmw}(b_{\bmw}^+) = \alpha_{\bmw}$
and $\psi_{\bmw}(\rho_{i,\bmw}) = \beta_{i,\bmw}$
when $\theta_i^2 > \tltheta_{\bmw}^2$,
where $\alpha_{\bmw}$ and $\beta_{i,\bmw}$ are
the largest real roots of $A_{\bmw}$ and $B_{i,\bmw}$, respectively.

%% file: content/supp,ext,signal,var.tex
%!TEX root = ../wpca,supp.tex

\section{Detailed verification of consistency for estimators in \cref{ex:var:est}}
\label{sec:supp:ext:signal:var}%
The simple noise variance estimator $\bhtv(\bmY)$
is well-known;
its consistency
follows straightforwardly by noting that
\begin{align*}
  \frac{\|\bmF \bmZ_\ell\|_{\frob}^2}{d n_\ell}
  &
  \leq \|\bmF\|_{\opnorm}^2 \frac{\|\bmZ_\ell\|_{\frob}^2}{d n_\ell}
  = \lambda_1 \frac{\|\bmZ_\ell\|_{\frob}^2}{k n_\ell} \frac{k}{d}
  \asto 0
  , &
  \frac{\|\bmE_\ell\|_{\frob}^2}{d n_\ell}
  &\asto v_\ell
  ,
\end{align*}
by the law of large numbers,
so
$\htv_\ell(\bmY)
= \|\bmY_\ell\|_{\frob}^2/(d n_\ell)
= \|\bmF \bmZ_\ell + \bmE_\ell\|_{\frob}^2/(d n_\ell)
\asto v_\ell$.

The signal variance estimator $\htlambda_i(\bmY)$
is essentially
a bias-corrected version of
the inverse noise variance weighted data eigenvalue $\htlambda^{\inv}_i(\bmY; \bmv)$.
As is well known,
such data eigenvalues are typically upwardly biased
for high-dimensional data;
see, e.g., \cite{hong2018apo,johnstone2018pih,leeb2021oss,nadakuditi2014oaa,passemier2015oeo}.
The shrinkage $\Xi$ simply corrects for the asymptotic bias.

In particular,
note that
it follows from \cref{p:eq:singval:limit}
in the proof of \cref{lem:rec:lim}
that
\begin{align*}
  \htlambda^{\inv}_i(\bmY; \bhtv(\bmY))
  &
  =
  \httheta_{i,\bhtw^{\inv}(\bmY)}^2
  \quad
  \asto
  \quad
  \begin{cases}
    \rho_{i,\bbrw^{\inv}}^2 & \text{if } \theta_i > \tltheta_{\bbrw^{\inv}} , \\
    b_{\bbrw^{\inv}}^2 & \text{otherwise} ,
  \end{cases}
\end{align*}
where
$\htw^{\inv}_\ell(\bmY)
\coloneqq
[\sum_{\ell=1}^L p_\ell / \htv_\ell(\bmY)]^{-1} / \htv_\ell(\bmY)
\asto
\brw^{\inv}_\ell
\coloneqq
\brv/v_\ell$
since $\bhtv(\bmY) \asto \bmv$.
To obtain the algebraic form,
recall
from \cref{proof:rec:lim:unif:alg}
that
\begin{align*}
  \psi_{\bbrw^{\inv}}(b_{\bbrw^{\inv}}^+) &= \alpha_{\bbrw^{\inv}}
  , &
  \psi_{\bbrw^{\inv}}(\rho_{i,\bbrw^{\inv}}) &= \beta_{i,\bbrw^{\inv}}
  , &
  \theta_i > \tltheta_{\bbrw^{\inv}} &\Leftrightarrow \alpha_{\bbrw^{\inv}} < \beta_{i,\bbrw^{\inv}}
  ,
\end{align*}
note that $\alpha_{\bbrw^{\inv}} = \brv + \brv \sqrt{c}$
and $\beta_{i,\bbrw^{\inv}} = \brv + c\lambda_i$,
and apply
\cref{alg:desc:prop:a} from \insupp{calc:alg:desc:rewrite}
to
invert $\psi_{\bbrw^{\inv}}$
and
conclude that
when $c ( \lambda_i / \brv )^2 > 1$
\begin{equation*}
  \htlambda^{\inv}_i(\bmY; \bhtv(\bmY))
  \quad
  \asto
  \quad
  \rho_{i,\bbrw^{\inv}}^2
  =
  \frac{(\lambda_i + \brv/c)(\lambda_i + \brv)}{\lambda_i}
  .
\end{equation*}
This gives the asymptotic bias for inverse noise variance weighted data eigenvalues;
a closely related result was derived in \cite{leeb2021oss}.
It
follows immediately that
\begin{equation*}
  \htlambda_i(\bmY)
  =
  \Xi\Bigg(
    \htlambda^{\inv}_i\big(\bmY; \bhtv(\bmY)\big);
    \Bigg[
      \sum_{\ell=1}^L \frac{p_\ell}{\htv_\ell(\bmY)}
    \Bigg]^{-1}
  \Bigg)
  \quad
  \asto
  \quad
  \Xi\Bigg(
    \frac{( \lambda_i + \brv/c )( \lambda_i + \brv )}{\lambda_i}
    ;
    \brv
  \Bigg)
  = \lambda_i
  ,
\end{equation*}
and so $\htlambda_i(\bmY)$ is consistent.

%% file: content/experiment,preprocess.tex
%!TEX root = ../wpca,supp.tex

\section{Preprocessing details for SDSS data}
\label{exp:preprocess}%
This section describes the details
of the subset selected and the preprocessing performed
on the SDSS data for illustrating optimally weighted PCA
in \cref{sec:exp:astro}.
In particular,
the dataset was formed via the following steps:
\begin{enumerate}
  \item Collect the spectra from DR16Q for which
  \begin{itemize}
    \item $\verb|SURVEY| = \verb|"eboss"|$,
    \item $\verb|PLATEQUALITY| = \verb|"good"|$,
    \item redshift: $2.0 < \verb|Z| < 2.1$,
    \item $\verb|BAL_PROB| < 0.2$,
    \item the measured rest frame wavelengths cover the range of 1480--1620
      without any missing entries
      (i.e., no entries in the range with $\verb|IVAR| = 0$).
  \end{itemize}
  \item Form data $\bmy_j \in \bbR^d$ and variance profile $\bmp_j \in \bbR^d$ vectors
    for each collected spectrum
    via linear interpolation of \verb|FLUX| and $1 \oslash \verb|IVAR|$
    on the grid of rest frame wavelengths
    $\verb|LAMREST| = (1480,1480.5,\dots,1620) \in \bbR^d$.
  \item Center each spectrum $\bmy_j \gets \bmy_j - (1/d)\bm1_{d \times d} \bmy_j$.
  \item Normalize each spectrum so that its mean flux
    for rest frame wavelengths in 1525--1575 is $\pm 1$.
    Namely,
  \begin{enumerate}
    \item $\sigma_j \gets |\verb|mean|(\bmy_j(1525 < \verb|LAMREST| < 1575))|$,
    \item $\bmy_j \gets \bmy_j / \sigma_j$,
    \item $\bmp_j \gets \bmp_j / \sigma_j^2$.
  \end{enumerate}
  \item Remove spectra with $\min(\bmp_j)/\max(\bmp_j) \leq 0.4$
    (i.e., keep only roughly homogeneous variance profiles).
  \item Compute average variances $v_j \gets \verb|mean|(\bmp_j)$ for each sample.
\end{enumerate}